\documentclass[a4paper,12pt,oneside]{article}
\usepackage[english]{babel}
\usepackage[T1]{fontenc} 
\usepackage[utf8]{inputenc}
\usepackage{amsthm}
\usepackage{bbm}
\usepackage{amsmath}
\usepackage{amssymb}  
\usepackage{indentfirst}
\usepackage{fancyhdr}
\usepackage{amsthm}
\usepackage{graphicx}
\usepackage{pdfpages}
\usepackage{esint}
\usepackage{url}

\definecolor{dg}{rgb}{0.01, 0.75, 0.24}

\numberwithin{equation}{section}

\theoremstyle{plain} 
\newtheorem{thm}{Theorem}[section]

\newtheorem{lem}[thm]{Lemma} 
\newtheorem{prop}[thm]{Proposition} 
\newtheorem{rmk}[thm]{Remark}

\def\ibr2 {\int_{B_{\frac{R}{2}}}}
\def\ibro2 {\int_{B_{\frac{\rho}{2}}}}

\def\R{\mathbb{R}}

\def\fibr2 {\displaystyle\fint_{B_{\frac{R}{2}}}}
\def\fibro2 {\fint_{B_{\frac{\rho}{2}}}}

\usepackage{geometry}
\allowdisplaybreaks[4]

\begin{document}

\author {   
\sc{Antonio Giuseppe Grimaldi}\thanks{Dipartimento di Matematica e Applicazioni "R. Caccioppoli", Università degli Studi di Napoli "Federico II", Via Cintia, 80126 Napoli
 (Italy). E-mail: \textit{antoniogiuseppe.grimaldi@unina.it}} 
  }

\title{Higher regularity for minimizers of very degenerate convex integrals}

\maketitle
\maketitle

\begin{abstract}
\noindent In this paper, we consider minimizers of integral functionals of the type 
  \begin{equation*}
    \mathcal{F}(u):= \int_\Omega \dfrac{1}{p} \bigl( |Du(x)|_{\gamma(x)}-1\bigr)_+^p \ \mathrm{d}x,
\end{equation*}
for $p >1$, where $u : \Omega \subset \mathbb{R}^n \to \mathbb{R}^N$, with $N \ge 1$, is a possibly vector-valued function. Here, $| \cdot |_\gamma$ is the associated norm of a bounded, symmetric and
coercive bilinear form on $\mathbb{R}^{nN}$.
We prove that $\mathcal{K}(x,Du)$ is continuous in $\Omega$, for any continuous function $\mathcal{K}: \Omega \times \mathbb{R}^{nN} \rightarrow \mathbb{R}$ vanishing on $\bigl\{ (x,\xi ) \in \Omega \times \mathbb{R}^{nN} : |\xi|_{\gamma(x)} \le 1  \bigr\}$.
\end{abstract}

\medskip
\noindent \textbf{Keywords}: Asymptotically convex functionals, continuity of the gradient, very degenerate elliptic systems.
\medskip \\
\medskip
\noindent \textbf{MSC 2020:} 49N60 $\cdot$ 35D10 $\cdot$ 35J70

\section{Introduction}
In this paper we consider minimizers $u : \Omega  \rightarrow \mathbb{R}^N$, $N \ge 1$, of integral functionals of the type
\begin{equation}\label{FEL}
    \mathcal{F}(u):= \int_\Omega \dfrac{1}{p} \bigl( |Du(x)|_{\gamma(x)}-1\bigr)_+^p \ \mathrm{d}x
\end{equation}
for $p>1$, where $\Omega \subset \mathbb{R}^n$ is a bounded domain, $n \ge 2$. Here, we denote for $x \in \Omega$
\begin{equation*}
    \langle \xi, \eta \rangle_{\gamma(x)} := \sum_{i=1}^{N} \sum_{\alpha,\beta=1}^n \gamma_{\alpha\beta}(x) \xi^i_\alpha \eta^i_\beta, \quad \forall \xi, \eta \in \mathbb{R}^{nN}
\end{equation*}
and $|\xi|_{\gamma(x)}^2:= \langle \xi, \xi \rangle_{\gamma(x)}$.
The $\mathcal{C}^1$ coefficients $\gamma_{\alpha\beta} : \Omega \rightarrow \mathbb{R}$ satisfy the following assumptions:
\begin{itemize}
    \item[(i)] $\gamma_{\alpha\beta}=\gamma_{\beta\alpha}$ (\textit{simmetry});
    
    \item[(ii)] $\exists \ C_0,C_1>0$ such that $$ C_0 |\xi|^2 \le \sum_{\alpha,\beta =1}^n \gamma_{\alpha\beta}(x) \xi_\alpha \xi_\beta \le C_1 |\xi|^2, \quad \forall \xi \in \mathbb{R}^n, \ \forall x \in \Omega;$$

    \item[(iii)] $\exists \ C_2 >0$ such that $$ \sup_{x \in \Omega} | \nabla \gamma_{\alpha\beta}(x)| \le C_2, \quad \forall \alpha,\beta \in \{ 1,...,n \}.$$
\end{itemize}
The Euler-Lagrange system of the functional $\mathcal{F}(u)$ is
\begin{equation}\label{EL}
    \sum_{\alpha,\beta=1}^n D_\alpha \bigl[ h\bigl(|Du|_\gamma\bigr)\gamma_{\alpha\beta} D_\beta u^i \bigr]=0, \quad \forall i =1,...,N
\end{equation}
where $h(t):=\frac{(t-1)_+^{p-1}}{t}$, for $t>0$. A function $u \in W^{1,p}_{loc}(\Omega, \mathbb{R}^N)$ is  a weak solution of the Euler-Lagrange system \eqref{EL} if, and only if,
$$\int_\Omega \langle \textbf{A}(x,Du), D \varphi \rangle_\gamma \ \mathrm{d}x=0$$
holds true for any testing function $\varphi \in \mathcal{C}_0^\infty(\Omega, \mathbb{R}^N)$, where we used the notation
$$\textbf{A}(x,\xi):= h \bigl( |\xi|_{\gamma(x)} \bigr)\xi, \quad \forall \xi \in \mathbb{R}^{nN}, \forall x \in \Omega.$$

In the special case $N=1$ and $\gamma_{\alpha \beta}=\delta^{\alpha \beta}$, the equation \eqref{EL} naturally arises as a model for optimal transport problems with congestion effects. In fact, minimizing \eqref{FEL} is equivalent to the dual minimization problem
\begin{equation}\label{dual}
    \min \biggl\{  \int_\Omega \mathcal{H}(\sigma) \ \mathrm{d} x : \sigma \in L^q(\Omega, \mathbb{R}^n), \ \text{div}\sigma =0, \ \sigma \cdot \nu_{\partial \Omega}=0  \biggr\},
\end{equation}
where the integrand 
$$\mathcal{H}(\sigma):=H(|\sigma|), \quad  \text{with} \ H(t)=t+\dfrac{t^q}{q} \ \text{and} \ \dfrac{1}{p}+\dfrac{1}{q}=1$$
is the convex conjugate of $F(t):= \frac{1}{p}(t-1)_+^p$, or equivalently $F=\mathcal{H}^*$, and $\sigma$ represents the traffic flow. The function $g(t)=H'(t)$ models the congestion effect. Note that $t \mapsto g(t)$ is increasing and $g(0)=1>0$, so that moving in an empty street has non zero cost. As shown in \cite{brasco} the unique minimizer $\sigma(x)$ of \eqref{dual} is given by $D_\xi F(Du(x))$. We refer to \cite{brasco1,brasco2,brasco,carlier,carlier1,santambrogio} and reference therein for a detailed derivation of the model and discussion about the physical interpretation of minimizers.

We note that our functional $\mathcal{F}$ falls into the category of \textit{asymptotically convex} functionals, i.e.\ integral functionals whose integrand is uniformly convex only at infinity. This
class of functionals was first studied by Chipot and Evans in \cite{chipot} and received a lot of attention since then.
Moreover, we observe that no more than Lipschitz regularity can be expected for solutions of equations or systems as in \eqref{EL}. Indeed, every $C_1^{-1/2}$-Lipschitz continuous function solves \eqref{EL}.
Here, we are interested in the $\mathcal{C}^1$-regularity of minimizers in the vectorial case $N \ge 1$. This kind of general structure as in \eqref{EL} was first considered by Tolksdorf \cite{tolk}, where $\mathcal{C}^{1,\alpha}$-regularity was achieved. This result generalized the well-known ones of Uraltseva \cite{ur} for scalar equation and of Uhlenbeck \cite{uh} for systems. However, in our setting we cannot expect to obtain $\mathcal{C}^{1,\alpha}$-regularity, since the solution $u(x)$ is in general no more that Lipschitz continuous. Instead, we will show that $\mathcal{K}(x,Du(x))$ is continuous for any continuous function $\mathcal{K}: \Omega \times \mathbb{R}^{nN} \rightarrow \mathbb{R}$ vanishing on the set $\bigl\{ (x,\xi ) \in \Omega \times \mathbb{R}^{nN} : |\xi|_{\gamma(x)} \le 1  \bigr\}$. For integrands without $x$-dependence such a regularity results for $Du$ was first proved in \cite{santambrogio} in the special case $N=1$, $n=2$ and extended to the case $N=1$, $n \ge 2$ in \cite{clop,colombo}. In the recent paper \cite{BoDuGiPdN}, this result was extended to the vectorial case $N \ge 2$. Regularity of minimizers of functionals with $x$-dependence of the form 
\begin{equation}\label{fxd}
    \tilde{\mathcal{F}}(u) := \int_\Omega \dfrac{a(x)}{p} ( |Du|-1  )^p_+ \ \mathrm{d} x
\end{equation}
was studied in \cite{mons}. In a parabolic setting, systems with structure as in \eqref{EL} was considered in \cite{BoDuLiSc}, where the local H\"older continuity of the spatial gradient has been established.

In this paper, we will go one step further and consider functionals as in \eqref{FEL} with integrands more general than the ones defined in \eqref{fxd}. To simplify the exposition we only consider the homogeneous case, but an analogous result is also expected to hold in the non-homogeneous case under suitable assumptions on the right-hand side.

The main result of this work is the following.
\begin{thm}\label{mainthm}
    Let $p>1$ and $u \in W^{1,p}_{loc}(\Omega, \mathbb{R}^N) \cap  W^{1,\infty}_{loc}(\Omega, \mathbb{R}^N) $ be a weak solution of \eqref{EL} in $\Omega$. Then, $\mathcal{K}(x,Du) \in \mathcal{C}^0(\Omega)$, for any continuous function $\mathcal{K}: \Omega \times \mathbb{R}^{nN} \rightarrow \mathbb{R}$ vanishing on the set $\bigl\{ (x,\xi ) \in \Omega \times \mathbb{R}^{nN} : |\xi|_{\gamma(x)} \le 1  \bigr\}$.
\end{thm}

We give a brief sketch of 
the proof of Theorem \ref{mainthm}, which is divided into several steps. In the first step, we approximate $u$ with a sequence of functions $(u_\varepsilon)_{\varepsilon >0}$, which are weak solutions to the Dirichlet problems
\begin{equation}\label{REL}
     \sum_{\alpha,\beta=1}^n D_\alpha \bigl[ \bigl(h\bigl(|Du_\varepsilon|_\gamma\bigr)+ \varepsilon \bigr)\gamma_{\alpha\beta}D_\beta u_\varepsilon \bigr]=0
\end{equation}
on a ball compactly contained in $\Omega$ and with Dirichlet boundary datum $u$. The approximating functions $u_\varepsilon$ are more regular than $u$. In particular, they possess weak second derivatives, which allow us to differentiate the regularized system \eqref{REL} with respect to $x_\beta$ and obtain some energy estimates for the second derivatives of $u_\varepsilon$. In the next step, we show uniform energy bounds, uniform quantitative interior $L^\infty$-gradient bounds and uniform quantitative higher differentiability estimates for the solutions $u_\varepsilon$ with respect to $\varepsilon$. Moreover, we derive strong $L^p$-convergence of
$\mathcal{G}_\delta(x,Du_\varepsilon) \rightarrow \mathcal{G}_\delta(x,Du) $ as $\varepsilon \to 0$. The map $\mathcal{G}_\delta : \Omega \times \mathbb{R}^{nN} \to \mathbb{R}^{nN}$, with $\delta \in (0,1]$, is defined by
$$ \mathcal{G}_\delta(x,\xi):= \dfrac{(|\xi|_{\gamma(x)}-1-\delta)_+}{|\xi|_{\gamma(x)}}\xi.$$
Observe that $\mathcal{G}_\delta$ vanishes on the larger set $\{ |\xi|_{\gamma(x)} \le 1+ \delta \}$ and the regularized system \eqref{REL} is uniformly elliptic with respect to $\varepsilon$ on the complement of this set. This allows us to prove the H\"older continuity of $\mathcal{G}_\delta (x, Du_\varepsilon)$ uniformly with respect to $\varepsilon$. Letting $\varepsilon \to 0$, we obtain the H\"older continuity of $\mathcal{G}_\delta(x,Du)$ by an application of the Ascoli-Arzelà Theorem. However, the constants in the quantitative estimates, i.e.\ the H\"older exponent and the H\"older norm, may blow up when $\delta \to 0$. Therefore, we cannot expect $\mathcal{G}(x,Du)$ to be H\"older continuous. 
To prove the H\"older continuity of $\mathcal{G}_\delta(x,Du_\varepsilon)$, we distinguish between two
different regimes: the \textit{degenerate} and \textit{non-degenerate regime}. The degenerate regime
is characterized by the fact that the measure of those points in a ball in which $|\mathcal{G}_\delta (x, Du_\varepsilon)|_{\gamma(x)}$
is far from its supremum is large, while the non-degenerate regime is characterized
by the opposite condition. In the non-degenerate regime we compare $u_\varepsilon$ with a solution of a
linearized system.
This allows us to derive a quantitative $L^2$-decay for the excess of $\mathcal{G}_{ \delta}(x,Du_\varepsilon)$
$$\Psi_\delta(x_0,r):=\fint_{B_{r}(x_0)}
    |\mathcal{G}_{ \delta}(x,Du_\varepsilon)- \bigl( \mathcal{G}_{ \delta}(x,Du_\varepsilon) \bigr)_{x_0,r}|^2 \ \mathrm{d}x $$
on some smaller ball (see Proposition \ref{ND}). This step is based on a suitable
comparison estimate and the higher integrability of $u_\varepsilon$. On the smaller ball we are
again in the non-degenerate regime, so that the argument can be iterated yielding a
Campanato-type estimate for the $L^2$-excess of $\mathcal{G}_\delta (x, Du_\varepsilon)$. In the degenerate regime we
establish that 
$$U_\varepsilon := (|Du_\varepsilon|_\gamma -1- \delta )^2_+$$
 is a sub-solution to a linear uniformly elliptic
equation with measurable coefficients. Then, a De Giorgi type argument allows a reduction of
the modulus of $\mathcal{G}_\delta (x, Du_\varepsilon)$ on some smaller ball (see Proposition \ref{D}). However, on this
smaller scale we are either in the degenerate or non-degenerate
regime. In the
non-degenerate regime we can conclude as above, while in the degenerate regime the
reduction of the modulus of $\mathcal{G}_\delta (x, Du_\varepsilon)$ applies again. This argument can be iterated as
long as we stay in the degenerate regime. However, if at a certain scale the switching
from degenerate to non-degenerate occurs, the above Campanato type decay applies.
If no switching occurs, we have at any scale of the iteration process a reduction of
the modulus of $\mathcal{G}_\delta (x, Du_\varepsilon)$. This, however, shows that the supremum of $\mathcal{G}_\delta (x, Du_\varepsilon)$ on shrinking concentric balls converges to $0$. Thus, we can conclude that 
$\mathcal{G}_\delta (x, Du_\varepsilon)$ is H\"older continuous. Using the uniform convergence of $\mathcal{G}_\delta(x,Du) \to \mathcal{G}(x,Du)$, we obtain that $\mathcal{G}(x,Du)$ is continuous. From this, we easily derive the continuity of $\mathcal{K}(x,Du)$.

\section{Notation and preliminary results}
In this paper we shall denote by $C$ or $c$ a general positive constant that may vary on different occasions, even within the same line of estimates. Relevant dependencies will be suitably emphasized using parentheses or subscripts. 
\\In what follows, $B_{r}(x_0)= \{ x \in \mathbb{R}^{n} : |x-x_0 | < r  \}$ will denote the ball centered at $x_0$ of radius $r$. For a function $v \in L^{1}(B_r(x_0),\mathbb{R}^k)$, the symbol
\begin{center}
$(v)_{x_0,r}:=\displaystyle\fint_{B_r(x_0)} v(x) \ \mathrm{d}x = \dfrac{1}{|B_r(x_0)|} \displaystyle\int_{B_r(x_0)} v(x) \ \mathrm{d}x $.
\end{center}
will denote the integral mean of the function $v$ over the ball $B_r(x_0)$.

Let define
\begin{equation}\label{fung}
    g(t):= \dfrac{(t-1)_+^p}{t}, \quad \forall t >0
\end{equation}
and \begin{equation}\label{funG}
    \mathcal{G}(x,\xi):= \dfrac{(|\xi|_{\gamma(x)}-1)_+}{|\xi|_{\gamma(x)}}\xi, \quad \text{for} \ \xi \in \mathbb{R}^k, k \in \mathbb{N}, x \in \Omega.
\end{equation}
Observe that
$$g \bigl( |\xi|_{\gamma(x)} \bigr) \xi = |\mathcal{G}(x,\xi)|^{p-1}_{\gamma(x)} \mathcal{G}(x,\xi), \quad \forall \xi \in \mathbb{R}^k, \forall x \in \Omega.$$
Moreover, we define for $\delta \in (0,1]$
\begin{equation}
    \mathcal{G}_\delta(x,\xi):= \dfrac{(|\xi|_{\gamma(x)}-1-\delta)_+}{|\xi|_{\gamma(x)}}\xi, \quad \forall \xi \in \mathbb{R}^k, k \in \mathbb{N}, \forall x \in \Omega.
\end{equation}
\\Let us recall that for every function $v : \R^n \rightarrow \R^k$ the finite difference operator is defined by 
$$\tau_{s,h}v(x):= v(x+h e_s)-v(x)$$
where $h \in \R^n$, $e_s$ is the unit vector in the $x_s$ direction and $s \in \{1,...,n\}$. 
The next result about finite difference operator is a kind of integral version of Lagrange Theorem (see for instance \cite[Lemma 8.1]{Giusti}).
\begin{lem}\label{ldiff}
If $0<\rho<R,$ $|h|<\frac{R-\rho}{2},$ $1<p<+\infty$, $s \in \{1,...,n\}$ and $v \in W^{1,p}(B_{R}, \R^{k })$, then
\begin{center}
$\displaystyle\int_{B_{\rho}} |\tau_{s,h}v(x)|^{p} \ \mathrm{d} x \leq c(n,p)|h|^{p} \displaystyle\int_{B_{R}} |Dv(x)|^{p} \ \mathrm{d}x$.
\end{center}
\end{lem}
Now, we recall the fundamental Sobolev embedding property (see e.g. \cite[Lemma 8.2]{Giusti}).
\begin{lem}\label{ldiff2}
Let $v : \mathbb{R}^{n} \rightarrow \mathbb{R}^{k}$, $v \in L^{p}(B_{R},\R^k)$, $1<p<+\infty$. Suppose that there exist $\rho \in (0,R)$ and $M>0$ such that
\begin{center}
$\displaystyle\sum_{s=1}^{n} \displaystyle\int_{B_{\rho}} |\tau_{s,h}v(x)|^{p} \ \mathrm{d}x \leq M^{p}|h|^{p}$,
\end{center}
for every $h \in \R^n $ with $|h|< \frac{R - \rho}{2}$. Then $v \in W^{1,p}(B_{\rho},\R^k)$ and it holds
\begin{center}
$\Vert Dv \Vert_{L^{p}(B_{\rho})} \leq M$.
\end{center}
\end{lem}

\subsection{Algebraic Lemmas}
In this section we state some algebraic inequalities that we will use later on.
\begin{lem}\label{lem1}
    For $\eta, \xi \in \mathbb{R}^k_{\neq 0}, k \in \mathbb{N}, x \in \Omega$, we have
    $$\bigg| \dfrac{\eta}{|\eta|_{\gamma(x)}}- \dfrac{\xi}{|\xi|_{\gamma(x)}} \bigg| \ \le c 
    \dfrac{|\eta-\xi|}{|\eta|} 
    ,$$
    for a constant $c:=c(C_0,C_1)$.
\end{lem}
\proof 
Using the property $(ii)$ of $\gamma$ we have
\begin{align*}
    \bigg| \dfrac{\eta}{|\eta|_{\gamma(x)}}- \dfrac{\xi}{|\xi|_{\gamma(x)}} \bigg|
\dfrac{|\eta|}{|\eta-\xi|} = & \ \dfrac{\big|\eta |\xi|_{\gamma(x)}-\xi |\eta|_{\gamma(x)} \big|}{|\eta|_{\gamma(x)} |\xi|_{\gamma(x)}} \dfrac{|\eta|}{|\eta-\xi|}\\
\le & \ c \  \dfrac{|\eta-\xi||\xi|_{\gamma(x)}+|\xi| \big| |\xi|_{\gamma(x)}- |\eta|_{\gamma(x)} \big| }{|\xi|_{\gamma(x)}|\eta-\xi|} \\
\le & \ c \ \biggl[  1+  \dfrac{|\xi|  |\xi -\eta|_{\gamma(x)}  }{|\xi|_{\gamma(x)}|\eta-\xi|} \biggr] =: c(C_0,C_1).
\end{align*}
This implies the desired inequality.
\endproof

The following lemma can be deduced as in \cite[Lemma 8.3]{Giusti}
\begin{lem}\label{lem2}
For any $\alpha >0$, there exists a constant $c:=c(\alpha, C_0,C_1)$ such that $\forall \eta,\xi \in \mathbb{R}^k_{\neq 0}, k \in \mathbb{N}, x \in \Omega$, we have
$$\dfrac{1}{c} \ \big| |\eta|^{\alpha-1}_{\gamma(x)}\eta- |\xi|^{\alpha-1}_{\gamma(x)}\xi \big| \le \bigl(  |\eta|_{\gamma(x)}+|\xi|_{\gamma(x)} \bigr)^{\alpha-1} |\eta-\xi| \le c  \ \big| |\eta|^{\alpha-1}_{\gamma(x)}\eta- |\xi|^{\alpha-1}_{\gamma(x)}\xi \big|.$$
\end{lem}

\begin{lem}\label{lem3}
Let $\delta \ge 0$ and $\eta, \xi \in \mathbb{R}^k, k \in \mathbb{N}, x\in \Omega$. Then, we have
$$\big| \mathcal{G}_\delta(x,\xi)-\mathcal{G}_\delta(x,\eta) \big| \le c |\eta- \xi |.$$
Moreover, if $\delta >0$ and $|\eta|_{\gamma(x)} \ge 1+\delta$ there holds
$$|\eta-\xi| \le c \biggl(1+\dfrac{1}{\delta} \biggr) |\mathcal{G}(x,\eta)-\mathcal{G}(x,\xi) |,$$
    where $c:=c(C_0,C_1)$.
\end{lem}
\proof The proof goes as that of \cite[Lemma 2.3]{BoDuGiPdN}  taking into account Lemma \ref{lem1}.
\endproof

For the proof of next lemmas we refer to \cite[Section 2.2]{BoDuGiPdN}. 
\begin{lem}\label{lem4}
There exists a constant $c:=c(p)$ such that for any $a>1$ and $b \ge 0$, we have
$$|h(b)-h(a)|b \le c(p) \dfrac{[a-1+(b-1)_+]^{p-1}}{a-1}|b-a|.$$
\end{lem}

\begin{lem}\label{lem5}
    For $a>1$, we have
    $$|h'(a)| \le \dfrac{p(a-1)^{p-2}}{a}.$$
    Moreover, for $a,b>1$ there holds
    $$|h'(b)b-h'(a)a|\le c(p) [(a-1)^{p-3}+(b-1)^{p-3}]|b-a|.$$
\end{lem}

\begin{lem}\label{lem6}
    For any $t \ge 0$, we have
    \begin{itemize}
        \item[1)] $g(t)^2 \le h(t) (t-1)^p_+$;

        \item[2)] $g(t)^2+g'(t)^2t^2 \le \dfrac{p^2}{p-1}[h(t)+h'(t)t](t-1)^p_+.$
    \end{itemize}
\end{lem}

\subsection{Bilinear forms}

For $\varepsilon \in [0,1]$, we define
$$h_\varepsilon(t):=h(t)+\varepsilon, \quad \forall t \ge 0.$$
Moreover, we let
$$\textbf{A}_\varepsilon(x,\xi):=h_\varepsilon \bigl(|\xi|_{\gamma(x)} \bigr)\xi, \quad \forall \xi \in \mathbb{R}^{nN}.$$
Now, we introduce bilinear forms on the spaces $\mathbb{R}^{n^2 N}$, $ \mathbb{R}^{n N}$ and $\mathbb{R}^{n } $ that will be useful in the following sections. For $x \in \Omega,\xi \in \mathbb{R}^{nN} \setminus \{0 \}$, with $|\xi|_{\gamma(x)} \neq 1$ if $1 < p <2$, we define the bilinear forms
$$\mathcal{A}_\varepsilon(x,\xi)(\eta,\zeta):= \biggl[ h_\varepsilon \bigl(|\xi|_{\gamma(x)}\bigr)\gamma_{\alpha\beta}(x)\delta^{ij}+h'_\varepsilon \bigl(|\xi|_{\gamma(x)} \bigr)|\xi|_{\gamma(x)} \dfrac{\gamma_{\alpha \delta}(x) \xi^i_\delta\gamma_{\beta\kappa}(x)\xi^j_\kappa}{|\xi|_{\gamma(x)}}   \biggr] \gamma_{\nu\sigma}(x) \eta^i_{\alpha\nu}\zeta^j_{\beta\sigma}, $$
for all $\eta,\zeta \in \mathbb{R}^{n^2N}$, and 
$$\mathcal{B}_\varepsilon(x,\xi)(\eta,\zeta):= \biggl[ h_\varepsilon \bigl(|\xi|_{\gamma(x)}\bigr)\gamma_{\alpha\beta}(x)\delta^{ij}+h'_\varepsilon \bigl(|\xi|_{\gamma(x)} \bigr)|\xi|_{\gamma(x)} \dfrac{\gamma_{\alpha \delta}(x) \xi^i_\delta\gamma_{\beta\kappa}(x)\xi^j_\kappa}{|\xi|^2_{\gamma(x)}}   \biggr]  \eta^i_{\alpha}\zeta^j_{\beta}, $$
for all $\eta,\zeta \in \mathbb{R}^{nN}$, and finally
$$\mathcal{C}_\varepsilon(x,\xi)(\eta,\zeta):= \biggl[ h_\varepsilon \bigl(|\xi|_{\gamma(x)}\bigr)\gamma_{\alpha\beta}(x)+h'_\varepsilon \bigl(|\xi|_{\gamma(x)} \bigr)|\xi|_{\gamma(x)} \dfrac{\gamma_{\alpha \delta}(x) \xi^i_\delta\gamma_{\beta\kappa}(x)\xi^j_\kappa}{|\xi|^2_{\gamma(x)}}   \biggr] \eta_{\alpha}\zeta_{\beta}, $$
for all $\eta,\zeta \in \mathbb{R}^{n}$. The next lemma provides the revelant ellipticity and boundedness properties of the bilinear forms $\mathcal{A}_\varepsilon(x,\xi), \mathcal{B}_\varepsilon(x,\xi)$ and $\mathcal{C}_\varepsilon(x,\xi)$. We will adopt the following abbreviations
\begin{equation*}
\lambda(t):=
    \begin{cases}
        \min \{ h(t),(p-1)(t-1)^{p-2} \} \ & \text{if} \ t>1 \\
        0 \ & \text{if} \ 0 \le t \le 1
    \end{cases}
\end{equation*}
and 
\begin{equation*}
\Lambda(t):=
    \begin{cases}
        \max \{ h(t),(p-1)(t-1)^{p-2} \} \ & \text{if} \ t>1 \\
        0 \ & \text{if} \ 0 \le t \le 1.
    \end{cases}
\end{equation*}

\begin{lem}\label{ellipticity}
    Let $\varepsilon \in [0,1]$ and $\xi \in \mathbb{R}^{nN}$ and $x \in \Omega$. The bilinear form $\mathcal{A}_\varepsilon(x,\xi)$ satisfies
    $$ \bigl[  \varepsilon + \lambda \bigl( |\xi|_{\gamma(x)} \bigr) \bigr]|\eta|^2_{\gamma(x)}\le \mathcal{A}_\varepsilon(x,\xi)(\eta,\eta) \le \bigl[  \varepsilon + \Lambda \bigl( |\xi|_{\gamma(x)} \bigr) \bigr]|\eta|^2_{\gamma(x)},$$
    for any $\eta \in \mathbb{R}^{n^2N}$. The analogous estimates hold for the bilinear form $\mathcal{B}_\varepsilon(x,\xi)$ for any  $\eta \in \mathbb{R}^{nN}$, as well as for $\mathcal{C}_\varepsilon(x,\xi)$ for any $\eta \in \mathbb{R}^{n}$.
\end{lem}
\proof If $|\xi|_{\gamma(x)} \le 1$, the claim follows trivially from the definition of $\mathcal{A}_\varepsilon(x,\xi)$. Therefore, we consider the case $|\xi|_{\gamma(x)} >1$.
We first note that the positive definiteness of $(\gamma_{\alpha\beta})$ and the Cauchy-Schwarz inequality imply
\begin{equation}\label{CSineq}
    0 \le \dfrac{\gamma_{\alpha \delta}(x) \xi^i_\delta\gamma_{\beta\kappa}(x)\xi^j_\kappa}{|\xi|^2_{\gamma(x)}} \gamma_{\nu\sigma}(x)\eta^i_{\alpha\nu}\eta^j_{\beta\sigma} \le |\eta|^2_{\gamma(x)}
\end{equation}
for every $\eta \in \mathbb{R}^{n^2N}$.
In order to prove the inequality we distinguish two cases. If $h_\varepsilon' \bigl( |\xi|_{\gamma(x)}   \bigr) \ge 0$, we can discard the second term in the definition of $\mathcal{A}_\varepsilon(x,\xi)$ due to the lower bound in \eqref{CSineq} and obtain
$$ \mathcal{A}_\varepsilon(x,\xi)(\eta, \eta) \ge  h_\varepsilon \bigl( |\xi|_{\gamma(x)}   \bigr) |\eta|^2_{\gamma(x)} \ge \bigl[  \varepsilon + \lambda \bigl( |\xi|_{\gamma(x)} \bigr) \bigr]|\eta|^2_{\gamma(x)}. $$
Otherwise, if $ h_\varepsilon' \bigl( |\xi|_{\gamma(x)}   \bigr) < 0$, we use the upper bound in \eqref{CSineq} and the equality
$$h(t)+h'(t)t=(p-1)(t-1)^{p-2} \quad \text{for every} \ t>1$$
to conclude
\begin{align*}
  \mathcal{A}_\varepsilon(x,\xi)(\eta, \eta) \ge &  \ \bigl[  \varepsilon + h \bigl( |\xi|_{\gamma(x)}   \bigr) +  h' \bigl( |\xi|_{\gamma(x)}   \bigr) |\xi|_{\gamma(x)}  \bigr]|\eta|^2_{\gamma(x)} \\
  = & \ \bigl[ \varepsilon + (p-1) \bigl( |\xi|_{\gamma(x)} -1 \bigr)^{p-2} \bigr] |\eta|^2_{\gamma(x)} \ge \bigl[  \varepsilon + \lambda \bigl( |\xi|_{\gamma(x)} \bigr) \bigr]|\eta|^2_{\gamma(x)}.
\end{align*}
This proves the asserted lower bound. For the derivation of the upper bound, we proceed analogously. If $h_\varepsilon' \bigl( |\xi|_{\gamma(x)}   \bigr) \ge 0$, we use the upper bound in \eqref{CSineq} to deduce
\begin{align*}
  \mathcal{A}_\varepsilon(x,\xi)(\eta, \eta) \le &  \ \bigl[  \varepsilon + h \bigl( |\xi|_{\gamma(x)}   \bigr) +  h' \bigl( |\xi|_{\gamma(x)}   \bigr) |\xi|_{\gamma(x)}  \bigr]|\eta|^2_{\gamma(x)} \\
  = & \ \bigl[ \varepsilon + (p-1) \bigl( |\xi|_{\gamma(x)} -1 \bigr)^{p-2} \bigr] |\eta|^2_{\gamma(x)} \le \bigl[  \varepsilon + \Lambda \bigl( |\xi|_{\gamma(x)} \bigr) \bigr]|\eta|^2_{\gamma(x)},
\end{align*}
while in the case $h_\varepsilon' \bigl( |\xi|_{\gamma(x)}   \bigr) < 0$, we use the lower bound in \eqref{CSineq} to get
$$ \mathcal{A}_\varepsilon(x,\xi)(\eta, \eta) \le  h_\varepsilon \bigl( |\xi|_{\gamma(x)}   \bigr) |\eta|^2_{\gamma(x)} \le \bigl[  \varepsilon + \Lambda \bigl( |\xi|_{\gamma(x)} \bigr) \bigr]|\eta|^2_{\gamma(x)}. $$
This proves the claim for the bilinear form $\mathcal{A}_\varepsilon(x,\xi)$. The estimates for $\mathcal{B}_\varepsilon(x,\xi)$
and $\mathcal{C}_\varepsilon(x,\xi)$ follow in the same way if we replace \eqref{CSineq} by
\begin{equation*}
    0 \le \dfrac{\gamma_{\alpha \delta}(x) \xi^i_\delta\gamma_{\beta\kappa}(x)\xi^j_\kappa}{|\xi|^2_{\gamma(x)}} \eta^i_{\alpha}\eta^j_{\beta} \le |\eta|^2_{\gamma(x)},
\end{equation*}
for every $\eta \in \mathbb{R}^{nN}$, and by
\begin{equation*}
    0 \le \dfrac{\gamma_{\alpha \delta}(x) \xi^i_\delta\gamma_{\beta\kappa}(x)\xi^j_\kappa}{|\xi|^2_{\gamma(x)}} \eta_{\alpha}\eta_{\beta} \le |\eta|^2_{\gamma(x)},
\end{equation*}
for every $\eta \in \mathbb{R}^{n}$, respectively.
\endproof

In the following lemma, we prove some monotonicity and growth properties of the vector field $\mathbf{A}_\varepsilon$.
\begin{lem}\label{monotonicity}
Let $\varepsilon \in [0,1]$ and $x,y \in \Omega$ and $\xi, \tilde{\xi} \in \mathbb{R}^{nN}$, with $|\xi|_{\gamma(x)}>1$. Then, we have 
\item[1)] $|\textbf{A}_\varepsilon(x,\tilde{\xi})-\textbf{A}_\varepsilon(x,\xi)|_{\gamma(x)} \le c(p,C_0,C_1) \biggl[  
 \varepsilon + \dfrac{\bigl[ |\xi|_{\gamma(x)}-1 + \bigl( |\tilde{\xi}|_{\gamma(x)}-1 \bigr)_+\bigr]^{p-1}}{|\xi|_{\gamma(x)}-1}\biggr] |\xi-\tilde{\xi}|;$
\item[2)] $\langle \textbf{A}_\varepsilon(x, \tilde{\xi})-\textbf{A}_\varepsilon(x,\xi),  \tilde{\xi}-\xi \rangle_{\gamma(x)}  \ge c(C_0,C_1) \biggl[ \varepsilon+ \dfrac{\min \{ 1,p-1\}}{2^{p+1}}  \dfrac{(|\xi|_{\gamma(x)}-1)^p}{|\xi|_{\gamma(x)}\bigl(|\xi|_{\gamma(x)}+|\tilde{\xi}|_{\gamma(x)} 
 \bigr)}\biggr]|\tilde{\xi}-\xi|^2;$
 \item[3)]  $|\langle \textbf{A}_\varepsilon (x,\xi), \tilde{\xi} \rangle_{\gamma(x)}- \langle \textbf{A}_\varepsilon (y,\xi), \tilde{\xi} \rangle_{\gamma(y)}| \\
   \text{ \quad \quad   \quad \quad  \quad \quad } \le c(p,C_0,C_1,C_2) \biggl[\varepsilon+ \dfrac{\bigl[  |\xi|_{\gamma(x)}-1+ \bigl( |\xi|_{\gamma(y)}-1 \bigr)_+\bigr]^{p-1}}{|\xi|_{\gamma(x)}-1}  \biggr]|x-y||\xi||\tilde{\xi}|$
\end{lem}
\proof
The first two inequalities can be proved similarly as in \cite[Lemma 2.8]{BoDuGiPdN}. Therefore, we take care of the remaining inequality. Due to hypotheses $(ii), (iii)$ for $\gamma$ and Lemma \ref{lem4}, we have
\begin{align*}
    &|\langle \textbf{A}_\varepsilon (x,\xi), \tilde{\xi} \rangle_{\gamma(x)}- \langle \textbf{A}_\varepsilon (y,\xi),  \tilde{\xi} \rangle_{\gamma(y)}| \\
    & \ \ \ =  \ |h_\varepsilon(|\xi|_{\gamma(x)})\langle \xi, \tilde{\xi} \rangle_{\gamma(x)}- h_\varepsilon(|\xi|_{\gamma(y)}) \langle \xi, \tilde{\xi} \rangle_{\gamma(y)}| \\
    & \ \ \ \le  \   h_\varepsilon(|\xi|_{\gamma(x)})  |\langle \xi, \tilde{\xi} \rangle_{\gamma(x)}-\langle \xi, \tilde{\xi} \rangle_{\gamma(y)}| +  |h(|\xi|_{\gamma(x)}) -h_\varepsilon(|\xi|_{\gamma(y)}) | | \xi|_{\gamma(y)} |\tilde{\xi} |_{\gamma(y)}\\
    & \ \ \ \le  \ c \biggl[\varepsilon + \dfrac{\bigl( |\xi|_{\gamma(x)}-1 \bigr)^{p-1}}{|\xi|_{\gamma(x)}} \biggr] |x-y| |\xi||\tilde{\xi}|  +   c \dfrac{\bigl[ |\xi|_{\gamma(x)}-1 +\bigl( |\xi|_{\gamma(y)}-1 \bigr)_+\bigr]^{p-1}}{|\xi|_{\gamma(x)}-1}   ||\xi|_{\gamma(x)}-|\xi|_{\gamma(y)}||\tilde{\xi}|\\
    & \ \ \ \le   \ c \biggl[\varepsilon + \dfrac{\bigl( |\xi|_{\gamma(x)}-1 \bigr)^{p-1}}{|\xi|_{\gamma(x)}-1} \biggr] |x-y| |\xi||\tilde{\xi}|  +   c \dfrac{\bigl[ |\xi|_{\gamma(x)}-1 +\bigl( |\xi|_{\gamma(y)}-1 \bigr)_+\bigr]^{p-1}}{|\xi|_{\gamma(x)}-1}  |x-y||\xi| |\tilde{\xi}|\\
     & \ \ \ \le  \ c  \biggl[\varepsilon+ \dfrac{\bigl[  |\xi|_{\gamma(x)}-1+ \bigl( |\xi|_{\gamma(y)}-1 \bigr)_+\bigr]^{p-1}}{|\xi|_{\gamma(x)}-1}  \biggr]|x-y||\xi||\tilde{\xi}|
\end{align*}
for a constant $c:=c(p,C_0,C_1,C_2)$.
This completes the proof of the lemma.
\endproof

\begin{lem}\label{bilformB}
   Let $\varepsilon \in [0,1]$ and $x,y \in \Omega$ and $\xi  \in \mathbb{R}^{nN}$, with $  1+\frac{1}{4} \mu \le |\xi|_{\gamma(x)}, |\xi|_{\gamma(y)} \le 1 +2 \mu$, for some $\mu >0$. Then, we have  
$$  \big| \bigl[  \mathcal{B}_\varepsilon(x,\xi) -  \mathcal{B}_\varepsilon(y,\xi)  \bigr] \bigl( \xi, \eta \bigr)   \big| \le c (\varepsilon+\mu^{p-2}) |x-y||\xi|^2|\eta|  $$
   for every $\eta \in \mathbb{R}^{nN}$ and for a constant $c:=c(p,C_0,C_1,C_2)$.
\end{lem}
\proof Recalling the definition of $\mathcal{B}_\varepsilon$, we compute
\begin{align}\label{i1i2}
     \big| \bigl[  \mathcal{B}_\varepsilon & (x,\xi) -  \mathcal{B}_\varepsilon(y,\xi)  \bigr] \bigl( \xi, \eta \bigr)   \big|  \notag \\
     \le  & \ \big|   h_\varepsilon \bigl(  |\xi|_{\gamma(x) }\bigl)\gamma_{\alpha\beta}(x) -  h_\varepsilon \bigl(  |\xi|_{\gamma(y) } \bigr)\gamma_{\alpha\beta}(y) \big|\xi_\alpha^i\eta^i_\beta \notag \\
     &+ \bigg| h'_\varepsilon \bigl(  |\xi|_{\gamma(x) } \bigr)   |\xi|_{\gamma(x) } \dfrac{\gamma_{\alpha \delta}(x)\xi^i_\delta \gamma_{\beta k} (x) \xi^j_{k}}{|\xi|_{\gamma(x)}^2} -
      h'_\varepsilon \bigl(  |\xi|_{\gamma(y) } \bigr)   |\xi|_{\gamma(y) } \dfrac{\gamma_{\alpha \delta}(y)\xi^i_\delta \gamma_{\beta k} (y) \xi^j_{k}}{|\xi|_{\gamma(y)}^2} \bigg| \xi_\alpha^i\eta^j_\beta \notag\\
      =& \ I_1+I_2.
\end{align}
We start by estimating the term $I_1$. Using properties $(ii)$ and $(iii)$ of $\gamma$ and Lemma \ref{lem4}, we have 
\begin{align*}
    I_1  \le & \ \sum_{\alpha, \beta =1}^n \bigl(h_\varepsilon \bigl(  |\xi|_{\gamma(x) }\bigl) \big|  \gamma_{\alpha \beta}(x)- \gamma_{\alpha \beta}(y)\big| |\xi||\eta|
    + \big|   h_\varepsilon \bigl(  |\xi|_{\gamma(x) }\bigl) -  h_\varepsilon \bigl(  |\xi|_{\gamma(y) } \bigr) \big| \gamma_{\alpha \beta}(y) |\xi||\eta| \bigr)\\
    \le & \ c \bigg[ \varepsilon + \dfrac{\bigl(  |\xi|_{\gamma(x)}-1 \bigr)^{p-1}}{|\xi|_{\gamma(x)}}   \biggr] |x-y| |\xi||\eta| +  c  \dfrac{\bigl[  |\xi|_{\gamma(x)}-1+ \bigl( |\xi|_{\gamma(y)}-1 \bigr)\bigr]^{p-1}}{|\xi|_{\gamma(x)}-1}  \big| |\xi|_{\gamma(x)}-|\xi|_{\gamma(y)}  \big| |\eta|\\
    \le & \ c  \bigg[ \varepsilon + \dfrac{\bigl[  |\xi|_{\gamma(x)}-1+ \bigl( |\xi|_{\gamma(y)}-1 \bigr)\bigr]^{p-1}}{|\xi|_{\gamma(x)}-1} \biggr]| x-y| |\xi||\eta| \\
    \le & \ c \bigl[ \varepsilon+\mu^{p-2}  \bigr] | x-y| |\xi||\eta|
\end{align*}
for a constant $c:=c(p,C_0,C_1,C_2)$.  
To derive the upper bound for $I_2$, we use Lemma \ref{lem5} and properties $(ii)$ and $(iii)$ of $\gamma$ and obtain
\begin{align*}
    I_2 \le & \ c  \big| h'_\varepsilon \bigl(  |\xi|_{\gamma(x) } \bigr)   |\xi|_{\gamma(x) } -
      h'_\varepsilon \bigl(  |\xi|_{\gamma(y) } \bigr)   |\xi|_{\gamma(y) } \big| |\xi||\eta| \\
      &+ \sum_{\alpha,\beta,\delta,k =1}^n \biggl( \big| h'_\varepsilon \bigl(  |\xi|_{\gamma(y) } \bigr)   \big| |\xi|_{\gamma(y) }     \bigg| \dfrac{\gamma_{\alpha \delta}(x)\gamma_{\beta k} (x) }{|\xi|_{\gamma(x)}^2} - \dfrac{\gamma_{\alpha \delta}(y)\gamma_{\beta k} (y) }{|\xi|_{\gamma(y)}^2}  \bigg| |\xi|^3 |\eta|
      \biggr)\\
      \le & \ c \bigl[  \bigl( |\xi|_{\gamma(x)}-1 \bigr)^{p-3}+ \bigl( |\xi|_{\gamma(y)}-1 \bigr)^{p-3}    \bigr] \big| |\xi|_{\gamma(x)}-|\xi|_{\gamma(y)}  \big| |\xi||\eta|\\
      & +   c \bigl( |\xi|_{\gamma(y)}-1 \bigr)^{p-2}
      \big| |\xi|_{\gamma(y)}^2 \gamma_{\alpha \delta}(x)\gamma_{\beta k} (x)-|\xi|_{\gamma(x)}^2\gamma_{\alpha \delta}(y)\gamma_{\beta k} (y)  \big| \dfrac{|\eta|}{|\xi|}\\
      \le & \ c \bigl[  \bigl( |\xi|_{\gamma(x)}-1 \bigr)^{p-2}+ \bigl( |\xi|_{\gamma(y)}-1 \bigr)^{p-2}    \bigr] \big| |x-y| |\xi|^2|\eta|\\
      &+ c \bigl( |\xi|_{\gamma(y)}-1 \bigr)^{p-2} \bigl[  |\gamma_{\alpha \delta}(x)-\gamma_{\alpha \delta}(y)| \gamma_{\beta k}(x)  + |\gamma_{\beta k}(x)-\gamma_{\beta k}(y)| \gamma_{\alpha \delta}(y) \bigr] |\xi|_{\gamma(y)} |\eta|\\
      &+ c \bigl( |\xi|_{\gamma(y)}-1 \bigr)^{p-2} \big| |\xi|^2_{\gamma(y)}-|\xi|^2_{\gamma(x)}  \big| \gamma_{\alpha \delta}(y)  \gamma_{\beta k}(y) \dfrac{|\eta|}{|\xi|}\\
      \le & \ c \mu^{p-2} |x-y||\xi|^2|\eta|,
\end{align*}
with $c:=c(p,C_0,C_1,C_2)$. Inserting the previous estimates in \eqref{i1i2}, we obtain the desired inequality.
\endproof

The proof of next two lemmas can be deduced similarly as in \cite[Lemmas 2.9 and 2.10]{BoDuGiPdN}.
\begin{lem}\label{lem9}
Let $\varepsilon , \delta \in (0,1]$ and $x \in \Omega$ and $\xi, \tilde{\xi}\in \mathbb{R}^{nN}$. Then, we have
$$\varepsilon^{\frac{1}{2}} |\xi-\tilde{\xi}|_{\gamma(x)} +|\mathcal{G}_\delta(x,\xi)-\mathcal{G}_\delta(x,\tilde{\xi})|^p_{\gamma(x)} \le \varepsilon^{\frac{1}{2}}|\xi|^2_{\gamma(x)}+c \varepsilon^{-\frac{1}{2}} \langle \textbf{A}_\varepsilon(x, \tilde{\xi})-\textbf{A}_\varepsilon(x,\xi), \tilde{\xi}-\xi \rangle_{\gamma(x)} ,$$
    for a constant $c:=c(p,\delta,C_0,C_1)$.
\end{lem}

In the following lemma we quantify the remainder term in the linearization of $\textbf{A}_\varepsilon$.

\begin{lem}\label{linearization}
    Let $\varepsilon \in [0,1]$ and $x \in \Omega$ and $\xi ,\tilde{\xi} \in \mathbb{R}^{nN}$, with $|\xi|_{\gamma(x)} \ge 1+\frac{1}{4} \mu$ and $|\xi|_{\gamma(x)}, |\tilde{\xi}|_{\gamma(x)} \le 1 +2 \mu$, for some $\mu >0$. Then, we have
    $$|\mathcal{B}_\varepsilon(x,\xi)(\tilde{\xi}-\xi,\eta)-  \langle \textbf{A}_\varepsilon(x, \tilde{\xi})-\textbf{A}_\varepsilon(x,\xi), \eta \rangle_{\gamma(x)}| \le c\mu^{p-3}|\tilde{\xi}-\xi|^2|\eta|,$$
    for every $\eta \in \mathbb{R}^{nN}$ and for a constant $c:=c(p,C_0,C_1)$.
\end{lem}

\begin{lem}\label{lem11}
    For any $\varepsilon \in [0,1]$, any ball $B_R \subset \mathbb{R}^n$ and any $v \in W^{2,2}_{loc}(B_R, \mathbb{R}^N)$, we have 
    $$|D \bigl[ g \bigl( |Dv|_{\gamma} \bigr) Dv \bigr]|^2 \le c \bigl[   \mathcal{A}_\varepsilon(x,Dv)(D^2v, D^2v)\bigl(  |Dv|_\gamma-1 \bigr)^p_+ + g'\bigl( |Dv|_\gamma \bigr)^2 |Dv|^4_{\gamma}\bigr],$$
    for a constant $c:=c(p,C_0,C_1,C_2)$.
\end{lem}
\proof For $\alpha,\beta \in \{ 1,...,n \}$ and $i \in \{ 1,...,N\}$, we compute
$$D_\alpha \bigl[  g \bigl( |Dv|_\gamma \bigr) D_\beta v^i \bigr] = g \bigl( |Dv|_\gamma \bigr)D_\alpha  D_\beta v^i + g' \bigl( |Dv|_\gamma \bigr) \dfrac{2F_\alpha+V_\alpha}{2 |Dv|_\gamma}  D_\beta v^i,$$
where we used the notation
$$F_\alpha:= \gamma_{k\lambda} D_k v^j D_\lambda D_\alpha v^j \quad \text{and} \quad V_\alpha :=  \dfrac{\partial \gamma_{k \lambda}}{\partial x_\alpha} D_k v^j D_\lambda v^j.$$
Then, taking square and summing over $\alpha,\beta \in \{ 1,...,n \}$ and $i \in \{ 1,...,N\}$, we get
\begin{align*}
    |D \bigl[  g \bigl( |Dv|_\gamma \bigr) D v |^2  \le 
    c \biggl[ &  g \bigl( |Dv|_\gamma \bigr)^2 |D^2v|^2+ g' \bigl( |Dv|_\gamma \bigr)^2 \dfrac{|Dv|^2|F|^2}{|Dv|^2_\gamma} 
    + g' \bigl( |Dv|_\gamma \bigr)^2 \dfrac{|Dv|^2|V|^2}{4|Dv|_\gamma^2}
    \biggr].
\end{align*}
Note that 
\begin{equation}\label{VF}
    |V| \le c(C_0,C_1,C_2) |Dv|^2 \quad \text{and} \quad |F| \le c(C_0,C_1) |Dv||D^2v|.
\end{equation}
Moreover, we can rewrite $|F|^2_\gamma$ as
\begin{align}
|F|^2_\gamma \ = \ \gamma_{\alpha\beta}F_\alpha F_\beta \ = & \ \gamma_{\alpha\beta}\gamma_{k\lambda}D_k v^j D_\lambda v^j \gamma_{\theta \rho}D_\theta v^i D_\beta D_\rho v^i \notag \\
=&  \ \bigl(  \gamma_{k \lambda}D_k v^j \gamma_{\theta \rho} D_\theta v^i \bigr) \bigl(  \gamma_{\alpha\beta} D_\alpha D_\lambda v^j D_\beta D_\rho v^j \bigr). \label{Fgamma}
\end{align}
An application of the first inequality in Lemma \ref{lem6}, property $(ii)$ of $\gamma$ and \eqref{VF} yield
\begin{align*}
    |D \bigl[  g \bigl( |Dv|_\gamma \bigr) D v |^2  \le 
    c \bigl[ &  h \bigl( |Dv|_\gamma \bigr) |D^2v|^2_\gamma(|Dv|_\gamma-1)^p_+ + g' \bigl( |Dv|_\gamma \bigr)^2 |F|^2 + g' \bigl( |Dv|_\gamma \bigr)^2 |Dv|_\gamma^4
    \bigr].
\end{align*}
We now take care of the quantity $g' \bigl( |Dv|_\gamma \bigr)^2 |F|^2$. Using the second inequality in Lemma \ref{lem6}, \eqref{VF} and \eqref{Fgamma} and recalling the definition of $\mathcal{A}_\varepsilon(x,Dv)$, we infer
\begin{align*}
    g' \bigl( |Dv|_\gamma \bigr)^2 |F|^2 \ = & \ g' \bigl( |Dv|_\gamma \bigr)^2 |Dv|_\gamma^2 \dfrac{|F|^2}{|Dv|_\gamma^2}\\
    \le & \ \dfrac{p^2}{p-1} \bigl[ h(|Dv|_\gamma ) +h'(|Dv|_\gamma ) |Dv|_\gamma   \bigr] \dfrac{|F|^2}{|Dv|_\gamma^2 } \bigl(  |Dv|_\gamma -1\bigr)^p_+ \\
    \le & \ c\biggl[ h(|Dv|_\gamma) |D^2v|_\gamma^2 \\
       & \quad \quad + h'(|Dv|_\gamma) |Dv|_\gamma \dfrac{  \gamma_{k \lambda}D_k v^j \gamma_{\theta \rho} D_\theta v^i }{|Dv|_\gamma^2}   \gamma_{\alpha\beta} D_\alpha D_\lambda v^j D_\beta D_\rho v^j   \biggr] \bigl( |Dv|_\gamma-1 \bigr)^p_+ \\
       = & \ c \mathcal{A}_\varepsilon(x,Dv)(D^2v, D^2v)\bigl( |Dv|_\gamma-1 \bigr)^p_+ ,
\end{align*}
for a constant $c:=c(p,C_0,C_1,C_2)$.
Inserting this above, we derive the desired estimate.
\endproof

\section{Approximating problems}
In order to prove Theorem \ref{mainthm}, we construct a sequence of functions $u_\varepsilon$ which are solutions to regularized problems.
We let $\varepsilon \in (0,1]$, $\mathbf{p}:=\max \{ p,2 \}$ and $u \in W^{1,p}_{loc}(\Omega, \mathbb{R}^N) \cap  W^{1,\infty}_{loc}(\Omega, \mathbb{R}^N) $ a weak solution of \eqref{EL}. By $u_\varepsilon \in u+ W^{1,\mathbf{p}}_0(B_R, \mathbb{R}^N) $ we denote the unique weak solution of the regularized
elliptic system
\begin{equation}\label{regsystem}
\begin{cases}
    \sum_{\alpha,\beta=1}^n D_\alpha \bigl[ h_\varepsilon\bigl(|Du_\varepsilon|_\gamma\bigr)\gamma_{\alpha\beta}D_\beta u_\varepsilon^i \bigr]=0, \quad \forall i =1,...,N \  &\text{in} \ B_R \\
    u_\varepsilon =u \ & \text{on} \ \partial B_R
    \end{cases}
\end{equation}
The weak formulation of \eqref{regsystem} is
\begin{equation}\label{regsystem2}
    \int_{B_R} \langle \textbf{A}_\varepsilon(x,Du_\varepsilon), D \varphi \rangle_\gamma \ \mathrm{d}x=0, 
\end{equation}
for every $\varphi \in W^{1,\mathbf{p}}_0(B_R, \mathbb{R}^N)$. 

Now, we will show that $u_\varepsilon$ is Lipschitz continuous and possesses weak second derivatives. Besides, we have a uniform energy bound for $Du_\varepsilon$ with respect to $\varepsilon$.

\begin{lem}\label{energyDu}
    For any $\varepsilon \in (0,1]$ we have $u_\varepsilon \in W^{1,\infty}_{loc} ( B_R, \mathbb{R}^N ) \cap  W^{2,2}_{loc} ( B_R, \mathbb{R}^N )$. Moreover, for any ball $B_{2 \rho}(x_0) \Subset B_R$ the uniform $W^{1,\infty}$-bound
    \begin{equation}
        \sup_{B_\rho(x_0)} |Du_\varepsilon| \le c \biggl(  \fint_{B_{2\rho}(x_0)} (1+|Du_\varepsilon|^p) \ \mathrm{d}x \biggr)^\frac{1}{p} \label{boundue}
    \end{equation}
    and the quantitative $W^{2,2}$-estimate
    $$\int_{B_{\rho/2}(x_0)} |D^2u_\varepsilon|^2 \ \mathrm{d}x \le \dfrac{c}{\rho^2} \int_{B_{2\rho}(x_0)} \biggl[  |Du_\varepsilon|^2+\dfrac{1}{\varepsilon^2} |Du_\varepsilon|^{2(p-1)} \biggr] \ \mathrm{d} x$$
    hold true, for a constant $c:=c(n,p,C_0,C_1,C_2)$.
\end{lem}
\proof The proof of the $W^{1,\infty}$-estimate for the gradient can be inferred from \cite[Theorem 2.7]{fonseca}.
The $W^{2,2}$-estimate follows in standard way from the difference quotient method. We leave the details to the reader.
\endproof

We have the following uniform energy bound for $Du_\varepsilon$, whose proof follows in a usual way (see e.g.\ \cite{BoDuGiPdN}).

\begin{lem}\label{uniformbound}
    There exists a constant $c:=c(p,C_0,C_1)$ such that for any $\varepsilon \in (0,1]$ we have
    $$\int_{B_R} \bigl(  |Du_\varepsilon|^p + \varepsilon |Du_\varepsilon|^2  \bigr) \ \mathrm{d} x \le c \int_{B_R} \bigl(  |Du|^p + \varepsilon |Du|^2 +1 \bigr) \ \mathrm{d} x. $$
\end{lem}

\begin{rmk}\label{rmkub}
    Combining Lemma \ref{uniformbound}, inequality \eqref{boundue} and the Lipschitz continuity of $u$, we have the boundedness of $Du_\varepsilon$ on any ball $B_{2R} \Subset \Omega$ uniformly with respect to $\varepsilon \in (0,1]$. More precisely, there exists a constant 
\begin{equation*}
    M:=M(n,N,p,C_0,C_1,C_2,R,\Vert Du \Vert_{L^p(B_{2R})}),
\end{equation*}
independent of $\varepsilon$, such that $$\sup_{x \in B_{R}}|Du_\varepsilon|_\gamma \le M.$$

\end{rmk}

The following lemma ensures strong convergence of the approximating solutions in the sense that $\mathcal{G}_\delta(x,Du_\varepsilon)$ converges to  $\mathcal{G}_\delta(x,Du)$ in $L^p$.
\begin{lem}\label{lem3.3}
    Let $\delta \in (0,1]$ and $u_\varepsilon$, with $\varepsilon \in (0,1]$, be the unique weak solution of \eqref{regsystem}. Then, we have
    $$ \mathcal{G}_\delta(x,Du_\varepsilon) \rightarrow \mathcal{G}_\delta(x,Du) \quad \text{in} \ L^p(B_R, \mathbb{R}^{nN}) \ \text{as} \ \varepsilon \rightarrow0^+.$$
\end{lem}
\proof Testing \eqref{regsystem2} and the weak formulation of \eqref{EL} with $\varphi= u_\varepsilon-u$, we have
$$\int_{B_R} \langle \textbf{A}_\varepsilon(x,Du_\varepsilon)-\textbf{A}_\varepsilon(x,Du), Du_\varepsilon-Du \rangle_{\gamma} \ \mathrm{d}x = \varepsilon \int_{B_R} \langle Du, Du-Du_\varepsilon \rangle_\gamma \ \mathrm{d}x.$$
Using Lemma \ref{lem9} and Young's inequality, we find 
\begin{align*}
    \int_{B_R} \bigl[ \varepsilon^{\frac{1}{2}} |Du_\varepsilon-Du|^2_\gamma + & |\mathcal{G}_\delta(x,Du_\varepsilon)-\mathcal{G}_\delta(x,Du)|^p_\gamma\bigr] \ \mathrm{d}x \\
    \le & \ \varepsilon^{\frac{1}{2}} \int_{B_R} |Du|^2_\gamma \ \mathrm{d}x
    +c \varepsilon^{-\frac{1}{2}} \int_{B_R}  \langle \textbf{A}_\varepsilon(x,Du_\varepsilon)-\textbf{A}_\varepsilon(x,Du), Du_\varepsilon-Du \rangle_{\gamma} \ \mathrm{d}x \\
    \le & \  \varepsilon^{\frac{1}{2}} \int_{B_R} |Du|^2_\gamma \ \mathrm{d}x +  c\varepsilon^{\frac{1}{2}} \int_{B_R} |Du|_\gamma |Du_\varepsilon-Du|_\gamma \ \mathrm{d}x\\
    \le & \  \varepsilon^{\frac{1}{2}} \int_{B_R} |Du_\varepsilon-Du|^2_\gamma \ \mathrm{d}x +c  \varepsilon^{\frac{1}{2}} \int_{B_R} |Du|^2_\gamma \ \mathrm{d}x.
\end{align*}
Hence, we obtain for a constant $c:=c(p,\delta,C_0,C_1)$
\begin{equation*}
    \int_{B_R}   |\mathcal{G}_\delta(x,Du_\varepsilon)-\mathcal{G}_\delta(x,Du)|^p_\gamma \ \mathrm{d}x \le c  \varepsilon^{\frac{1}{2}} \int_{B_R} |Du|^2_\gamma \ \mathrm{d}x,
\end{equation*}
which is finite since $Du \in L^\infty_{loc}(\Omega,\mathbb{R}^{nN})$. Therefore, we can conclude that $\mathcal{G}_\delta(x,Du_\varepsilon) \rightarrow \mathcal{G}_\delta(x,Du) $ strongly in $ L^p(B_R, \mathbb{R}^{nN}) $ as $\varepsilon \rightarrow0^+$.
\endproof

\section{Main result}
In this section we will prove Theorem \ref{mainthm} assuming that Propositions \ref{ND} and \ref{D} below hold true. We refer to Sections \ref{proofND} and \ref{secD} for the proof of those two propositions.

We denote by $u \in W^{1,p}_{loc}(\Omega, \mathbb{R}^N) $ a weak solution of \eqref{EL}. In the following we will assume that $u \in W^{1,\infty}_{loc}(\Omega, \mathbb{R}^N) $.

\subsection{H\"older continuity of $\mathcal{G}_\delta(x,Du_\varepsilon)$}
 In this subsection we will prove that $\mathcal{G}_\delta(x,Du_\varepsilon)$ is locally H\"older continuous for any $\delta \in (0,1]$.
\\We consider a ball $B_R \Subset \Omega$ and we denote by $u_\varepsilon$ the unique weak solution of the regularized problem \eqref{regsystem}. For $0< r_0 <R$, we let $r_1:= \frac{1}{2}\bigl( R+r_0 \bigr)$. Then, by Remark \ref{rmkub} we have that there exists a constant 
\begin{equation}\label{constant}
    M:=M(n,N,p,C_0,C_1,C_2,R-r_0,\Vert Du \Vert_{L^p(B_R)}),
\end{equation}
independent of $\varepsilon $, such that $\sup_{x \in B_{r_1}}|Du_\varepsilon|_\gamma \le M$. Now, we consider $x_0 \in B_{r_0}$ and a radius $\rho \le r_1$ such that $B_{2 \rho}(x_0) \subset B_{r_1}$. Then, we have
\begin{equation}\label{bounded}
    \sup_{x \in B_{2\rho}(x_0)}|Du_\varepsilon|_\gamma \le 1+\delta+\mu
\end{equation}
for some $\mu >0$ such that
\begin{equation}\label{delta}
    1+\delta+\mu \le M.
\end{equation}
Next, for $\nu \in (0,1)$ we define the super-level set of $|Du_\varepsilon|_\gamma$ by
\begin{equation*}
    E^\nu_\rho(x_0):= \bigl\{  x \in B_\rho(x_0): |Du_\varepsilon(x)|_{\gamma(x)} -1-\delta >(1-\nu)\mu \bigr\}.
\end{equation*}

\begin{prop}\label{ND}
    Let $\beta \in (0,1)$ and $\varepsilon, \delta \in (0,1]$ and 
    \begin{equation}\label{mu}
        0 < \delta < \mu.
    \end{equation}
    Then, there exist $\nu:=\nu(n,N,p,M,\delta,C_0,C_1,C_2)\in \bigl(0, \frac{1}{8} \bigr]$ and $\hat{\rho}:=\hat{\rho}(n,N,p,M,\delta,C_0,C_1,C_2)\in (0,1]$ such that there holds: whenever $B_{2\rho}(x_0) \subset B_{r_1}$ is a ball with radius $\rho \le \hat{\rho}$ and center $x_0 \in B_{r_0}$, and $u_\varepsilon$ is the unique weak solution to the regularized system \eqref{regsystem} and hypothesis \eqref{bounded} and \eqref{delta} and the measure condition
    \begin{equation}\label{measureND}
        |B_\rho(x_0) \setminus E^\nu_\rho(x_0)| < \nu |B_\rho(x_0)|
    \end{equation}
    are satisfied, then the limit
    \begin{equation}\label{limit}
        \Gamma_{x_0} := \lim_{r \rightarrow 0} \bigl( \mathcal{G}_{2 \delta}(x,Du_\varepsilon) \bigr)_{x_0,r}
    \end{equation}
    exists, and the excess decay estimate
    \begin{equation}\label{decay} 
        \fint_{B_r(x_0)} |\mathcal{G}_{2 \delta}(x,Du_\varepsilon)-\Gamma_{x_0}|^2 \ \mathrm{d}x
        \le \biggl(  \dfrac{r}{\rho} \biggr)^{2 \beta} \mu^2
 \end{equation}
 holds true for any $0 < r \le \rho$. Moreover, we have
 \begin{equation}\label{gammanormx0}
     |\Gamma_{x_0}|_{\gamma(x_0)} \le \mu.
 \end{equation}
\end{prop}
The statement for the degenerate regime is as follows.
\begin{prop}\label{D}
    Let $\varepsilon, \delta \in (0,1], \mu >0$ and $\nu \in \bigl(0, \frac{1}{4} \bigr]$. Then, there exist constants $k\in \bigl[ 2^{-\frac{1}{2}},1 \bigr)$ and $\hat{c} \ge 1$ depending at most on $n,p,M,\delta, \nu , C_0,C_1,C_2$ such that there holds: whenever $B_{2\rho}(x_0) \subset B_{r_1}$ is a ball with radius $\rho \le \hat{\rho}$ and center $x_0 \in B_{r_0}$, and $u_\varepsilon$ is the unique weak solution to the regularized system \eqref{regsystem} and hypothesis \eqref{bounded} and \eqref{delta} and the measure condition
    \begin{equation}\label{measureD}
        |B_\rho(x_0) \setminus E^\nu_\rho(x_0)| \ge \nu |B_\rho(x_0)|
    \end{equation}
    are satisfied, then, either $\mu^2 < \hat{c} \rho$ or 
    $$\sup_{B_{2\rho}(x_0)}|\mathcal{G}_{ \delta}(x,Du_\varepsilon)|_\gamma \le k \mu$$
    hold true.
\end{prop}
Now, we are able to prove the main result of this subsection.
\begin{thm}\label{thmholder}
    Let $\varepsilon,\delta \in (0,1]$ and $u_\varepsilon$ be the unique solution of the regularized system \eqref{regsystem}. Then, $\mathcal{G}_\delta(x,Du_\varepsilon)$ is H\"older continuous in $B_{r_0}$, for any $0< r_0 <R$, with H\"older exponent $\alpha_\delta \in (0,1)$ and a H\"older constant $c_\delta$ both depending at most on $n,N,p,M,\delta,C_0,C_1$ and $C_2$.
\end{thm}
\proof By $\nu \in \bigl( 0, \frac{1}{8} \bigr]$ and $\hat{\rho} \in (0,1]$ we denote the constants from Proposition \ref{ND} applied with $\beta = \frac{1}{2}$ and by $k \in \bigl[ 2^{-\frac{1}{2}},1 \bigr)$ and $\hat{c} \ge 1$ we denote the ones from Proposition \ref{D}. We let $\mu= M-1-\delta$ and
$$\rho_*:= \min \biggl\{  \hat{\rho}, \dfrac{(k \mu)^2}{\hat{c}} \biggr\}.$$
We consider a ball $B_{2 \rho}(x_0) \subset B_{r_1}$ with center $x_0 \in B_{r_0}$ and $\rho \le \rho_*$. We have that \eqref{constant}, \eqref{bounded} and \eqref{delta} are satisfied on this ball. Our aim is to prove that $\mathcal{G}_{2\delta}(x,Du_\varepsilon)$ is H\"older continuous in $B_{r_0}$ with H\"older exponent
$$\alpha:= -\dfrac{\text{log} \ k}{\text{log} \ 2} \in \biggl(0, \dfrac{1}{2} \biggr].$$
We divide the proof  of the theorem in two parts.\\
\textbf{Step 1.} We prove that
 the limit
$$  \Gamma_{x_0} := \lim_{r \rightarrow 0} \bigl( \mathcal{G}_{2 \delta}(x,Du_\varepsilon) \bigr)_{x_0,r}$$
exists and 
    \begin{equation}\label{DEC}
        \fint_{B_r(x_0)} |\mathcal{G}_{2 \delta}(x,Du_\varepsilon)-\Gamma_{x_0}|^2 \ \mathrm{d}x
        \le c \biggl(  \dfrac{r}{\rho} \biggr)^{2 \alpha} \mu^2
    \end{equation}
holds true for a constant $c:=c(C_0)$.

We define for $i \in \mathbb{N}_0$ radii
$$\rho_i:= 2^{-i} \rho \quad \text{and} \quad \mu_i:=k^i\mu$$
and observe that
\begin{equation}\label{mui}
    \mu_i=k^i\mu \le 2^{-\alpha i}\mu = \biggl( \dfrac{\rho_i}{\rho} \biggr)^\alpha \mu
    \end{equation}
for any $i \in \mathbb{N}_0$. Now, suppose that assumption \eqref{measureD} holds on $B_\rho(x_0)$. Then, Proposition \ref{D} yields that either $\mu^2 < \hat{c} \rho$ or 
$$\sup_{B_{\rho_1}(x_0)}|\mathcal{G}_\delta(x,Du_\varepsilon)|_\gamma \le k \mu = \mu_1.$$
Note that the first alternative cannot happen, since it would imply
$$\mu^2 < \hat{c} \rho \le \hat{c} \rho_* \le k^2 \mu^2.$$
Hence, we conclude that \eqref{bounded} holds on $B_{\rho_1}(x_0)$ with $\mu=\mu_1$. If the measure theoretic condition \eqref{measureD} is satisfied with $\rho=\rho_1$ and $\mu=\mu_1$, then a second application of Proposition \ref{D} yields that either $\mu_1^2 < \hat{c}\rho_1$ or
$$\sup_{B_{\rho_2}(x_0)}|\mathcal{G}_\delta(x,Du_\varepsilon)|_\gamma \le k \mu_1 = \mu_2.$$
As before,  the first alternative cannot happen, since it would imply
$$\mu^2_1 < \hat{c} \rho_1 \le 2^{-1} (k \mu)^2 \le k^2 \mu^2_1.$$
Assume that \eqref{measureD} is satisfied for $i=1,...,i_0-1$ up to some $i_0 \in \mathbb{N}$, i.e.\ that \eqref{measureD} holds true on the balls $B_{\rho_i}(x_0)$ with $\mu=\mu_i$. Then, we can conclude that 
\begin{equation}\label{sup}
    \sup_{B_{\rho_i}(x_0)}|\mathcal{G}_\delta(x,Du_\varepsilon)|_\gamma \le \mu_i
\end{equation}
for $i=0,...,i_0$. Now assume that \eqref{measureD} fails to hold for some $i_0 \in \mathbb{N}_0$. If $\mu_{i_0} > \delta$, the hypothesis of Proposition \ref{ND} are satisfied on $B_{\rho_{i_0}}(x_0)$ and we can conclude that the limit
    \begin{equation}
        \Gamma_{x_0} := \lim_{r \rightarrow 0} \bigl( \mathcal{G}_{2 \delta}(x,Du_\varepsilon) \bigr)_{x_0,r}
    \end{equation}
    exists and that 
    \begin{equation}\label{Gdelta} 
        \fint_{B_r(x_0)} |\mathcal{G}_{2 \delta}(x,Du_\varepsilon)-\Gamma_{x_0}|^2 \ \mathrm{d}x
        \le   \dfrac{r}{\rho_{i_0}}  \mu^2_{i_0}
 \end{equation}
for any $0< r \le \rho_{i_0}$. Moreover, we have
\begin{equation}\label{Gamma}
    |\Gamma_{x_0}|_{\gamma(x_0)} \le \mu_{i_0}.
\end{equation}
Therefore, we obtain from \eqref{Gdelta} and \eqref{mui} that 
 \begin{equation}
        \fint_{B_r(x_0)} |\mathcal{G}_{2 \delta}(x,Du_\varepsilon)-\Gamma_{x_0}|^2 \ \mathrm{d}x
        \le   \dfrac{r}{\rho_{i_0}}  \biggl( \dfrac{\rho_{i_0}}{\rho} \biggr)^{2\alpha} \mu^2 \le \biggl( \dfrac{r}{\rho} \biggr)^{2\alpha} \mu^2
 \end{equation}
holds true for any $0 < r \le \rho_{i_0} $. For a radius $r \in ( \rho_{i_0}, \rho ]$, there exists $i \in \{  0,...,i_0 \}$ such that $\rho_{i+1}< r \le \rho_i$. Using  \eqref{mui}, \eqref{sup} and \eqref{Gamma}, we obtain
\begin{align*}
    \fint_{B_r(x_0)} |\mathcal{G}_{2 \delta}(x,Du_\varepsilon)-\Gamma_{x_0}|^2 \ \mathrm{d}x \le & \ 2C_0^{-1} \sup_{B_{\rho_1}(x_0)}|\mathcal{G}_\delta(x,Du_\varepsilon)|^2_\gamma +2 C^{-1}_0 |\Gamma_{x_0}|^2_{\gamma(x_0)}\\
    \le & \ 4 C_0^{-1} \mu_i^2 \le 4 C_0^{-1} \biggl(  \dfrac{\rho_i}{\rho} \biggr)^{2 \alpha} \mu^2 \le 8 C_0^{-1} \biggl(  \dfrac{r}{\rho} \biggr)^{2 \alpha} \mu^2.
\end{align*}
Combining the preceding two inequalities, we have shown \eqref{DEC} provided $\mu_{i_0} > \delta$. 
\\In the case $\mu_{i_0} \le \delta$, thanks to \eqref{sup} we have that $\mathcal{G}_{2\delta}(x,Du_\varepsilon)=0$ on the ball $B_{\rho_{i_0}}(x_0)$. Combining this with \eqref{sup} and keeping in mind that
$$|\mathcal{G}_{2 \delta}(x,Du_\varepsilon)|_{\gamma(x)} \le |\mathcal{G}_{ \delta}(x,Du_\varepsilon)|_{\gamma(x)}, $$
we obtain
\begin{equation}\label{supgdelta}
    \sup_{B_{\rho_i}(x_0)}|\mathcal{G}_{2\delta}(x,Du_\varepsilon)|_\gamma \le \mu_i, \quad \forall i \in \mathbb{N}_0.
\end{equation}
In the last case when \eqref{measureD} holds for any $i \in \mathbb{N}$, then \eqref{sup} is satisfied for any $i \in \mathbb{N}$ and hence we obtain \eqref{supgdelta} also in this case. Moreover, \eqref{supgdelta} implies
$$\Gamma_{x_0} = \lim_{r \rightarrow 0} \bigl( \mathcal{G}_{2 \delta}(x,Du_\varepsilon) \bigr)_{x_0,r}=0.$$
For $r \in (0,\rho]$, we find $i \in \{0,1,...\}$ such that $\rho_{i+1} < r \le \rho_i$. Then, \eqref{supgdelta} and \eqref{mui} imply
\begin{align*}
    \fint_{B_r(x_0)} |\mathcal{G}_{2 \delta}(x,Du_\varepsilon)-\Gamma_{x_0}|^2 \ \mathrm{d}x \le & \ C_0^{-1} \sup_{B_{\rho_1}(x_0)}|\mathcal{G}_\delta(x,Du_\varepsilon)|^2_\gamma \\
    \le \ &  C_0^{-1} \mu_i^2 \le  C_0^{-1} \biggl(  \dfrac{\rho_i}{\rho} \biggr)^{2 \alpha} \mu^2 \le 2 C_0^{-1} \biggl(  \dfrac{r}{\rho} \biggr)^{2 \alpha} \mu^2.
\end{align*}
Thus, \eqref{DEC} is satisfied also in the remaining case.
\vspace{0.5cm}
\\
\textbf{Step 2.} The excess-decay estimate \eqref{DEC} implies the H\"older continuity in $B_{r_0}$ of the Lebesgue representative  $x \mapsto \Gamma_x$ of $\mathcal{G}_{2 \delta}(x,Du_\varepsilon)$. The proof is standard and it will be not presented here.
\endproof

\subsection{Proof of Theorem \ref{mainthm}}
Now we are in position to give the proof of our main
result.

\begin{proof}[{\it Proof of Theorem \ref{mainthm}.}]
Let $\varepsilon \in (0,1]$ and consider a fixed ball $B_R:=B_R(y_0) \subset \Omega$. By $u_\varepsilon$ we denote the weak solution to \eqref{regsystem}. Next, we fix $\delta \in (0,1]$ and $r \in (0,R)$. From Theorem \ref{mainthm} we know that $\mathcal{G}_\delta(x,Du_\varepsilon)$ is H\"older continuous in $\bar{B}_r$ with H\"older exponent $\alpha_\delta \in (0,1)$ and constant $c_\delta$ both depending at most on $n,N,p,M,\delta,C_0,C_1$ and $C_2$. From Lemma \ref{lem3.3}, we know that $\mathcal{G}_\delta(x,Du_\varepsilon) \rightarrow \mathcal{G}_\delta(x,Du)$ in $L^p(B_R, \mathbb{R}^{nN})$ as $\varepsilon \rightarrow 0^+$. This implies that there exists a subsequence $(\varepsilon_i)_{i \in \mathbb{N}}$ such that $\varepsilon_i \rightarrow 0^+$ as $i \rightarrow + \infty$ and $\mathcal{G}_\delta(x,Du_\varepsilon) \rightarrow \mathcal{G}_\delta(x,Du)$ a.e.\ in $B_R$. On the other hand, by Ascoli-Arzelà Theorem we can conclude that $\mathcal{G}_\delta(x,Du_\varepsilon) $ converges uniformly on compact subsets of $B_R$. Therefore, the limit function $\mathcal{G}_\delta(x,Du)$ is H\"older continuous in $\bar{B}_r$ for any $\delta \in (0,1]$, with H\"older exponent $\alpha_\delta \in (0,1)$ and constant $c_\delta$.
\\Moreover, we have $\mathcal{G}_\delta(x,Du) \rightarrow \mathcal{G}(x,Du)$ uniformly in $\bar{B}_r$ as $\delta \rightarrow 0^+$. Indeed, it holds
\begin{align*}
    |\mathcal{G}_\delta(x,Du)- & \mathcal{G}(x,Du)| \\
    \le & \ c |\mathcal{G}_\delta(x,Du)-\mathcal{G}(x,Du)|_{\gamma(x)} \\
    = & \ c \biggl|  \bigl( |Du(x)|_{\gamma(x)}-1-\delta \bigr)_+ \dfrac{Du}{|Du(x)|_{\gamma(x)}} -  \bigl( |Du(x)|_{\gamma(x)}-1 \bigr)_+ \dfrac{Du}{|Du(x)|_{\gamma(x)}}   \biggr|_{\gamma(x)} \\
    \le & \ c |  \bigl( |Du(x)|_{\gamma(x)}-1-\delta \bigr)_+  -  \bigl( |Du(x)|_{\gamma(x)}-1 \bigr)_+   |_{\gamma(x)}  \le \ c \delta
\end{align*}
for every $x \in \bar{B}_r$ and for a constant $c:=c(C_0,C_1)$. As the uniform limit of a sequence of continuous functions, $\mathcal{G}(x,Du)$ is continuous on $\bar{B}_r$. Moreover, $\mathcal{G}(x,Du)$ is uniformly continuous on $\bar{B}_r$.

Now, let $\mathcal{K}: \Omega \times \mathbb{R}^{nN} \rightarrow \mathbb{R}$ be any continuous function vanishing on the set $\bigl\{ (x,\xi ) \in \Omega \times \mathbb{R}^{nN} : |\xi|_{\gamma(x)} \le 1  \bigr\}$. Since $u \in W^{1,\infty}_{loc}(\Omega, \mathbb{R}^N)$, we find $M>0$ such that $\sup_{\bar{B}_r}|Du|_\gamma \le M$. By $\omega : \mathbb{R}_+ \rightarrow \mathbb{R}_+$ we denote the modulus of continuity of $\mathcal{K}$ on $\bigl\{ (x,\xi ) \in \bar{B}_r \times \mathbb{R}^{nN} : |\xi|_{\gamma(x)} \le M  \bigr\}$, i.e. for any $(x,\xi), (y,\eta)\in \bar{B}_r \times \mathbb{R}^{nN}$ with $|\xi|_{\gamma(x)}, |\eta|_{\gamma(y)} \le M$ we have 
$$|\mathcal{K}(x,\xi)-\mathcal{K}(y,\eta)| \le \ \omega \bigl(  |x-y|+|\xi-\eta| \bigr).$$
Next, given $\varepsilon \in (0,1)$ we choose $\delta >0$ such that
\begin{equation}\label{uniformG}
    |\mathcal{G}(x,Du(x))-\mathcal{G}(y,Du(y))| < \varepsilon, \quad \forall x,y, \in \bar{B}_r \ \text{with} \ |x-y | < \delta.
\end{equation}
Now, we distinguish two cases. Let assume that $|Du(x)|_{\gamma(x)} \le 1 + \sqrt{\varepsilon}$. If $|Du(x)|_{\gamma(x)} > 1$, using the fact that $\mathcal{K}=0$ on the set $\bigl\{ (x,\xi ) \in \Omega \times \mathbb{R}^{nN} : |\xi|_{\gamma(x)} \le 1  \bigr\}$, we conclude that
\begin{align*}
    | \mathcal{K} \bigl(  x,Du(x) \bigr)|  = & \ \bigg|  \mathcal{K} \bigl(  x,Du(x) \bigr)- \mathcal{K} \biggl(  x,\dfrac{Du(x)}{|Du(x)|_{\gamma(x)}} \biggr)  \bigg| \\
    \le & \ \bigg| \omega \biggl(  \bigg|  Du(x)- \dfrac{Du(x)}{|Du(x)|_{\gamma(x)}}  \bigg| \biggr)   \bigg| \le \ \omega \bigl(  \sqrt{C_0^{-1} \varepsilon} \bigr).
\end{align*}
The previous inequality holds trivially if $|Du(x)|_{\gamma(x)} \le 1$. Moreover, by \eqref{uniformG} we have
\begin{align*}
    \bigl(  |Du(y)|_{\gamma(y)}-1 \bigr)_+ = & \ |\mathcal{G}(y,Du(y))|_{\gamma(y)} \\
    \le & \  \sqrt{C_1}|\mathcal{G}(y,Du(y))| \\
    \le & \  \sqrt{C_1}|\mathcal{G}(y,Du(y))-\mathcal{G}(x,Du(x))| + \sqrt{C_0^{-1}C_1}|\mathcal{G}(x,Du(x))|_{\gamma(x)}\\
    \le & \ \sqrt{C_1} \varepsilon + \sqrt{C_0^{-1}C_1} \bigl(  |Du(x)|_{\gamma(x)}-1 \bigr)_+ \\
    \le & \ C \bigl(  \varepsilon + \sqrt{\varepsilon} \bigr) \le \ 2 C \sqrt{\varepsilon},
\end{align*}
for a constant $C:=C(C_0,C_1)$. This implies $|Du(y)|_{\gamma(y)} \le 1+2C \sqrt{\varepsilon}$. Similarly as above we conclude
$$ | \mathcal{K} \bigl(  y,Du(y) \bigr)| \le \omega \bigl(  C \sqrt{\varepsilon} \bigr). $$
Combining the estimates from above, we end up with
$$ | \mathcal{K} \bigl(  x,Du(x) \bigr) - \mathcal{K} \bigl(  y,Du(y) \bigr)| \le 2 \omega \bigl(  C \sqrt{\varepsilon} \bigr), $$
where $C:=C(C_0,C_1)$.
On the other hand, if $|Du(x)|_{\gamma(x)} > 1 + \sqrt{\varepsilon}$,
Lemma \ref{lem3} and \eqref{uniformG} imply
\begin{align*}
    |Du(x)-Du(y)| \le & \ c \biggl( 1 + \dfrac{2}{\sqrt{\varepsilon}} \biggr)|\mathcal{G}(x,Du(x))-\mathcal{G}(y,Du(y))|  \\
    \le & \ c \bigl(  \varepsilon + 2 \sqrt{\varepsilon} \bigr) \le \ 3 c \sqrt{\varepsilon}.
\end{align*}
Diminishing $\delta $ if necessary, we get
$$| \mathcal{K} \bigl(  x,Du(x) \bigr) - \mathcal{K} \bigl(  y,Du(y) \bigr)| \le 2 \omega \bigl(  c \sqrt{\varepsilon} \bigr),$$
for a constant $c:=c(C_0,C_1)$. 
\\Hence, $\mathcal{K}(x,Du)$ is continuous on $\bar{B}_r$, for any $\bar{B}_r \Subset B_R \Subset \Omega$. This implies that $\mathcal{K}(x,Du)$ is continuous in $\Omega$. 
\end{proof}

\section{Estimates for second order derivatives}

In the following we will differentiate the regularized elliptic system \eqref{regsystem}.
The resulting integral estimate will subsequently
be used in different directions. On the one hand, it implies several variants of energy estimates
for second derivatives of $u_\varepsilon$. On the other hand, it implies that a certain function of $|Du_\varepsilon|_\gamma$ is a sub-solution
of an elliptic equation with measurable coefficients.

\subsection{The main integral inequality}

Let $u_\varepsilon \in u+ W^{1,\mathbf{p}}_0(B_R, \mathbb{R}^N) $ be the unique weak solution of \eqref{regsystem}. In the following, we write $u$ instead of $u_\varepsilon$ for the sake of simplicity. Let 
$ \Phi \in W^{1,\infty}_{loc} \bigl ( \mathbb{R}_{\ge 0}, \mathbb{R}_{\ge 0} \bigr) $ be a non-negative increasing function and $\xi \in \mathcal{C}^1_0(B_R)$ be a non-negative testing function. Testing \eqref{regsystem2} with the function
$$\varphi^i:= D_\nu \bigl[  \xi \gamma_{\nu \sigma} D_\sigma u^i \Phi \bigl( |Du|_\gamma\bigr)\bigr], \quad i=1,...,N$$
and performing an integration by parts we obtain
\begin{align}
    0  = &  \int_{B_R} h_\varepsilon \bigl( |Du|_\gamma \bigr) \gamma_{\alpha\beta} D_\alpha u^i D_\beta D_\nu\bigl[  \xi \gamma_{\nu \sigma} D_\sigma u^i \Phi \bigl( |Du|_\gamma  \bigr) \bigr] \ \mathrm{d}x \notag \\
    =& - \int_{B_R} D_\nu \bigl[ h_\varepsilon \bigl( |Du|_\gamma \bigr) \gamma_{\alpha\beta} D_\alpha u^i \bigr] D_\beta \bigl[  \xi \gamma_{\nu \sigma} D_\sigma u^i \Phi \bigl( |Du|_\gamma  \bigr) \bigr] \ \mathrm{d}x \label{6.1}
\end{align}
For every $i \in \{  1,...,N \}$ and $\beta,\nu \in \{  1,...,n \}$, we have
\begin{align*}
    D_\nu \bigl[ h_\varepsilon \bigl( |Du|_\gamma \bigr)\gamma_{\alpha\beta} D_\alpha u^i \bigr] 
    &=  h_\varepsilon \bigl( |Du|_\gamma \bigr)  \gamma_{\alpha\beta} D_\nu D_\alpha u^i +  h'_\varepsilon \bigl( |Du|_\gamma \bigr) |Du|_\gamma  \dfrac{\gamma_{\alpha\beta} D_\alpha u^i}{|Du|^2_\gamma} F_\nu\\
     & \ \  \ +  h'_\varepsilon \bigl( |Du|_\gamma \bigr) |Du|_\gamma \dfrac{\gamma_{\alpha\beta}D_\alpha u^i}{2|Du|_\gamma^2}  \dfrac{\partial\gamma_{k \lambda}}{\partial x_\nu} D_k u^j D_\lambda u^j\\
    & \ \ \ + h_\varepsilon \bigl( |Du|_\gamma \bigr) \dfrac{\partial \gamma_{\alpha\beta}}{\partial x_\nu} D_\alpha u^i \\
    & =:  I_1+I_2+I_3+I_4
\end{align*}
where we denoted $F_\nu := \gamma_{k\lambda} D_k u^j  D_\nu D_\lambda u^j$. Next, we consider the second factor in \eqref{6.1}. For every $i \in \{ 1,...,N \} $ and $\beta,\nu \in \{ 1,...,n \}$, we obtain
\begin{align*}
    D_\beta \bigl[  \xi \gamma_{\nu \sigma} D_\sigma u^i \Phi \bigl( |Du|_\gamma  \bigr) \bigr] 
    &= \gamma_{\nu\sigma} D_\beta D_\sigma u^i \Phi \bigl( |Du|_\gamma \bigr) \xi + \gamma_{\nu \sigma} D_\sigma u^i F_\beta \dfrac{\Phi' \bigl( |Du|_\gamma \bigr)}{|Du|_\gamma} \xi\\
    & \ \ \ + \gamma_{\nu \sigma} D_\sigma u^i \dfrac{\Phi' \bigl( |Du|_\gamma \bigr)}{2|Du|_\gamma} V_\beta \xi + \dfrac{\partial \gamma_{\nu\sigma}}{\partial x_\beta} D_\sigma u^i \Phi \bigl( |Du|_\gamma \bigr) \xi\\
    & \ \ \ + \gamma_{\nu \sigma} D_\sigma u^i \Phi \bigl( |Du|_\gamma \bigr)  \dfrac{ \partial \xi}{\partial x_\beta}\\
    &=: J_1 + J_2 +J_3+J_4+J_5
\end{align*}
where we denoted $V_\beta := \dfrac{\partial \gamma_{k \lambda}}{\partial x_\beta} D_k u^j D_\lambda u^j$.
Using these notations, we can write \eqref{6.1} as follows
\begin{equation}\label{6.2}
    \sum_{k=1}^4 \sum_{l=1}^5 \int_{B_R} I_k \cdot J_k \  \mathrm{d}x =0.
\end{equation}
In the sequel we evaluate or estimate separately the terms $I_k \cdot J_l$. We start with $\bigl( I_1 + I_2 \bigr) \cdot J_1$, which we rewrite in the form
\begin{align*}
    \dfrac{\bigl( I_1 + I_2 \bigr) \cdot J_1}{\Phi \bigl( |Du|_\gamma \bigr)\xi} = & \biggl[  h_\varepsilon \bigl( |Du|_\gamma \bigr) \gamma_{\alpha\beta}D_\nu D_\alpha u^i +h'_\varepsilon \bigl( |Du|_\gamma \bigr)|Du|_\gamma \dfrac{\gamma_{\alpha\beta} D_\alpha u^i}{|Du|_\gamma^2}F_\nu \biggr] \gamma_{\nu \sigma} D_\beta D_\sigma u^i\\
    =& \  \mathcal{A}_\varepsilon(x,Du)\bigl( D^2u, D^2u \bigr).
\end{align*}
This yields the identity
$$\bigl( I_1 + I_2 \bigr) \cdot J_1 = \mathcal{A}_\varepsilon(x,Du)\bigl( D^2u, D^2u\bigr) \Phi \bigl( |Du|_\gamma \bigr) \xi. $$
Now, we consider the term $\bigl( I_1+I_2 \bigr) \cdot J_2$. First, we compute
\begin{equation*}
    \dfrac{\bigl( I_1 \cdot J_2 \bigr) |Du|_\gamma}{\Phi' \bigl( |Du|_\gamma \bigr) \xi} = h_\varepsilon \bigl( |Du|_\gamma \bigr) \gamma_{\alpha\beta} D_\nu D_\alpha u^i \gamma_{\nu \sigma} D_\sigma u^i F_\beta = h_\varepsilon \bigl( |Du|_\gamma \bigr) \gamma_{\alpha \beta} F_\alpha F_\beta.
\end{equation*}
For the product $I_2 \cdot J_2$, we have
\begin{equation*}
    \dfrac{\bigl( I_2 \cdot J_2 \bigr) |Du|_\gamma}{\Phi' \bigl( |Du|_\gamma \bigr) \xi} = h'_\varepsilon \bigl( |Du|_\gamma\bigr) |Du|_\gamma \dfrac{\gamma_{\alpha\beta}D_\alpha u^i}{|Du|^2_\gamma} F_\nu \gamma_{\nu \sigma}D_\sigma u^i F_\beta.
\end{equation*}
Combining these identities, we find
\begin{align*}
    \bigl( I_1+I_2 \bigr) \cdot J_2 = & \ \dfrac{\Phi' \bigl( |Du|_\gamma \bigr) \xi}{|Du|_\gamma} \biggl[ h_\varepsilon \bigl( |Du|_\gamma \bigr) \gamma_{\alpha\beta}+ h'_\varepsilon \bigl( |Du|_\gamma\bigr) |Du|_\gamma \dfrac{\gamma_{\alpha\sigma} D_\sigma u^i \gamma_{\beta\delta}D_\delta u^i}{|Du|_\gamma^2}  \biggr] F_\alpha F_\beta\\
    = & \ \dfrac{\Phi' \bigl( |Du|_\gamma \bigr) \xi}{|Du|_\gamma}  \mathcal{C}_\varepsilon (x,Du)(F,F).
\end{align*}
Next, we turn to $\bigl(  I_1+I_2 \bigr) \cdot J_5$. We get
\begin{align*}
    \dfrac{\bigl(  I_1+I_2 \bigr) \cdot J_5}{\Phi \bigl( |Du|_\gamma \bigr)}
    = & \ h_\varepsilon \bigl( |Du|_\gamma \bigr) \gamma_{\alpha\beta} \bigl(  \gamma_{\nu\sigma} D_\sigma u^i D_\nu D_\alpha u^i\bigr) \dfrac{\partial \xi}{\partial x_\beta}  + h'_\varepsilon  \bigl( |Du|_\gamma\bigr)  |Du|_\gamma\dfrac{\gamma_{\alpha\beta}D_\alpha u^i}{|Du|_\gamma^2} F_\nu \gamma_{\nu \sigma} D_\sigma u^i \dfrac{\partial \xi}{\partial x_\beta}\\
    = & \ \mathcal{C}_\varepsilon (x,Du) (F,\nabla \xi).
\end{align*}
This yields
$$ \bigl(  I_1+I_2 \bigr) \cdot J_5 = \Phi \bigl( |Du|_\gamma \bigr) \mathcal{C}_\varepsilon (x,Du) (F,\nabla \xi).$$
Now, we take care of the terms $I_1 \cdot J_3$ and $I_2 \cdot J_3$. Note that
$$\dfrac{\bigl(I_1 \cdot J_3 \bigr) 2 |Du|_\gamma}{\Phi' \bigl(  |Du|_\gamma \bigr) \xi}= h_\varepsilon \bigl( |Du|_\gamma \bigr) \gamma_{\alpha\beta} D_\nu D_\alpha u^i \gamma_{\nu \sigma}D_\sigma u^i V_\beta= h_\varepsilon \bigl(  |Du|_\gamma \bigr) \gamma_{\alpha \beta} F_\alpha V_\beta.$$
Similarly, we obtain
$$\dfrac{\bigl(I_2 \cdot J_3 \bigr) 2 |Du|_\gamma}{\Phi' \bigl(  |Du|_\gamma \bigr) \xi} = h'_\varepsilon \bigl( |Du|_\gamma \bigr) |Du|_\gamma \dfrac{\gamma_{\alpha\beta}D_\alpha u^i}{|Du|_\gamma^2}F_\nu \gamma_{\nu\sigma}D_\sigma u^i V_\beta.$$
Hence, we derive that
$$\bigl( I_1+I_2 \bigr) \cdot J_3 = \dfrac{\Phi' \bigl( |Du|_\gamma \bigr) \xi}{2 |Du|_\gamma} \mathcal{C}_\varepsilon(x,Du)(F,V).$$
Using previous identities in \eqref{6.2}, we obtain
\begin{align}\label{6.3}
    0 = & \int_{B_R} \biggl[  \mathcal{A}_\varepsilon(x,Du) \bigl( D^2u, D^2u \bigr) \Phi \bigl(  |Du|_\gamma \bigr) + \mathcal{C}_\varepsilon(x,Du)(F,F) \dfrac{\Phi' \bigl( |Du|_\gamma  \bigr)}{|Du|_\gamma} \biggr] \xi  \ \mathrm{d} x \notag\\
    & + \int_{B_R} \mathcal{C}_\varepsilon(x,Du)(F,\nabla \xi) \Phi \bigl( |Du|_\gamma \bigr) \ \mathrm{d} x + \int_{B_R} \mathcal{C}_\varepsilon(x,Du)(F,V) \dfrac{\Phi' \bigl( |Du|_\gamma \bigr)}{2 |Du|_\gamma} \xi \ \mathrm{d}x \notag \\
    &+ \int_{B_R} \biggl[ \bigl(  I_1+I_2  \bigr) \cdot J_4 + \sum_{k=3}^4 \sum_{l=1}^5 I_k \cdot J_l\biggr] \ \mathrm{d}x.
\end{align}
Next, we estimate $I_1 \cdot J_4$ and $I_2 \cdot J_4$. We infer
$$\dfrac{I_1 \cdot J_4}{\Phi \bigl( |Du|_\gamma \bigr) \xi} = h_\varepsilon \bigl(  |Du|_\gamma \bigr) \gamma_{\alpha \beta} D_\nu D_\alpha u^i \dfrac{\partial \gamma_{\nu \sigma}}{\partial x_\beta} D_\sigma u^i$$
and
$$  \dfrac{I_2 \cdot J_4}{\Phi \bigl( |Du|_\gamma \bigr) \xi} = h'_\varepsilon \bigl(  |Du|_\gamma \bigr) |Du|_\gamma \dfrac{\gamma_{\alpha\beta}D_\alpha u^i}{|Du|_\gamma^2}F_\nu \dfrac{\partial \gamma_{\nu \sigma}}{\partial x_\beta} D_\sigma u^i.$$
Combining the previous identities, we get
$$\dfrac{\bigl( I_1+ I_2 \bigr) \cdot J_4}{\Phi \bigl( |Du|_\gamma \bigr)
\xi} = \biggl[ h_\varepsilon \bigl(  |Du|_\gamma \bigr) \gamma_{\alpha \beta} \delta^{ij} + h'_\varepsilon \bigl(  |Du|_\gamma \bigr) |Du|_\gamma \dfrac{\gamma_{k\beta}D_k u^i \gamma_{\alpha\lambda}D_\lambda u^i}{|Du|_\gamma^2} \biggr] D_\nu D_\alpha u^j \dfrac{\partial \gamma_{\nu \sigma}}{\partial x_\beta} D_\sigma u^i.$$
Denoting
$$\bigl(  \tau_\nu \bigr)^i_{\beta} := \dfrac{\partial \gamma_{\nu \sigma}}{\partial x_\beta} D_\sigma u^i$$
and recalling the definition of $\mathcal{B}_\varepsilon$, we rewrite the previous equality as
\begin{equation}\label{6.4}
  \dfrac{\bigl( I_1+ I_2 \bigr) \cdot J_4}{\Phi \bigl( |Du|_\gamma \bigr)
\xi} = \sum_{\nu=1}^n \mathcal{B}_\varepsilon (x,Du) \biggl( \dfrac{\partial Du}{\partial x_\nu}, \tau_\nu  \biggr) .
\end{equation}
Now, in order to estimate \eqref{6.4}, we define two matrices $\bigl( \delta_{\alpha\beta} \bigr)_{\alpha,\beta}$ and $\bigl(  \tilde{\gamma}_{\alpha\beta} \bigr)_{\alpha,\beta}$ such that
$$ \bigl( \delta_{\alpha\beta} \bigr)_{\alpha,\beta}^2 := \bigl( \gamma_{\alpha\beta} \bigr)_{\alpha,\beta} \quad \text{and} \quad \bigl( \tilde{\gamma}_{\alpha\beta} \bigr)_{\alpha,\beta}:= \bigl( \gamma_{\alpha\beta} \bigr)^{-1}_{\alpha,\beta}.$$
These matrices are well defined thanks to the properties of $\bigl( \gamma_{\alpha\beta} \bigr)_{\alpha,\beta}$. Moreover, we define $\bigl(   \tilde{\delta}_{\alpha\beta} \bigr)_{\alpha,\beta}$ so that
$$\bigl( \tilde{\delta}_{\alpha\beta} \bigr)_{\alpha,\beta}:=\bigl( \delta_{\alpha\beta} \bigr)_{\alpha,\beta}^{-1}.$$
We estimate $\mathcal{B}_\varepsilon (x,Du) \bigl( \frac{\partial Du}{\partial x_\nu}, \tau_\nu  \bigr)$ for every $\nu \in \{ 1,..., n \}$ as follows
\begin{align*}
   \mathcal{B}_\varepsilon (x,Du) \biggl( \dfrac{\partial Du}{\partial x_\nu}, \tau_\nu  \biggr) = &  \biggl[ h_\varepsilon \bigl(  |Du|_\gamma \bigr) \gamma_{\alpha \beta} \delta^{ij} + h'_\varepsilon \bigl(  |Du|_\gamma \bigr) |Du|_\gamma \dfrac{\gamma_{k\beta}D_k u^i \gamma_{\alpha\lambda}D_\lambda u^i}{|Du|_\gamma^2} \biggr] D_\nu D_\alpha u^j \bigl(  \tau_\nu \bigr)^i_\beta \\
   = & \bigl[ . . . \bigr] \delta_{\nu \lambda} \tilde{\delta}_{\lambda \sigma}  D_\nu D_\alpha u^j \bigl(  \tau_\nu \bigr)^i_\beta \\
   = & \bigl[ . . . \bigr] \bigl( \delta_{\nu \lambda}  D_\nu D_\alpha u^j \bigr) \bigl(  \tilde{\delta}_{\lambda \sigma} \bigl(  \tau_\nu \bigr)^i_\beta \bigr) \\
   = & \sum_{\lambda =1}^n \mathcal{B}_\varepsilon (x,Du) \bigl( \eta_\lambda, \tilde{\tau}_{\lambda} \bigr),
\end{align*}
where we used the notations $$ \bigl( \eta_\lambda \bigr)^j_\alpha := \delta_{\nu \lambda}  D_\nu D_\alpha u^i \quad \text{and} \quad \bigl( \tilde{\tau}_\lambda \bigr)^i_\beta := \tilde{\delta}_{\lambda \sigma} \bigl(  \tau_\nu \bigr)^i_\beta. $$
By applying Cauchy-Schwarz inequality and Young's inequality and recalling the definition of $\mathcal{A}_\varepsilon$, we arrive at
\begin{align*}
    \sum_{\lambda =1}^n \mathcal{B}_\varepsilon (x,Du) \bigl( \eta_\lambda, \tilde{\tau}_{\lambda} \bigr) \le &  \ \dfrac{\theta}{n} \sum_{\lambda =1}^n \mathcal{B}_\varepsilon (x,Du) \bigl( \eta_\lambda, \eta_{\lambda} \bigr) + \dfrac{n}{\theta} \sum_{\lambda =1}^n \mathcal{B}_\varepsilon (x,Du) \bigl( \tilde{\tau}_\lambda, \tilde{\tau}_{\lambda} \bigr)\\
    \le & \ \dfrac{\theta}{n} \bigl[... \bigr] \bigl( \eta_\lambda \bigr)^i_\alpha \bigl( \eta_\lambda \bigr)^j_\beta + \dfrac{n}{\theta} \bigl[... \bigr] \bigl( \tilde{\tau}_\lambda \bigr)^i_\alpha \bigl( \tilde{\tau}_\lambda \bigr)^j_\beta\\
    \le & \ \dfrac{\theta}{n} \bigl [... \bigr] \delta_{\nu \lambda} D_\nu D_\alpha u^i \delta_{k \lambda} D_k D_\beta u^j 
     + \dfrac{n}{\theta} \bigl[... \bigr] \delta_{\nu\lambda} \tau^i_{\alpha\nu} \delta_{k\lambda} \tau^j_{\beta k}\\
    = & \ \dfrac{\theta}{n} \bigl [... \bigr] \gamma_{\nu k} D_\nu D_\alpha u^i D_k D_\beta u^j + \dfrac{n}{\theta} \bigl[... \bigr] \gamma_{\nu k} \tau^i_{\alpha \nu} \tau^j_{\beta k}\\
    = & \ \dfrac{\theta}{n} \mathcal{A}_\varepsilon(x,Du) \bigl( D^2u,D^2u \bigr) + \dfrac{n}{\theta}\mathcal{A}_\varepsilon(x,Du) \bigl( \tau,\tau   \bigr).
\end{align*}
Hence, we eventually get
$$\bigl( I_1+I_2 \bigr) \cdot J_4\le \theta \mathcal{A}_\varepsilon(x,Du) \bigl( D^2u,D^2u \bigr) \Phi \bigl( |Du|_\gamma \bigr) \xi + \dfrac{2n}{\theta}\mathcal{A}_\varepsilon(x,Du) \bigl( \tau,\tau   \bigr)  \Phi \bigl( |Du|_\gamma \bigr) \xi,$$
where we denoted 
\begin{equation}\label{deftau}
     \tau^i_{\nu\beta} := \dfrac{\partial \gamma_{\nu \sigma}}{\partial x_\beta} D_\sigma u^i.
\end{equation}

Now, we estimate the terms $I_4 \cdot J_1$ and $I_3 \cdot J_1$. We compute
\begin{align*}
    \dfrac{I_4 \cdot J_1}{\Phi \bigl(  |Du|_\gamma \bigr) \xi} = & \ h_{\varepsilon} \bigl(  |Du|_\gamma \bigr) \dfrac{\partial \gamma_{\alpha \beta}}{\partial x_\nu} D_\alpha u^i \gamma_{\nu \sigma} D_\beta D_\sigma u^i \\
    = & \ h_{\varepsilon} \bigl(  |Du|_\gamma \bigr) \gamma_{\beta \theta} \tilde{\gamma}_{\theta \lambda} \dfrac{\partial \gamma_{\alpha \lambda}}{\partial x_\nu} D_\alpha u^i \gamma_{\nu \sigma} D_\beta D_\sigma u^i \\
    = & \ h_{\varepsilon} \bigl(  |Du|_\gamma \bigr) \gamma_{\beta \theta} \gamma_{\nu \sigma } \pi^i_{\theta \nu} D_\beta D_\sigma u^i
\end{align*}
and
\begin{align*}
   \dfrac{I_3 \cdot J_1}{\Phi \bigl(  |Du|_\gamma \bigr) \xi} = & \ 
   h'_{\varepsilon} \bigl(  |Du|_\gamma \bigr) \dfrac{\gamma_{\alpha \beta} D_\alpha u^i}{2 |Du|_\gamma^2} \dfrac{\partial \gamma_{k \lambda}}{\partial x_\nu} D_k u^j D_\lambda u^j \gamma_{\nu \sigma} D_\beta D_\sigma u^i\\
   = & \ h'_{\varepsilon} \bigl(  |Du|_\gamma \bigr) \dfrac{\gamma_{\alpha \beta} D_\alpha u^i \gamma_{\theta \rho} D_\rho u^i}{2 |Du|_\gamma^2} \tilde{\gamma}_{k \theta} \dfrac{\partial \gamma_{k \lambda}}{\partial x_\nu} D_\lambda u^j \gamma_{\nu \sigma} D_\beta D_\sigma u^i\\
   = & \ h'_{\varepsilon} \bigl(  |Du|_\gamma \bigr) \dfrac{\gamma_{\alpha \beta} D_\alpha u^i \gamma_{\theta \rho} D_\rho u^i}{2 |Du|_\gamma^2}  \gamma_{\nu \sigma} \pi^j_{\theta \nu} D_\beta D_\sigma u^i,
\end{align*}
where we denoted 
\begin{equation}\label{defpi}
    \pi^i_{\theta\nu}:= \tilde{\gamma}_{k \theta} \dfrac{\partial \gamma_{k \lambda}}{\partial x_\nu} D_\lambda u^i.
\end{equation}
Combining the preceding identities and applying Young's inequality, we deduce
\begin{align*}
    \dfrac{\bigl( I_3+ I_4 \bigr) \cdot J_1}{\Phi \bigl(  |Du|_\gamma \bigr) \xi} = & \ \dfrac{1}{2} h_\varepsilon \bigl( |Du|_\gamma \bigr) \langle \pi, D^2 u \rangle_\gamma + \dfrac{1}{2} \mathcal{A}_\varepsilon (x,Du) \bigl( \pi, D^2u \bigr)\\
    \ge &  -\beta  h_\varepsilon \bigl( |Du|_\gamma \bigr) |D^2u|^2_\gamma- \dfrac{1}{4\beta}  h_\varepsilon \bigl( |Du|_\gamma \bigr) |\pi|^2_\gamma\\
    &- \theta \mathcal{A}_\varepsilon (x,Du) \bigl( D^2u, D^2u \bigr)
    - \dfrac{1}{4\theta} \mathcal{A}_\varepsilon (x,Du) \bigl( \pi, \pi \bigr).
\end{align*}
Inserting the previous estimates in \eqref{6.3}, we get the inequality
\begin{align}\label{6.5}
     \int_{B_R} & \biggl[(1-2\theta)   \mathcal{A}_\varepsilon(x,Du) \bigl( D^2u, D^2u \bigr) \Phi \bigl(  |Du|_\gamma \bigr) + \mathcal{C}_\varepsilon(x,Du)(F,F) \dfrac{\Phi' \bigl( |Du|_\gamma  \bigr)}{|Du|_\gamma} \biggr] \xi  \ \mathrm{d} x \notag\\
    & + \int_{B_R} \mathcal{C}_\varepsilon(x,Du)(F,\nabla \xi) \Phi \bigl( |Du|_\gamma \bigr) \ \mathrm{d} x + \int_{B_R} \mathcal{C}_\varepsilon(x,Du)(F,V) \dfrac{\Phi' \bigl( |Du|_\gamma \bigr)}{2 |Du|_\gamma} \xi \ \mathrm{d}x \notag \\
    &- \beta \int_{B_R}h_\varepsilon \bigl( |Du|_\gamma \bigr) |D^2u|^2_\gamma \Phi \bigl(  |Du|_\gamma \bigr) \xi \ \mathrm{d}x
    + \int_{B_R}  \sum_{k=3}^4 \sum_{l=2}^5 I_k \cdot J_l \ \mathrm{d}x \notag  \\
    \le & \ c(n,\theta) \int_{B_R} \bigl[\mathcal{A}_\varepsilon (x,Du) \bigl( \tau, \tau \bigr) + \mathcal{A}_\varepsilon (x,Du) \bigl( \pi, \pi \bigr)   \bigr] \Phi \bigl( |Du|_\gamma \bigr) \xi \ \mathrm{d} x \notag\\
    & + \dfrac{1}{\beta} \int_{B_R} h_\varepsilon \bigl( |Du|_\gamma \bigr) |\pi|^2_\gamma \Phi \bigl(  |Du|_\gamma \bigr) \xi \ \mathrm{d}x.
\end{align}
Now, we choose $\theta =\frac{1}{4}$ and estimate the term
\begin{align*}
    (1  -2 & \theta)   \mathcal{A}_\varepsilon(x,Du) \bigl( D^2u, D^2u \bigr) - \beta h_\varepsilon \bigl( |Du|_\gamma \bigr) |D^2u|^2_\gamma \\
    = & \ \biggl( \dfrac{1}{2}-\beta \biggr)  h_\varepsilon \bigl( |Du|_\gamma \bigr) |D^2u|^2_\gamma+ \dfrac{1}{2}  h'_\varepsilon \bigl( |Du|_\gamma \bigr) |Du|_\gamma \dfrac{\gamma_{\alpha\beta}D_\alpha u^i \gamma_{k\lambda}D_k u^j}{|Du|_\gamma^2} \gamma_{\nu \sigma} D_\beta D_\nu u^i D_\lambda D_\sigma u^j.
\end{align*}
Choosing $\beta = \frac{1}{4} \min \{ 1, p-1 \}$, we have
\begin{align}\label{6.6}
    \biggl( \dfrac{1}{2}- & \beta \biggr)  h_\varepsilon \bigl( |Du|_\gamma \bigr) |D^2u|^2_\gamma+ \dfrac{1}{2}  h'_\varepsilon \bigl( |Du|_\gamma \bigr) |Du|_\gamma \dfrac{\gamma_{\alpha\beta}D_\alpha u^i \gamma_{k\lambda}D_k u^j}{|Du|_\gamma^2} \gamma_{\nu \sigma} D_\beta D_\nu u^i D_\lambda D_\sigma u^j \notag\\
    & - \dfrac{1}{4} \mathcal{A}_\varepsilon (x,Du) \bigl( D^2u, D^2u \bigr) \notag\\
    = & \ \dfrac{1}{4} \max \{0,2-p \}  h_\varepsilon \bigl( |Du|_\gamma \bigr) |D^2u|^2_\gamma \notag\\
    &+  \dfrac{1}{4} h'_\varepsilon \bigl( |Du|_\gamma \bigr) |Du|_\gamma \dfrac{\gamma_{\alpha\beta}D_\alpha u^i \gamma_{k\lambda}D_k u^j}{|Du|_\gamma^2} \gamma_{\nu \sigma} D_\beta D_\nu u^i D_\lambda D_\sigma u^j =: D_1+D_2 \ge 0.
\end{align}
Indeed, if $h'_\varepsilon \bigl( |Du|_\gamma \bigr)\ge 0$, the inequality follows trivially by \eqref{CSineq}.  Otherwise, if $h'_\varepsilon \bigl( |Du|_\gamma \bigr)\ < 0$, which can only happen if $p<2$ and $|Du|_\gamma >1$, the result follows again by an application of \eqref{CSineq}. Indeed, we have
\begin{align*}
    D_1 + D_2 \ge & \ \dfrac{1}{4} \biggl[ \max \{ 0, 2-p \}  h_\varepsilon \bigl( |Du|_\gamma \bigr)  +  h'_\varepsilon \bigl( |Du|_\gamma \bigr) |Du|_\gamma    \biggr] |D^2u|^2_\gamma \\
    = & \  \dfrac{1}{4} \biggl[ (2-p ) \dfrac{\bigl( |Du|_\gamma-1 \bigr)^{p-1}_+}{|Du|_\gamma}  + \dfrac{\bigl( |Du|_\gamma -1)^{p-2}_+  \bigl[  (p-2)|Du|_\gamma+1\bigr]}{|Du|_\gamma}    \biggr] |D^2u|^2_\gamma\\
    = & \ \dfrac{\bigl( |Du|_\gamma-1 \bigr)^{p-2}_+}{4|Du|_\gamma}
    \bigl[ (2-p) \bigl( |Du|_\gamma-\ \bigr) +(p-2)|Du|_\gamma+1  \bigr] |D^2u|^2_\gamma\\
    = & \ \dfrac{\bigl( |Du|_\gamma-1 \bigr)^{p-2}_+}{4|Du|_\gamma}
     |D^2u|^2_\gamma \ge 0
\end{align*}
Next, we compute the terms $I_3 \cdot J_2$ and $I_4 \cdot J_2$. We have
$$\dfrac{ \bigl(  I_3 \cdot J_2 \bigr) |Du|_\gamma}{\Phi' \bigl |Du|_\gamma \bigr) \xi} = h'_\varepsilon \bigl(  |Du|_\gamma \bigr) |Du|_\gamma \dfrac{\gamma_{\alpha\beta}D_\alpha u^i \gamma_{\nu\sigma}D_\sigma u^i}{2|Du|^2_\gamma} V_\nu F_\beta$$
and
$$ \dfrac{ \bigl(  I_3 \cdot J_2 \bigr) |Du|_\gamma}{\Phi' \bigl |Du|_\gamma \bigr) \xi} =  h_\varepsilon \bigl(  |Du|_\gamma \bigr) 
\gamma_{\beta \theta} \biggl(  \tilde{\gamma}_{\theta\lambda} \gamma_{\nu \sigma} D_\sigma u^i \dfrac{\partial \gamma_{\alpha \lambda}}{\partial x_\nu} D_\alpha u^i \biggr) F_\beta= h_\varepsilon \bigl(  |Du|_\gamma \bigr) 
\gamma_{\beta \theta} \omega_\theta F_\beta, $$
where we set 
\begin{equation}\label{defomega}
    \omega_\theta:= \tilde{\gamma}_{\theta\lambda} \gamma_{\nu \sigma} D_\sigma u^i \dfrac{\partial \gamma_{\alpha \lambda}}{\partial x_\nu} D_\alpha u^i .
\end{equation}
Combining the preceding identities, we get
$$  \bigl(  I_3 +I_4 \bigl) \ \cdot \ J_2 = \biggl[  h_\varepsilon \bigl(  |Du|_\gamma \bigr) \langle \omega ,F \rangle_\gamma + h'_\varepsilon \bigl(  |Du|_\gamma \bigr) |Du|_\gamma \dfrac{\gamma_{\alpha\beta}D_\alpha u^i \gamma_{\nu\sigma}D_\sigma u^i}{2|Du|^2_\gamma} V_\nu F_\beta   \biggr] \dfrac{\Phi' \bigl( |Du|_\gamma \bigr) \xi}{|Du|_\gamma}.$$
Now, using Young's inequality and recalling the definition of $\mathcal{C}_\varepsilon$, we derive the following estimate
\begin{align*}
     h_\varepsilon \bigl(  & |Du|_\gamma \bigr)   \langle \omega ,F \rangle_\gamma + h'_\varepsilon \bigl(  |Du|_\gamma \bigr) |Du|_\gamma \dfrac{\gamma_{\alpha\beta}D_\alpha u^i \gamma_{\nu\sigma}D_\sigma u^i}{2|Du|^2_\gamma} V_\nu F_\beta \\
     & \ \ \ + \mathcal{C}_\varepsilon(x,Du)(F,F) + \dfrac{1}{2} \mathcal{C}_\varepsilon(x,Du)(F,V)\\
     = & \ h_\varepsilon \bigl(   |Du|_\gamma \bigr)   \langle \omega ,F \rangle_\gamma  -\dfrac{1}{2} h_\varepsilon \bigl(   |Du|_\gamma \bigr)   \langle F ,V \rangle_\gamma 
      + \mathcal{C}_\varepsilon(x,Du)(F,V)+ \mathcal{C}_\varepsilon(x,Du)(F,F)\\
      \ge & - \eta h_\varepsilon \bigl(   |Du|_\gamma \bigr) |F|_\gamma^2-  \dfrac{1}{\eta}  h_\varepsilon \bigl(   |Du|_\gamma \bigr) |\omega|^2_\gamma - \dfrac{1}{2 \eta} h_\varepsilon \bigl(   |Du|_\gamma \bigr)  |V|^2_\gamma\\
      &+\dfrac{1}{2} \mathcal{C}_\varepsilon(x,Du)(F,F)- \dfrac{1}{2} \mathcal{C}_\varepsilon(x,Du)(V,V)\\
      = & \  \biggl(  \dfrac{1}{2} - \eta \biggr) h_\varepsilon \bigl(   |Du|_\gamma \bigr) |F|_\gamma^2 +\dfrac{1}{2}  h'_\varepsilon \bigl(  |Du|_\gamma \bigr) |Du|_\gamma \dfrac{\gamma_{\alpha\beta}D_\alpha u^i \gamma_{k\lambda}D_k u^i}{|Du|^2_\gamma} F_\beta F_\lambda\\
      & - \dfrac{1}{2\eta}  h_\varepsilon \bigl(   |Du|_\gamma \bigr)  |V|^2_\gamma -  \dfrac{1}{\eta}  h_\varepsilon \bigl(   |Du|_\gamma \bigr) |\omega|^2_\gamma - \dfrac{1}{2} \mathcal{C}_\varepsilon(x,Du)(V,V).
\end{align*}
Choosing $\eta:= \frac{1}{4} \min \{  1, p-1 \}$ and arguing as above, we infer
\begin{align*}
    \biggl(  \dfrac{1}{2} - \eta \biggr) h_\varepsilon \bigl(   |Du|_\gamma \bigr) |F|_\gamma^2 +\dfrac{1}{2}  h'_\varepsilon \bigl(  |Du|_\gamma \bigr) |Du|_\gamma \dfrac{\gamma_{\alpha\beta}D_\alpha u^i \gamma_{k\lambda}D_k u^i}{|Du|^2_\gamma} F_\beta F_\lambda - \dfrac{1}{4} \mathcal{C}_\varepsilon(x,Du)(F,F) \ge 0.
\end{align*}
Hence, inserting \eqref{6.6} and the previous estimates in \eqref{6.5} we get
\begin{align}\label{6.7}
     \dfrac{1}{4} \int_{B_R} & \biggl[   \mathcal{A}_\varepsilon(x,Du) \bigl( D^2u, D^2u \bigr) \Phi \bigl(  |Du|_\gamma \bigr) + \mathcal{C}_\varepsilon(x,Du)(F,F) \dfrac{\Phi' \bigl( |Du|_\gamma  \bigr)}{|Du|_\gamma} \biggr] \xi  \ \mathrm{d} x \notag\\
    & + \int_{B_R} \mathcal{C}_\varepsilon(x,Du)(F,\nabla \xi) \Phi \bigl( |Du|_\gamma \bigr) \ \mathrm{d} x 
    + \int_{B_R}  \sum_{k=3}^4 \sum_{l=2}^5 I_k \cdot J_l \ \mathrm{d}x \notag  \\
    \le & \ c(n) \int_{B_R} \bigl[\mathcal{A}_\varepsilon (x,Du) \bigl( \tau, \tau \bigr) + \mathcal{A}_\varepsilon (x,Du) \bigl( \pi, \pi \bigr)   \bigr] \Phi \bigl( |Du|_\gamma \bigr) \xi \ \mathrm{d} x \notag\\
    & + c \int_{B_R} h_\varepsilon \bigl( |Du|_\gamma \bigr) |\pi|^2_\gamma \Phi \bigl(  |Du|_\gamma \bigr) \xi \ \mathrm{d}x \notag\\
    &+ c \int_{B_R} \biggl[   h_\varepsilon \bigl(   |Du|_\gamma \bigr)  |V|^2_\gamma+  h_\varepsilon \bigl(   |Du|_\gamma \bigr) |\omega|^2_\gamma +\mathcal{C}_\varepsilon(x,Du)(V,V)  \biggl]  \dfrac{\Phi' \bigl( |Du|_\gamma \bigr) \xi}{|Du|_\gamma}.
\end{align}
Now, we consider the terms $I_3 \cdot J_3$ and $I_4 \cdot J_3$. Recalling the definition of $\omega$ in \eqref{defomega}, we derive
\begin{align*}
    \dfrac{\bigl( I_4 \cdot J_3 \bigr) |Du|_\gamma}{\Phi' \bigl( |Du|_\gamma \bigr) \xi} = & \ \dfrac{1}{2} h_\varepsilon \bigl( |Du|_\gamma \bigr) \dfrac{\partial \gamma_{\alpha\beta}}{\partial x_\nu} D_\alpha u^i \gamma_{\nu \sigma} D_\sigma u^i V_\beta\\
    = & \  \dfrac{1}{2} h_\varepsilon \bigl( |Du|_\gamma \bigr) \gamma_{\beta \lambda} \biggl( \tilde{\gamma}_{\lambda\rho}\dfrac{\partial \gamma_{\alpha\rho}}{\partial x_\nu} D_\alpha u^i \gamma_{\nu \sigma} D_\sigma u^i \biggr) V_\beta= \dfrac{1}{2} h_\varepsilon \bigl( |Du|_\gamma \bigr) \langle \omega, V \rangle_\gamma
\end{align*}
and
\begin{align*}
    \dfrac{\bigl( I_3 \cdot J_3 \bigr)  |Du|_\gamma}{\Phi' \bigl( |Du|_\gamma \bigr) \xi} = & \ h'_\varepsilon \bigl( |Du|_\gamma \bigr) |Du|_\gamma\dfrac{\gamma_{\alpha\beta}D_\alpha u^i \gamma_{\nu\sigma}D_\sigma u^i}{4 |Du|_\gamma^2} V_\nu V_\beta. 
\end{align*}
Combining the preceding identities, applying Young's inequality and recalling the definition of $\mathcal{C}_\varepsilon$, we obtain
\begin{align}\label{6.8}
    \dfrac{\bigl( I_3 +I_4 \bigr)  \cdot  J_3   |Du|_\gamma}{\Phi' \bigl( |Du|_\gamma \bigr) \xi} 
    =&  \ \dfrac{1}{4} \mathcal{C}_\varepsilon (x,Du)(V,V)  -\dfrac{1}{4} h_\varepsilon \bigl( |Du|_\gamma \bigr) |V|^2_\gamma  + \dfrac{1}{2} h_\varepsilon \bigl( |Du|_\gamma \bigr) \langle \omega, V \rangle_\gamma \notag \\
    \ge & \ \dfrac{1}{4} \mathcal{C}_\varepsilon (x,Du)(V,V)  -c h_\varepsilon \bigl( |Du|_\gamma \bigr) \bigl(|V|^2_\gamma +|\omega|^2_\gamma \bigr).
\end{align}
Next, we take care of the terms $I_3 \cdot J_4$ and $I_4 \cdot J_4$. Recalling the definition of $\tau$ in \eqref{deftau} and $\pi$ in \eqref{defpi}, we have 
\begin{align*}
    \dfrac{I_3 \cdot J_4}{\Phi \bigl(  |Du|_\gamma \bigr) \xi} = & \ h'_\varepsilon \bigl(  |Du|_\gamma \bigr) |Du|_\gamma \dfrac{\gamma_{\alpha \beta}D_\alpha u^i}{2 |Du|_\gamma^2} \dfrac{\partial \gamma_{k \lambda}}{\partial x_\nu} D_k u^j D_\lambda u^j \dfrac{\partial \gamma_{\nu \sigma}}{\partial x_\beta} D_\sigma u^i\\
    = & \  h'_\varepsilon \bigl(  |Du|_\gamma \bigr) |Du|_\gamma \dfrac{\gamma_{\alpha \beta}D_\alpha u^i \gamma_{\theta \rho}D_\rho u^j}{2 |Du|_\gamma^2} \biggl(  \tilde{\gamma}_{\lambda \theta} \dfrac{\partial \gamma_{k \lambda}}{\partial x_\nu} D_k u^j  \biggr) \biggl(\dfrac{\partial \gamma_{\nu \sigma}}{\partial x_\beta} D_\sigma u^i \biggr)\\
    = & \ h'_\varepsilon \bigl(  |Du|_\gamma \bigr) |Du|_\gamma \dfrac{\gamma_{\alpha \beta}D_\alpha u^i \gamma_{\theta \rho}D_\rho u^j}{2 |Du|_\gamma^2} \pi^j_{\theta \nu} \tau^i_{\nu \beta}
\end{align*}
and
\begin{align*}
    \dfrac{I_4 \cdot J_4}{\Phi \bigl(  |Du|_\gamma \bigr) \xi} = & \
    h_\varepsilon \bigl( |Du|_\gamma \bigr) \dfrac{\partial \gamma_{\alpha\beta}}{\partial x_\nu} D_\alpha u^i \dfrac{\partial \gamma_{\nu \sigma}}{\partial x_\beta} D_\sigma u^i\\
    = & \ h_\varepsilon \bigl( |Du|_\gamma \bigr) \gamma_{\beta \theta} \biggl( \tilde{\gamma}_{\theta\rho} \dfrac{\partial \gamma_{\alpha\rho}}{\partial x_\nu} D_\alpha u^i  \biggr) \dfrac{\partial \gamma_{\nu \sigma}}{\partial x_\beta} D_\sigma u^i \\
    =& \ h_\varepsilon \bigl( |Du|_\gamma \bigr) \gamma_{\beta \theta} \pi^i_{\theta \nu} \tau^i_{\beta \nu}.
\end{align*}
Combining the previous identities, applying Young's inequality and recalling the definition of $\mathcal{B}_\varepsilon$, we get
\begin{align}\label{6.9}
\dfrac{I_4 \cdot J_4}{\Phi \bigl(  |Du|_\gamma \bigr) \xi} =& \
\dfrac{1}{2}\sum_{\nu=1}^n \bigl[  \mathcal{B}_\varepsilon (x,Du) \bigl( \pi_\nu, \tau_\nu \bigr) +  h_\varepsilon \bigl( |Du|_\gamma \bigr)  \langle, \pi_\nu, \tau_\nu \rangle_\gamma\bigr]\notag \\
\ge & - c  \sum_{\nu=1}^n \bigl[  \mathcal{B}_\varepsilon (x,Du) \bigl( \tau_\nu, \tau_\nu \bigr) + \mathcal{B}_\varepsilon (x,Du) \bigl( \pi_\nu, \pi_\nu \bigr)  \bigr]  \notag\\
& -c h_\varepsilon \bigl( |Du|_\gamma \bigr) \bigl( |\pi|^2_\gamma+|\tau|_\gamma^2  \bigr).
\end{align}
Now, we consider the remaining terms $I_3 \cdot J_5$ and $I_4 \cdot J_5$. Recalling the definition of $\omega$ in \eqref{defomega}, we compute
\begin{align*}
    \dfrac{I_3 \cdot J_5}{\Phi \bigl(  |Du|_\gamma \bigr) } =&  \ h'_\varepsilon  \bigl(  |Du|_\gamma \bigr)  |Du|_\gamma\dfrac{\gamma_{\alpha \beta}D_\alpha u^i}{2 |Du|_\gamma^2}
    \dfrac{\partial \gamma_{k\lambda}}{\partial x_\nu}D_k u^j
D_\lambda u^j \gamma_{\nu \sigma}D_\sigma u^i
\dfrac{\partial \xi}{\partial x_\beta}\\
=& \ h'_\varepsilon  \bigl(  |Du|_\gamma \bigr)  |Du|_\gamma\dfrac{\gamma_{\alpha \beta}D_\alpha u^i \gamma_{\nu \sigma }D_\sigma u^i}{2 |Du|_\gamma^2}
    V_\nu
\dfrac{\partial \xi}{\partial x_\beta}
\end{align*}
and
\begin{align*}
   \dfrac{I_4 \cdot J_5}{\Phi \bigl(  |Du|_\gamma \bigr) } =& \ 
   h_\varepsilon  \bigl(  |Du|_\gamma \bigr)  \gamma_{\beta\theta} \biggl(   \tilde{\gamma}_{\theta \lambda} D_\sigma u^i \dfrac{\partial \gamma_{\alpha\lambda}}{\partial x_nu} D_\alpha u^i\biggr) \dfrac{\partial \xi}{\partial x_\beta}\\
   =& \ h_\varepsilon  \bigl(  |Du|_\gamma \bigr)  \gamma_{\beta\theta} \omega_\theta \dfrac{\partial \xi}{\partial x_\beta}.
\end{align*}
Combining the previous identities and recalling the definition of $\mathcal{C}_\varepsilon$, we obtain
\begin{align}\label{6.10}
    \dfrac{\bigl( I_3 +I_4 \bigr) \cdot J_5}{\Phi \bigl(  |Du|_\gamma \bigr) } =&  \ h_\varepsilon  \bigl(  |Du|_\gamma \bigr) \langle \omega, \nabla \xi \rangle_\gamma + h'_\varepsilon  \bigl(  |Du|_\gamma \bigr)  |Du|_\gamma\dfrac{\gamma_{\alpha \beta}D_\alpha u^i \gamma_{\nu \sigma }D_\sigma u^i}{2 |Du|_\gamma^2}
    V_\nu
\dfrac{\partial \xi}{\partial x_\beta}\notag\\
=& \ \dfrac{1}{2} \mathcal{C}_\varepsilon(x,Du) ( V ,\nabla \xi ) - \dfrac{1}{2}  h_\varepsilon  \bigl(  |Du|_\gamma \bigr) \langle V, \nabla \xi \rangle_\gamma +  h_\varepsilon  \bigl(  |Du|_\gamma \bigr) \langle \omega, \nabla \xi \rangle_\gamma \notag\\
\ge & \ \dfrac{1}{2} \mathcal{C}_\varepsilon(x,Du) ( V ,\nabla \xi ) - \dfrac{1}{2}  h_\varepsilon  \bigl(  |Du|_\gamma \bigr) |V|_\gamma |\nabla \xi |_\gamma  -   h_\varepsilon  \bigl(  |Du|_\gamma \bigr) |\omega|_\gamma |\nabla \xi |_\gamma.
\end{align}
Inserting estimates \eqref{6.8}, \eqref{6.9} and \eqref{6.10} in \eqref{6.7}, we derive
\begin{align}
     \dfrac{1}{4} \int_{B_R} & \biggl[   \mathcal{A}_\varepsilon(x,Du) \bigl( D^2u, D^2u \bigr) \Phi \bigl(  |Du|_\gamma \bigr) + \mathcal{C}_\varepsilon(x,Du)(F,F) \dfrac{\Phi' \bigl( |Du|_\gamma  \bigr)}{|Du|_\gamma} \biggr] \xi  \ \mathrm{d} x \notag\\
    & + \int_{B_R} \mathcal{C}_\varepsilon(x,Du)(F,\nabla \xi) \Phi \bigl( |Du|_\gamma \bigr) \ \mathrm{d} x 
    \notag  \\
    \le & \ c\sum_{\nu=1}^n \int_{B_R} \bigl[  \mathcal{B}_\varepsilon (x,Du) \bigl( \tau_\nu, \tau_\nu \bigr) + \mathcal{B}_\varepsilon (x,Du) \bigl( \pi_\nu, \pi_\nu \bigr)  \bigr] \Phi \bigl(  |Du|_\gamma \bigr) \xi \ \mathrm{d} x \notag\\
    & +c(n) \int_{B_R} \bigl[\mathcal{A}_\varepsilon (x,Du) \bigl( \tau, \tau \bigr) + \mathcal{A}_\varepsilon (x,Du) \bigl( \pi, \pi \bigr)   \bigr] \Phi \bigl( |Du|_\gamma \bigr) \xi \ \mathrm{d} x \notag\\
    & + \dfrac{1}{4} \int_{B_R } |\mathcal{C}_\varepsilon(x,Du) ( V ,\nabla \xi )| \Phi \bigl(  |Du|_\gamma \bigr) \ \mathrm{d}x 
    + \dfrac{1}{2} \int_{B_R} h_\varepsilon  \bigl(  |Du|_\gamma \bigr) |V|_\gamma |\nabla \xi |_\gamma \Phi \bigl(  |Du|_\gamma \bigr) \ \mathrm{d}x \notag\\
    &+ \int_{B_R}  h_\varepsilon   \bigl(  |Du|_\gamma \bigr) |\omega|_\gamma |\nabla \xi |_\gamma \Phi \bigl(  |Du|_\gamma \bigr) \ \mathrm{d}x \notag + c \int_{B_R} h_\varepsilon \bigl( |Du|_\gamma \bigr) \bigl( |\pi|^2_\gamma+|\tau|_\gamma^2  \bigr) \Phi \bigl(  |Du|_\gamma \bigr) \xi \ \mathrm{d}x \notag\\
    &+ c \int_{B_R} \biggl[   h_\varepsilon \bigl(   |Du|_\gamma \bigr)  |V|^2_\gamma+  h_\varepsilon \bigl(   |Du|_\gamma \bigr) |\omega|^2_\gamma +\mathcal{C}_\varepsilon(x,Du)(V,V)  \biggl]  \dfrac{\Phi' \bigl( |Du|_\gamma \bigr) \xi}{|Du|_\gamma}.
\end{align}
Observe that
\begin{align}\label{gammanorm}
\begin{cases}
|\tau|_\gamma \le c(C_2)|Du|_\gamma \\
|\pi|_\gamma\le c(C_0,C_1,C_2)|Du|_\gamma   \\
|\omega|_\gamma \le c(C_0,C_1,C_2)|Du|_\gamma^2  \\
|V|_\gamma \le c(C_2)|Du|_\gamma^2 
\end{cases}
\end{align}
Using \eqref{gammanorm}, we have proved the following
\begin{lem}\label{mainintineq}
Let $\varepsilon \in (0,1]$ and $u:=u_\varepsilon $ be the unique weak solution of \eqref{regsystem}. Then, for any increasing  function $ \Phi \in W^{1,\infty}_{loc} \bigl ( \mathbb{R}_{\ge 0}, \mathbb{R}_{\ge 0} \bigr) $ and any non-negative testing function $\xi \in \mathcal{C}^1_0(B_R)$, we have
\begin{align*}
     \dfrac{1}{4} \int_{B_R} & \biggl[   \mathcal{A}_\varepsilon(x,Du) \bigl( D^2u, D^2u \bigr) \Phi \bigl(  |Du|_\gamma \bigr) + \mathcal{C}_\varepsilon(x,Du)(F,F) \dfrac{\Phi' \bigl( |Du|_\gamma  \bigr)}{|Du|_\gamma} \biggr] \xi  \ \mathrm{d} x \notag\\
    & + \int_{B_R} \mathcal{C}_\varepsilon(x,Du)(F,\nabla \xi) \Phi \bigl( |Du|_\gamma \bigr) \ \mathrm{d} x 
    \notag  \\
    \le &  \ c\sum_{\nu=1}^n \int_{B_R} \bigl[  \mathcal{B}_\varepsilon (x,Du) \bigl( \tau_\nu, \tau_\nu \bigr) + \mathcal{B}_\varepsilon (x,Du) \bigl( \pi_\nu, \pi_\nu \bigr)  \bigr] \Phi \bigl(  |Du|_\gamma \bigr) \xi \ \mathrm{d} x \notag\\
    & + c \int_{B_R} \bigl[\mathcal{A}_\varepsilon (x,Du) \bigl( \tau, \tau \bigr) + \mathcal{A}_\varepsilon (x,Du) \bigl( \pi, \pi \bigr)   \bigr] \Phi \bigl( |Du|_\gamma \bigr) \xi \ \mathrm{d} x \notag\\
    & + \dfrac{1}{4} \int_{B_R } |\mathcal{C}_\varepsilon(x,Du) ( V ,\nabla \xi )| \Phi \bigl(  |Du|_\gamma \bigr) \ \mathrm{d}x \notag\\
    &+c\int_{B_R}  h_\varepsilon   \bigl(  |Du|_\gamma \bigr) |Du|^2_\gamma \bigl(|\nabla \xi | + \xi \bigr) \Phi \bigl(  |Du|_\gamma \bigr) \ \mathrm{d}x  \notag\\
    &+ c \int_{B_R} \bigl[   h_\varepsilon \bigl(   |Du|_\gamma \bigr)  |Du|^4_\gamma + \mathcal{C}_\varepsilon(x,Du)(V,V)  \bigl]  \dfrac{\Phi' \bigl( |Du|_\gamma \bigr) \xi}{|Du|_\gamma}.
\end{align*}
    for a constant $c:=c(n,p,C_0,C_1,C_2)$, where $\pi^i_{\theta\nu}:= \tilde{\gamma}_{k \theta} \dfrac{\partial \gamma_{k \lambda}}{\partial x_\nu} D_\lambda u^i$, $\bigl( \tilde{\gamma}_{\alpha\beta} \bigr)_{\alpha,\beta}:= \bigl( \gamma_{\alpha\beta} \bigr)^{-1}_{\alpha,\beta} $, $\tau^i_{\nu\beta} := \dfrac{\partial \gamma_{\nu \sigma}}{\partial x_\beta} D_\sigma u^i$, $F_\nu := \gamma_{k\lambda} D_k u^j  D_\nu D_\lambda u^j$ and $V_\beta := \dfrac{\partial \gamma_{k \lambda}}{\partial x_\beta} D_k u^j D_\lambda u^j$.
\end{lem}

\subsection{Sub-solution to an elliptic equation}
We start by showing that $U_\varepsilon:= \bigl(  |Du_\varepsilon|_\gamma-1-\delta \bigr)^2_+$ is a sub-solution of a certain elliptic equation. More precisely, we prove
\begin{lem}\label{elleq}
    Let $\varepsilon \in (0,1]$ and $u_\varepsilon \in W^{1,\mathbf{p}}\bigl( B_R, \mathbb{R}^N \bigr)$ be the weak solution of the regularized system \eqref{regsystem} satisfying \eqref{bounded} and \eqref{delta} on $B_\rho(x_0) \subset B_{r_0}\Subset B_R$. Then, the function 
    \begin{equation}\label{Ue}
        U_\varepsilon:= \bigl(  |Du_\varepsilon|_\gamma-1-\delta \bigr)^2_+
    \end{equation}
satisfies the integral inequality
\begin{equation}\label{intellin}
    \int_{B_\rho(x_0)} A_{\alpha \beta}(x) D_\alpha U_\varepsilon D_\beta \xi \ \mathrm{d}x \le c \int_{B_\rho(x_0)} \bigl(  \xi + |\nabla \xi| \bigr) \ \mathrm{d}x
\end{equation}
for any non-negative testing function $\xi \in \mathcal{C}^1_0(B_R)$ and for a constant $c:=c(M,\delta,p,C_0,C_1,C_2)$. Here, the coefficients $A_{\alpha \beta}$ are defined by
\begin{equation}\label{coeffA}
    A_{\alpha \beta}(x) \eta_\alpha \zeta_\beta:= |Du_\varepsilon(x)|_{\gamma(x)}\mathcal{C}_\varepsilon(x,Du_\varepsilon(x)) (\eta, \zeta ), \quad \quad \forall \eta, \zeta \in \mathbb{R}^n, \  x \in B_R.
\end{equation}
\end{lem}
\proof We apply Lemma \ref{mainintineq} with $\Phi(t):= ( t-1-\delta )_+$. Due to Lemma \ref{ellipticity}, the first two integrals on the left-hand side are non-negative and therefore can be discarded. Moreover, using once again Lemma \ref{ellipticity} and \eqref{gammanorm} in order to bound the right-hand side, we obtain
\begin{align}\label{6.14}
    L:= & \int_{B_{\rho}(x_0)} \mathcal{C}_\varepsilon(x,Du_\varepsilon)(F,\nabla \xi) \Phi \bigl( |Du_\varepsilon|_\gamma \bigr) \ \mathrm{d} x 
    \notag  \\
    \le & \  c \int_{B_{\rho}(x_0)} \bigl[  \varepsilon + \Lambda \bigl( |Du_\varepsilon|_\gamma \bigr) \bigr] |Du_\varepsilon|_\gamma^3 \Phi' \bigl( |Du_\varepsilon|_\gamma \bigr) \xi \ \mathrm{d} x \notag\\
    &+ c \int_{B_{\rho}(x_0)} \bigl[  \varepsilon + \Lambda \bigl( |Du_\varepsilon|_\gamma \bigr)  \bigl] |Du_\varepsilon|_\gamma^2 \bigl( \xi +|\nabla \xi| \bigr)\Phi \bigl( |Du_\varepsilon|_\gamma \bigr) \ \mathrm{d}x ,
\end{align}
for any non-negative function $\xi \in \mathcal{C}^1_0(B_{\rho}(x_0))$ and for a constant $c:=c(p,C_0,C_1,C_2)$. Observe that for every $\nu \in \{ 1,...,n \}$
$$ \dfrac{\partial U_\varepsilon}{\partial x_\nu} = \dfrac{ \Phi \bigl( |Du_\varepsilon |_\gamma \bigr)}{|Du_\varepsilon|_\gamma} \bigl( V_\nu + 2 F_\nu \bigr) .$$
Hence, we have
$$ \nabla U_\varepsilon = \dfrac{ \Phi \bigl( |Du_\varepsilon |_\gamma \bigr)}{|Du_\varepsilon|_\gamma} \bigl( V + 2 F \bigr) .$$
By the linearity of $\mathcal{C}_\varepsilon(x,Du_\varepsilon)$ with respect to the first variable, we get
$$\mathcal{C}_\varepsilon(x,Du_\varepsilon) \bigl( \nabla U_\varepsilon, \nabla \xi  \bigr)= \dfrac{ \Phi \bigl( |Du_\varepsilon |_\gamma \bigr)}{|Du_\varepsilon|_\gamma} \bigl[  \mathcal{C}_\varepsilon(x,Du_\varepsilon) \bigl( V, \nabla \xi  \bigr) + 2 \mathcal{C}_\varepsilon(x,Du_\varepsilon) \bigl( F, \nabla \xi  \bigr) \bigr].$$
Therefore, we can write the left-hand side in \eqref{6.14} as follows
$$L= \dfrac{1}{2} \int_{B_{\rho}(x_0)} \mathcal{C}_\varepsilon(x,Du_\varepsilon) \bigl(\nabla U_\varepsilon,\nabla \xi \bigr )   |Du_\varepsilon|_\gamma  \ \mathrm{d} x - \dfrac{1}{2} \int_{B_{\rho}(x_0)} \mathcal{C}_\varepsilon(x,Du_\varepsilon)(V,\nabla \xi) \Phi \bigl( |Du_\varepsilon|_\gamma \bigr) \ \mathrm{d} x.  $$
Moreover, by the definition of the coefficients $A_{\alpha\beta}$, it holds
$$  \dfrac{1}{2} \int_{B_{\rho}(x_0)} \mathcal{C}_\varepsilon(x,Du_\varepsilon) \bigl(\nabla U_\varepsilon,\nabla \xi \bigr )   |Du_\varepsilon|_\gamma  \ \mathrm{d} x = \dfrac{1}{2}  \int_{B_\rho(x_0)} A_{\alpha \beta}(x) D_\alpha U_\varepsilon D_\beta \xi \ \mathrm{d}x $$
and using Cauchy-Schwarz inequality and Lemma \ref{ellipticity} we derive the estimate
$$ \mathcal{C}_\varepsilon(x,Du_\varepsilon)(V,\nabla \xi) \Phi \bigl( |Du_\varepsilon|_\gamma \bigr) \le c \bigl[  \varepsilon+ \Lambda \bigl(  |Du_\varepsilon|_\gamma \bigr) \bigr] |Du_\varepsilon|_\gamma^2 |\nabla \xi | \Phi \bigl(   |Du_\varepsilon|_\gamma \bigr).$$
Note that due to assumptions \eqref{bounded} and \eqref{delta} we have
$$ |Du_\varepsilon|_\gamma \le M \quad \text{and} \quad \Phi \bigl(  |Du_\varepsilon|_\gamma \bigr) \le \mu.$$
Moreover, $\Phi (t)=0$ for every $t \in [0,1+\delta]$. Therefore, in the case $|Du_\varepsilon|_\gamma \le 1+\delta$, it holds
$$ \bigl[  \varepsilon + \Lambda \bigl( |Du_\varepsilon|_\gamma \bigr)  \bigl] |Du_\varepsilon|_\gamma \Phi \bigl( |Du_\varepsilon|_\gamma \bigr) = 0 = \bigl[  \varepsilon + \Lambda \bigl( |Du_\varepsilon|_\gamma \bigr)  \bigl] |Du_\varepsilon|_\gamma \Phi' \bigl( |Du_\varepsilon|_\gamma \bigr).   $$
On the other hand, if $|Du_\varepsilon|_\gamma > 1+\delta$, we infer
\begin{align*}
  \bigl[  \varepsilon + \Lambda  & \bigl( |Du_\varepsilon|_\gamma \bigr)  \bigl] |Du_\varepsilon|_\gamma \Phi \bigl( |Du_\varepsilon|_\gamma \bigr)   \\
  \le & \ \bigl[ \varepsilon|Du_\varepsilon|_\gamma + \max \bigl\{  \bigl( |Du_\varepsilon|_\gamma-1 \bigr)^{p-1}, (p-1) |Du_\varepsilon|_\gamma \bigl( |Du_\varepsilon|_\gamma-1 \bigr)^{p-2}  \bigr\}   \bigr] \mu\\
  \le & \ \biggl[  M +\dfrac{pM^p}{\delta}\biggr] \mu =: c(p,M,\delta).
\end{align*}
Similarly, we get
\begin{align*}
     \bigl[  \varepsilon + \Lambda  & \bigl( |Du_\varepsilon|_\gamma \bigr)  \bigl] |Du_\varepsilon|_\gamma \Phi' \bigl( |Du_\varepsilon|_\gamma \bigr) \le c(p,M,\delta).
\end{align*}
So, we eventually arrive at
\begin{align*}
   \dfrac{1}{2}  \int_{B_\rho(x_0)} & A_{\alpha \beta}(x)  D_\alpha U_\varepsilon D_\beta \xi \ \mathrm{d}x  \\
   = & \ \dfrac{1}{2} \int_{B_{\rho}(x_0)} \mathcal{C}_\varepsilon(x,Du_\varepsilon) \bigl(\nabla U_\varepsilon,\nabla \xi \bigr )   |Du_\varepsilon|_\gamma  \ \mathrm{d} x  \\
   \le &  \ c \int_{B_{\rho}(x_0)} \bigl[  \varepsilon + \Lambda \bigl( |Du_\varepsilon|_\gamma \bigr) \bigr] |Du_\varepsilon|_\gamma^3 \Phi' \bigl( |Du_\varepsilon|_\gamma \bigr) \xi \ \mathrm{d} x \notag\\
    & + c \int_{B_{\rho}(x_0)} \bigl[  \varepsilon + \Lambda \bigl( |Du_\varepsilon|_\gamma \bigr)  \bigl] |Du_\varepsilon|_\gamma^2 \bigl( \xi +|\nabla \xi| \bigr)\Phi \bigl( |Du_\varepsilon|_\gamma \bigr) \ \mathrm{d}x \\
    \le &  \ c \int_{B_{\rho}(x_0)} \bigl( \xi + |\nabla \xi |\bigr) \ \mathrm{d}x,
\end{align*}
with a constant $c:=c(p,M,\delta,C_0,C_1,C_2)$.
\endproof
The coefficients $A_{\alpha \beta}$ defined in \eqref{coeffA} can be written as
$$A_{\alpha \beta }(x) := h_\varepsilon  \bigl(  |Du_\varepsilon|_\gamma \bigr) |Du_\varepsilon|_\gamma \biggl[  \gamma_{\alpha\beta}(x)+ \dfrac{h'_\varepsilon  \bigl(  |Du_\varepsilon|_\gamma \bigr) \gamma_{\alpha k}D_k u^i \gamma_{\beta \nu }D_\nu u^i}{h_\varepsilon  \bigl(  |Du_\varepsilon|_\gamma \bigr) |Du_\varepsilon|_\gamma}  \biggr].$$
They are degenerate elliptic  for $\varepsilon=0$, since $h_\varepsilon  \bigl(  |Du_\varepsilon|_\gamma \bigr) |Du_\varepsilon|_\gamma =0$ on the set $\bigl\{   |Du_\varepsilon|_\gamma \le 1 \bigr\}$. On the other hand, the function $U_\varepsilon$ has support in the set $B_R \cap \bigl\{   |Du_\varepsilon|_\gamma \ge 1 +\delta\bigr\}$. 
This fact will allow us to modify the coefficients $A_{\alpha\beta}$ on $B_R \cap \bigl\{   |Du_\varepsilon|_\gamma \le 1 +\delta\bigr\}$ and to prove an energy estimate for the function $U_\varepsilon$.

\begin{lem}\label{energyestUe}
    Let $\varepsilon \in (0,1]$ and $u_\varepsilon \in W^{1,\mathbf{p}}\bigl( B_R, \mathbb{R}^N \bigr)$ be the weak solution of the regularized system \eqref{regsystem} satisfying \eqref{bounded} and \eqref{delta} on $B_\rho(x_0) \subset B_{r_0}\Subset B_R$. Then, for any $k >0$ and any $\tau \in (0,1)$ we have
    $$\int_{B_{\tau \rho}(x_0)} |D \bigl(  U_\varepsilon-k \bigr)_+|^2 \ \mathrm{d}x \le \dfrac{c}{(1-\tau)^2\rho^2} \int_{B_{ \rho}(x_0)} \bigl( U_\varepsilon-k \bigr)^2_+ \ \mathrm{d}x +c |B_\rho(x_0) \cap \bigl\{ U_\varepsilon>k \bigr\}|,$$
    for a constant $c:=c(n,p,M,\delta,C_0,C_1,C_2)$.
\end{lem}
\proof We define the new coefficients $\tilde{A}_{\alpha\beta}$ by
\begin{equation*}
    \tilde{A}_{\alpha\beta}(x):=
    \begin{cases}
        \delta_{\alpha\beta} \quad \quad &\text{on} \ \ \ \{  x \in B_R : |Du_\varepsilon|_\gamma \le 1+ \delta \}\\
        A_{\alpha\beta} \quad \quad &\text{on} \ \  \ \{  x \in B_R : |Du_\varepsilon|_\gamma > 1+ \delta \}
    \end{cases}
\end{equation*}
From this definition, Lemma \ref{elleq} ad the fact that $ \mathrm{supp} \  U_\varepsilon \subset B_R \cap \bigl\{   |Du_\varepsilon|_\gamma \ge 1 +\delta\bigr\} $, we have that $U_\varepsilon$ satisfies the inequality
\begin{equation}\label{6.15}
    \int_{B_\rho(x_0)} \tilde{A}_{\alpha \beta}(x) D_\alpha U_\varepsilon D_\beta \xi \ \mathrm{d}x \le c \int_{B_\rho(x_0)} \bigl(  \xi + |\nabla \xi| \bigr) \ \mathrm{d}x
\end{equation}
for any non negative function $\xi \in \mathcal{C}^1_c \bigl(  B_R \bigr)$.
\\ We now investigate the upper bound and ellipticity of the coefficients $\tilde{A}_{\alpha\beta}$. We will show that there exist $0 < l \le L < + \infty$ both depending at most on $p,M,\delta,C_0,C_1$ and $C_2$ such that
$$l |\zeta|^2 \le \tilde{A}_{\alpha\beta} \zeta_\alpha \zeta_\beta \le L |\zeta|^2$$
for any $x \in B_{r_0}$ and $\zeta \in \mathbb{R}^n$. We start with the former one. On the set where $|Du_\varepsilon|_\gamma \le 1+ \delta$, the upper bound and ellipticity hold for $l=1=L$. On the set where $|Du_\varepsilon|_\gamma > 1+ \delta$, we have from Lemma \ref{ellipticity} that
\begin{align*}
    A_{\alpha \beta}(x) \zeta_\alpha \zeta_\beta= & \ |Du_\varepsilon(x)|_{\gamma}\mathcal{C}_\varepsilon(x,Du_\varepsilon(x)) (\zeta, \zeta )\\
    \le &   \ \bigl[  \varepsilon + \Lambda  \bigl( |Du_\varepsilon|_\gamma \bigr)  \bigl] |Du_\varepsilon|_\gamma |\zeta|^2_\gamma  \\
  \le & \ C_1 \bigl[ \varepsilon|Du_\varepsilon|_\gamma + \max \bigl\{  \bigl( |Du_\varepsilon|_\gamma-1 \bigr)^{p-1}, (p-1) |Du_\varepsilon|_\gamma \bigl( |Du_\varepsilon|_\gamma-1 \bigr)^{p-2}  \bigr\}   \bigr]  |\zeta|^2\\
  \le & \ C_1 \biggl[  M +\dfrac{pM^p}{\delta}\biggr] |\zeta|^2,
\end{align*}
hence $L=L(p,M,\delta,C_1)= C_1 \bigl[  M +\frac{pM^p}{\delta}\bigr]$.
On the other hand, it holds
\begin{align*}
     A_{\alpha \beta}(x) \zeta_\alpha \zeta_\beta \ge &
      \ C_0 \bigl[ \varepsilon|Du_\varepsilon|_\gamma + \min \bigl\{  \bigl( |Du_\varepsilon|_\gamma-1 \bigr)^{p-1}, (p-1) |Du_\varepsilon|_\gamma \bigl( |Du_\varepsilon|_\gamma-1 \bigr)^{p-2}  \bigr\}   \bigr]  |\zeta|^2\\
     \ge & \ C_0 \min \{  1, p-1 \} \delta^{p-1} |\zeta|^2,
\end{align*}
thus the ellipticity holds for $l=l(p,\delta,C_0)=  C_0 \min \{  1, p-1 \} \delta^{p-1} $.

Now, choosing $\xi = \eta^2 \bigl( U_\varepsilon -k \bigr)_+$ as test function in \eqref{6.15}, where $\eta \in \mathcal{C}^1_c \bigl(  B_\rho(x_0) \bigr)$ is a cut-off function with $\eta=1$ on $B_{\tau \rho(x_0)}$ and $|\nabla \eta | \le \frac{2}{(1-\tau )\rho}$, we get
\begin{align*}
    \int_{B_\rho(x_0)} &  \eta^2\tilde{A}_{\alpha \beta}(x) D_\alpha U_\varepsilon  D_\beta \bigl( U_\varepsilon -k \bigr)_+ \ \mathrm{d}x + 2 \int_{B_\rho(x_0)} \eta \tilde{A}_{\alpha \beta}(x) D_\alpha U_\varepsilon  D_\beta \eta \bigl( U_\varepsilon -k \bigr)_+  \ \mathrm{d}x \\
     & \le c \int_{B_\rho(x_0)} \bigl[\eta^2 \bigl( U_\varepsilon -k \bigr)_+ + |2\eta \nabla\eta  \bigl( U_\varepsilon -k \bigr)_+ + \eta^2 D \bigl( U_\varepsilon -k \bigr)_+  \bigr]\ \mathrm{d}x .
\end{align*}
By applying the Young's inequality, we can estimate the left-hand side of the previous inequality as follows
\begin{align*}
    \int_{B_\rho(x_0)} & \eta^2 \tilde{A}_{\alpha \beta}(x) D_\alpha U_\varepsilon  D_\beta \bigl( U_\varepsilon -k \bigr)_+ \ \mathrm{d}x + 2 \int_{B_\rho(x_0)} \eta \tilde{A}_{\alpha \beta}(x) D_\alpha U_\varepsilon  D_\beta \eta \bigl( U_\varepsilon -k \bigr)_+  \ \mathrm{d}x \\
    & \ge \dfrac{1}{2} \int_{B_\rho(x_0)}   \eta^2\tilde{A}_{\alpha \beta}(x) D_\alpha U_\varepsilon  D_\beta \bigl( U_\varepsilon -k \bigr)_+ \ \mathrm{d}x 
     -c \int_{B_\rho(x_0)} |\nabla \eta|^2 \bigl( U_\varepsilon -k \bigr)_+ \mathrm{d}x.
\end{align*}
On the other hand, we can bound the right-hand side as follows
\begin{align*}
    c \int_{B_\rho(x_0)} & \bigl[\eta^2 \bigl( U_\varepsilon -k \bigr)_+ + |2\eta \nabla\eta  \bigl( U_\varepsilon -k \bigr)_+ + \eta^2 D \bigl( U_\varepsilon -k \bigr)_+  \bigr]\ \mathrm{d}x \\
    \le & \ c \int_{B_\rho(x_0)}  \eta^2 \bigl( U_\varepsilon -k \bigr)_+^2 \ \mathrm{d}x + c |B_\rho(x_0) \cap \bigl\{  U_\varepsilon > k \bigr\}| \\
    & + \ \dfrac{c}{(1-\tau)^2\rho^2} \int_{B_\rho(x_0)}  \eta^2 \bigl( U_\varepsilon -k \bigr)_+^2 \ \mathrm{d}x  + \dfrac{l}{4} \int_{B_\rho(x_0)}  \eta^2 | D    \bigl( U_\varepsilon -k \bigr)_+ |^2 \ \mathrm{d}x .
\end{align*}
Combining the previous estimates, using the ellipticity bound for the coefficients $A_{\alpha\beta}$ and reabsorbing the term containing $ D    \bigl( U_\varepsilon -k \bigr)_+$ from the right-hand side to the left-hand side, we obtain the desired estimate. 
\endproof

\subsection{Energy estimates}
In this subsection, we turn our attention back to the result from Lemma \ref{mainintineq}, this time keeping the terms containing $\mathcal{A}_\varepsilon$ and $\mathcal{C}_\varepsilon$ on the left-hand of the integral inequality. We assume that hypothesis of Proposition \ref{ND} are in force. Using Lemmas \ref{mainintineq} and \ref{ellipticity}, for any increasing function $\Phi \in W^{1,\infty}_{loc} \bigl( \mathbb{R}_{\ge 0}, \mathbb{R}_{\ge 0} \bigr)$ and any non-negative function $\xi =\eta^2 \in \mathcal{C}^1_0 \bigl(  B_{\rho}(x_0) \bigr) $ we have
\begin{align*}
     \int_{B_{\rho}(x_0)}  \biggl[  \mathcal{A}_\varepsilon & (x,Du_\varepsilon)  \bigl( D^2u_\varepsilon, D^2u_\varepsilon \bigr) \Phi \bigl(  |Du_\varepsilon|_\gamma \bigr) + \mathcal{C}_\varepsilon(x,Du_\varepsilon)(F,F) \dfrac{\Phi' \bigl( |Du_\varepsilon|_\gamma  \bigr)}{|Du_\varepsilon|_\gamma} \biggr] \eta^2  \ \mathrm{d} x \notag\\
    \le & \ 8 \int_{B_{\rho}(x_0)} |\mathcal{C}_\varepsilon(x,Du_\varepsilon)(F,\nabla \eta) | \Phi \bigl( |Du_\varepsilon|_\gamma \bigr) \eta \ \mathrm{d} x \\
     & + c \int_{B_{\rho}(x_0)} \bigl[  \varepsilon + \Lambda \bigl( |Du_\varepsilon|_\gamma \bigr) \bigr] |Du_\varepsilon|_\gamma^3 \Phi' \bigl( |Du_\varepsilon|_\gamma \bigr) \eta^2 \ \mathrm{d} x \notag\\
    &+ c \int_{B_{\rho}(x_0)} \bigl[  \varepsilon + \Lambda \bigl( |Du_\varepsilon|_\gamma \bigr)  \bigl] |Du_\varepsilon|_\gamma^2 \Phi \bigl( |Du_\varepsilon|_\gamma \bigr) \eta^2 \ \mathrm{d}x \\
     &+c \int_{B_{\rho}(x_0)} \bigl[  \varepsilon + \Lambda \bigl( |Du_\varepsilon|_\gamma \bigr)  \bigl] |Du_\varepsilon|_\gamma^2 \Phi \bigl( |Du_\varepsilon|_\gamma \bigr) \eta |\nabla \eta| \ \mathrm{d}x\\
     =:& \ I_1+I_2+I_3+I_4.
\end{align*}
We estimate $I_1$ by Cauchy-Schwarz inequality and the upper bound from Lemma \ref{ellipticity} thus obtaining
\begin{align*}
    I_1 \le & \int_{B_{\rho}(x_0)} \mathcal{C}_\varepsilon(x,Du_\varepsilon)(F,F) \dfrac{ \Phi' \bigl( |Du_\varepsilon|_\gamma \bigr)}{|Du_\varepsilon|_\gamma} \eta^2 \ \mathrm{d} x \\
    &+  c  \int_{B_{\rho}(x_0)} \mathcal{C}_\varepsilon(x,Du_\varepsilon)(\nabla \eta,\nabla \eta) \dfrac{ \Phi^2 \bigl( |Du_\varepsilon|_\gamma \bigr) |Du_\varepsilon|_\gamma}{\Phi' \bigl(|Du_\varepsilon|_\gamma \bigr)}  \ \mathrm{d} x \\
    \le &  \int_{B_{\rho}(x_0)} \mathcal{C}_\varepsilon(x,Du_\varepsilon)(F,F) \dfrac{ \Phi' \bigl( |Du_\varepsilon|_\gamma \bigr)}{|Du_\varepsilon|_\gamma} \eta^2 \ \mathrm{d} x \\
    &+  c  \int_{B_{\rho}(x_0)} \bigl[ \varepsilon + \Lambda \bigl( |Du_\varepsilon|_\gamma \bigr)  \bigr]   \dfrac{ \Phi^2 \bigl( |Du_\varepsilon|_\gamma \bigr) |Du_\varepsilon|_\gamma}{\Phi' \bigl(|Du_\varepsilon|_\gamma \bigr)} |\nabla \eta |^2 \ \mathrm{d} x .
\end{align*}
Moreover, for $I_4$ we have
\begin{align*}
    I_4 \le & \  c \ \int_{B_{\rho}(x_0)} \bigl[  \varepsilon + \Lambda \bigl( |Du_\varepsilon|_\gamma \bigr)  \bigl] |Du_\varepsilon|_\gamma^3 \Phi' \bigl( |Du_\varepsilon|_\gamma \bigr) \eta^2 \ \mathrm{d}x \\
    & +  c  \int_{B_{\rho}(x_0)} \bigl[ \varepsilon + \Lambda \bigl( |Du_\varepsilon|_\gamma \bigr)  \bigr]   \dfrac{ \Phi^2 \bigl( |Du_\varepsilon|_\gamma \bigr) |Du_\varepsilon|_\gamma}{\Phi' \bigl(|Du_\varepsilon|_\gamma \bigr)} |\nabla \eta |^2 \ \mathrm{d} x .
\end{align*}
Inserting the estimates above and re-absorbing the term containing $\mathcal{C}_\varepsilon(x,Du_\varepsilon)(F,F)$ into the left-hand side, we find
\begin{align}\label{6.16}
     \int_{B_{\rho}(x_0)}    \mathcal{A}_\varepsilon & (x,Du_\varepsilon)  \bigl( D^2u_\varepsilon, D^2u_\varepsilon \bigr) \Phi \bigl(  |Du_\varepsilon|_\gamma \bigr)  \eta^2  
 \ \mathrm{d} x \notag \\
 \le & \ c \int_{B_{\rho}(x_0)} \bigl[  \varepsilon + \Lambda \bigl( |Du_\varepsilon|_\gamma \bigr)  \bigl] |Du_\varepsilon|_\gamma^2 \bigl(  \Phi \bigl( |Du_\varepsilon|_\gamma \bigr) + \Phi' \bigl( |Du_\varepsilon|_\gamma \bigr) |Du_\varepsilon|_\gamma  \bigr) \eta^2 \ \mathrm{d}x \notag\\
  & + \ c  \int_{B_{\rho}(x_0)} \bigl[ \varepsilon + \Lambda \bigl( |Du_\varepsilon|_\gamma \bigr)  \bigr]   \dfrac{ \Phi^2 \bigl( |Du_\varepsilon|_\gamma \bigr) |Du_\varepsilon|_\gamma}{\Phi' \bigl(|Du_\varepsilon|_\gamma \bigr)} |\nabla \eta |^2 \ \mathrm{d} x .
\end{align}
Now, we choose $\Phi (t):= (t-1)^p_+ \tilde{\Phi}(t)$, where $\tilde{\Phi}(t) \in W^{1,\infty}_{loc}\bigl(  \mathbb{R}_{\ge 0},  \mathbb{R}_{\ge 0} \bigr)$ is non-decreasing. Note that
$$\Phi' (t) = (t-1)^{p-1}_+ \bigl[  p  \tilde{\Phi}(t) + (t-1) \tilde{\Phi}'(t) \bigr].$$
For $t \in [0, 1 + 2 \mu ]$ we have
\begin{align*}
    \bigl[ \varepsilon + \Lambda \bigl( t \bigr) \bigr] \dfrac{\Phi^2(t)t}{\Phi'(t)} = & \   \bigl[ \varepsilon + \Lambda \bigl( t \bigr) \bigr] \ \dfrac{(t-1)^{p+1}_+  t  \tilde{\Phi}^2(t)}{p \tilde{\Phi}(t)+(t-1)_+\tilde{\Phi}'(t)} \\
    \le & \  \bigl[  \varepsilon (t-1)^{p+1}_+  t+ \max \bigl\{ (t-1)^{2p}_+, (p-1) (t-1)_+^{2p-1}t   \bigr\} \bigr]     \dfrac{ \tilde{\Phi}^2(t)}{p \tilde{\Phi}(t)+(t-1)_+\tilde{\Phi}'(t)} \\
    \le & \ c (\delta)\bigl[  \varepsilon\mu^{p+2} + \mu^{2p} \bigr]  \dfrac{ \tilde{\Phi}^2(t)}{p \tilde{\Phi}(t)+(t-1)_+\tilde{\Phi}'(t)} \\
    \le & \ c(p,M,\delta)  \dfrac{ \mu^{2p} \tilde{\Phi}^2(t)}{p \tilde{\Phi}(t)+(t-1)_+\tilde{\Phi}'(t)},
\end{align*}
where we used $\delta \le \mu \le M$, which implies $t \le 1 + 2 \mu \le \bigl(  2+ \frac{1}{\delta} \bigr) \mu$ and on the other hand $\mu^2 \le \max \bigl\{ \delta^{2-p}, M^{2-p} \bigr\} \mu^p   $. Next, we estimate
\begin{align*}
       \Phi(t)+t \Phi'(t) = & \   (t-1)_+^p \tilde{\Phi}(t) +pt (t-1)_+^{p-1} \tilde{\Phi}(t) + t(t-1)_+^p \tilde{\Phi}'(t)  \\
    \le & \ (1+2 \mu)^{p} \tilde{\Phi}(t) + p (1+2 \mu)^{p} \bigl[  \tilde{\Phi}(t)+ t \tilde{\Phi}'(t) \bigr] \\
    \le & \ c(p,M,\delta) \bigl[  \tilde{\Phi}(t)+ t \tilde{\Phi}'(t) \bigr].
\end{align*}
Hence, we derive
\begin{align*}
    \bigl[ \varepsilon + \Lambda \bigl( t \bigr) \bigr] t^2 \bigl[  \Phi(t)+t \Phi'(t) \bigr] \le & \ c \bigl[  \varepsilon t^2  + \max \bigl\{ (t-1)^{p-1}_+t, (p-1) (t-1)_+^{p-2}t   \bigr\} \bigr] \bigl[  \tilde{\Phi}(t)+ t \tilde{\Phi}'(t) \bigr] \\
    \le & \ c \bigl[   \mu^2  + \max \bigl\{ \mu^p, (p-1) \mu^{p-1}   \bigr\} \bigr] \bigl[  \tilde{\Phi}(t)+ t \tilde{\Phi}'(t) \bigr]\\
    \le & \ c(p,M,\delta)  \bigl[  \tilde{\Phi}(t)+ t \tilde{\Phi}'(t) \bigr].
\end{align*}
Due to assumption \eqref{bounded} and \eqref{mu}, we have that $|Du_\varepsilon|_\gamma \le 1+ 2 \mu \le \bigl( 2+ \frac{1}{\delta}  \bigr) \mu$ on the ball $B_\rho(x_0)$. This allows us to use the preceding estimates in \eqref{6.16} to bound the right-hand side from above. Thus, we get
\begin{align*}
     \int_{B_{\rho}(x_0)}    \mathcal{A}_\varepsilon & (x,Du_\varepsilon)  \bigl( D^2u_\varepsilon, D^2u_\varepsilon \bigr) \Phi \bigl(  |Du_\varepsilon|_\gamma \bigr)  \eta^2  
 \ \mathrm{d} x  \\
 \le & \ c \int_{B_{\rho}(x_0)} \bigl[  \tilde{\Phi}\bigl( |Du_\varepsilon|_\gamma \bigr)+ t \tilde{\Phi}' \bigl(|Du_\varepsilon|_\gamma \bigr) \bigr] \eta^2 \ \mathrm{d}x \\
 + & \ c \int_{B_{\rho}(x_0)}  \dfrac{ \mu^{2p} \tilde{\Phi}^2 \bigl(|Du_\varepsilon|_\gamma \bigr)}{p \tilde{\Phi}\bigl(|Du_\varepsilon|_\gamma \bigr)+\bigl(|Du_\varepsilon|_\gamma -1 \bigr)_+\tilde{\Phi}'(t)} |\nabla \eta|^2    \ \mathrm{d}x.
\end{align*}
Then, applying Lemma \ref{lem11}, we obtain
\begin{align} \label{gradg}
  \int_{B_{\rho}(x_0)}  \big|D \bigl[ & g \bigl( |Du_\varepsilon|_\gamma \bigr)  Du_\varepsilon \bigr]\big|^2  \tilde{\Phi} \bigl(  |Du_\varepsilon|_\gamma \bigr)  \eta^2  
 \ \mathrm{d} x \notag \\
  \le & \ c \int_{B_{\rho}(x_0)}    g'\bigl( |Du_\varepsilon|_\gamma \bigr)^2  |Du_\varepsilon|^4_\gamma       \tilde{\Phi} \bigl(  |Du_\varepsilon|_\gamma \bigr)  \eta^2  
 \ \mathrm{d} x  \notag\\
  + & \ c \int_{B_{\rho}(x_0)} \bigl[  \tilde{\Phi}\bigl( |Du_\varepsilon|_\gamma \bigr)+ t \tilde{\Phi}' \bigl(|Du_\varepsilon|_\gamma \bigr) \bigr] \eta^2 \ \mathrm{d}x \notag \\
 + & \ c \int_{B_{\rho}(x_0)}  \dfrac{ \mu^{2p} \tilde{\Phi}^2 \bigl(|Du_\varepsilon|_\gamma \bigr)}{p \tilde{\Phi}\bigl(|Du_\varepsilon|_\gamma \bigr)+\bigl(|Du_\varepsilon|_\gamma -1 \bigr)_+\tilde{\Phi}'(t)} |\nabla \eta|^2    \ \mathrm{d}x .
\end{align}
Using the second inequality in Lemma \ref{lem6}, we can estimate the first integral on the right-hand side of the previous inequality as follows
\begin{align*}
     & \int_{B_{\rho}(x_0)}    g' \bigl( |Du_\varepsilon|_\gamma \bigr)^2  |Du_\varepsilon|^4_\gamma       \tilde{\Phi} \bigl(  |Du_\varepsilon|_\gamma \bigr)  \eta^2 \ \mathrm{d} x  \\
     &  \le    \int_{B_{\rho}(x_0)}  \dfrac{p^2}{p-1}\bigl[h\bigl(  |Du_\varepsilon|_\gamma 
     \bigr)+h' \bigl(  |Du_\varepsilon|_\gamma 
     \bigr) |Du_\varepsilon|_\gamma \bigr]\bigl( |Du_\varepsilon|_\gamma-1 \bigr)^p_+ |Du_\varepsilon|_\gamma^2        \tilde{\Phi} \bigl(  |Du_\varepsilon|_\gamma \bigr)  \eta^2
 \ \mathrm{d} x.
\end{align*}
Now, since $|Du_\varepsilon|_\gamma  \le \bigl( 2+ \frac{1}{\delta}  \bigr) \mu$ on $B_\rho(x_0)$ and the function $t \mapsto [h(t)+h'(t)t](t-1)_+^p  = (p-1) (t-1)_+^{2p-2} $ is increasing, we derive that
\begin{align*}
      \int_{B_{\rho}(x_0)}  &  g' \bigl( |Du_\varepsilon|_\gamma \bigr)^2  |Du_\varepsilon|^4_\gamma       \tilde{\Phi} \bigl(  |Du_\varepsilon|_\gamma \bigr)  \eta^2 \ \mathrm{d} x   
    \le  c(p,\delta,M)  \int_{B_{\rho}(x_0)}         \tilde{\Phi} \bigl(  |Du_\varepsilon|_\gamma \bigr)  \eta^2 \ \mathrm{d} x .
\end{align*}
Inserting the preceding estimate in \eqref{gradg}, we eventually get
\begin{align} \label{gradg2}
  \int_{B_{\rho}(x_0)}  \big|D \bigl[ & g \bigl( |Du_\varepsilon|_\gamma \bigr)  Du_\varepsilon \bigr]\big|^2  \tilde{\Phi} \bigl(  |Du_\varepsilon|_\gamma \bigr)  \eta^2  
 \ \mathrm{d} x \notag \\
  \le & \ c \int_{B_{\rho}(x_0)} \bigl[  \tilde{\Phi}\bigl( |Du_\varepsilon|_\gamma \bigr)+ t \tilde{\Phi}' \bigl(|Du_\varepsilon|_\gamma \bigr) \bigr] \eta^2 \ \mathrm{d}x \notag \\
 + & \ c \int_{B_{\rho}(x_0)}  \dfrac{ \mu^{2p} \tilde{\Phi}^2 \bigl(|Du_\varepsilon|_\gamma \bigr)}{p \tilde{\Phi}\bigl(|Du_\varepsilon|_\gamma \bigr)+\bigl(|Du_\varepsilon|_\gamma -1 \bigr)_+\tilde{\Phi}'(t)} |\nabla \eta|^2    \ \mathrm{d}x ,
\end{align}
for a constant $c:=c(n,p,M,\delta,C_0,C_1,C_2)$.

Choosing $\tilde{\Phi} \equiv 1$, we derive the following energy inequality.
\begin{lem}\label{enineq1}
    Let $\varepsilon \in (0,1]$ and $u_\varepsilon \in W^{1,p} \bigl(   B_R, \mathbb{R}^N\bigr)$ be a weak solution of the regularized system \eqref{regsystem} such that hypotheses \eqref{bounded}, \eqref{delta} and \eqref{mu} are in force on $B_\rho (x_0) \subset B_{r_0} \Subset B_R$. Then, for any $\tau \in (0,1)$ there holds
        \begin{equation*}
           \int_{B_{\tau\rho}(x_0)}  \big|D \bigl[  g \bigl( |Du_\varepsilon|_\gamma \bigr)  Du_\varepsilon \bigr]\big|^2  
 \ \mathrm{d} x  \le c \biggl[ \dfrac{\mu^{2p}}{\rho^2(1-\tau)^2}+1   \biggr]|B_\rho|,
        \end{equation*}
        for a constant $c:=c(n,p,M,\delta,C_0,C_1,C_2)$.
\end{lem}
On the other hand, choosing $ \tilde{\Phi}(t):=\bigl(t-1-\delta-(1-2\nu)\mu \bigr)^2_+$ for some $\nu \in \bigl(0, \frac{1}{8} \bigr]$, it holds
\begin{lem}\label{enineq2}
 Let $\nu \in \bigl(0, \frac{1}{8} \bigr]$, $\varepsilon \in (0,1]$ and $u_\varepsilon \in W^{1,p} \bigl(   B_R, \mathbb{R}^N\bigr)$ be a weak solution of the regularized system \eqref{regsystem} such that hypotheses \eqref{bounded}, \eqref{delta}, \eqref{mu} and \eqref{measureND} are in force on $B_\rho (x_0) \subset B_{r_0} \Subset B_R$. Then, for any $\tau \in (0,1)$ we have
    \begin{equation*}
           \int_{E^\nu_{\tau\rho}(x_0)}  \big|D \bigl[  g \bigl( |Du_\varepsilon|_\gamma \bigr)  Du_\varepsilon \bigr]\big|^2   
 \ \mathrm{d} x  \le c \biggl[ \dfrac{\mu^{2p} \nu}{\rho^2(1-\tau)^2}+\dfrac{1}{\nu}   \biggr]|B_\rho|,
        \end{equation*}
        for a constant $c:=c(n,p,M,\delta,C_0,C_1,C_2)$.
\end{lem}

\section{The non-degenerate regime}\label{secND}

Here we give the proof of Proposition \ref{ND}. Throughout this section we will assume the following assumption. For given $\varepsilon \in (0,1]$ we denote by $u_\varepsilon \in W^{1,p} \bigl( B_R, \mathbb{R}^N \bigr)$ the unique weak solution of the regularized system \eqref{regsystem}. Moreover, we assume that for some $\delta \in (0,1]$ and $\mu > \delta$ and a ball $B_{2 \rho}(x_0) \subset B_{r_1}$ with $\rho \le 1$ assumptions \eqref{bounded}, \eqref{delta} and \eqref{mu} are in force. We denote by
\begin{equation}\label{intmean}
    \mathbf{\Phi} (x_0,\rho):= \fint_{B_\rho(x_0)} \big|Du_\varepsilon - ( D u_\varepsilon )_{x_0, \rho} \bigr|^2 \ \mathrm{d}x
\end{equation}
the excess function of $Du_\varepsilon$.
\subsection{Higher integrability}
The following higher integrability result holds.

\begin{lem}\label{higherint}
    Under the general assumptions of Section \ref{secND}, there exist $\vartheta \in \bigl( 0, \frac{1}{2} \bigr]$ and $c>0$, depending at most on $n,p,M,\delta,C_0$ and $C_1$, such that for any $\xi \in \mathbb{R}^{nN}$ satifying 
    $$1+ \delta + \dfrac{1}{4} \mu \le  |\xi|_{\gamma(x)} \le 1+ \delta +\mu$$
    for any $x \in B_\rho(x_0)$,
    we have
    $$ \fint_{B_{\rho/2}(x_0)} |D u_\varepsilon - \xi|^{2(1+\vartheta)} \ \mathrm{d}x \le c \biggl[  \fint_{B_{\rho}(x_0)} |D u_\varepsilon - \xi|^{2} \ \mathrm{d}x   \biggr]^{1+ \vartheta}+c \rho^{2(1+\vartheta)}.$$
\end{lem}
\proof 
We consider a ball $B_s(z_0) \subset B_\rho(x_0)$ and we test the weak form \eqref{regsystem2} by the testing function 
$\varphi := \eta^2 w$, where $w:= u_\varepsilon - ( u_\varepsilon )_{z_0,s}-\xi (x-z_0)$ and $\eta \in \mathcal{C}_0^1 \bigl( B_s(z_0) \bigr)$ is  a cut-off function with $\eta=1$ in $B_{s/2}(z_0)$, $0 \le \eta \le 1$ and $|\nabla \eta | \le \frac{4}{s}$. We obtain
\begin{align*}
    0= & \int_{B_s(z_0)} \langle \mathbf{A}_\varepsilon \bigl( x, Du_\varepsilon \bigr) , D \varphi \rangle_{\gamma(x)}\ \mathrm{d} x \\
    = & \int_{B_s(z_0)} \langle \mathbf{A}_\varepsilon \bigl( x, Du_\varepsilon \bigr) - \mathbf{A}_\varepsilon (x, \xi) , D \varphi \rangle_{\gamma(x)}\ \mathrm{d} x  \\
    &+ \int_{B_s(z_0)} \bigl[\langle \mathbf{A}_\varepsilon ( x, \xi) , D \varphi \rangle_{\gamma(x)}- \langle \mathbf{A}_\varepsilon (x_0, \xi) , D \varphi \rangle_{\gamma(x_0)} \bigr] \mathrm{d} x \\
    =&  \int_{B_s(z_0)} \langle \mathbf{A}_\varepsilon \bigl( x, Du_\varepsilon \bigr) - \mathbf{A}_\varepsilon (x, \xi) , \eta^2 D w \rangle_{\gamma(x)}\ \mathrm{d} x \\
    &+ \int_{B_s(z_0)} \langle \mathbf{A}_\varepsilon \bigl( x, Du_\varepsilon \bigr) - \mathbf{A}_\varepsilon (x, \xi) , 2 \eta \nabla \eta \otimes w  \rangle_{\gamma(x)}\ \mathrm{d} x \\
    & + \int_{B_s(z_0)} \bigl[\langle \mathbf{A}_\varepsilon ( x, \xi) , D \varphi \rangle_{\gamma(x)}- \langle \mathbf{A}_\varepsilon (x_0, \xi) , D \varphi \rangle_{\gamma(x_0)}\ \bigr] \mathrm{d} x \\
    = & I+II+III,
\end{align*}
which yields
\begin{equation}\label{hi1}
    I \le |II|+ |III|.
\end{equation}
We use the monotonicity formula of $\mathbf{A}_\varepsilon$ from Lemma \ref{monotonicity} in order to estimate the integral $I$ from below. Due to our assumption on $\xi$ and \eqref{bounded}, we have
$$4 |\xi|_{\gamma(x)} \ge 4(1+\delta )+ \mu \ge |Du_\varepsilon|_{\gamma(x)}$$
which implies
$$5 |\xi|_{\gamma(x)} \ge |Du_\varepsilon|_{\gamma(x)}+|\xi|_{\gamma(x)}.$$
Therefore, we obtain
\begin{equation*}
    I \ge \hat{c} \int_{B_s(z_0)} \biggl[  \varepsilon + \lambda \dfrac{\bigl(  |\xi|_{\gamma(x)}-1 \bigr)^p}{|\xi|_{\gamma(x)}^2}  \biggr] \eta^2 |Dw|^2 \ \mathrm{d} x
\end{equation*}
for a constant $\hat{c}:= \hat{c}(C_0,C_1)$, where we denoted $\lambda:= \dfrac{\min \{ 1,p-1 \}}{5 \cdot 2^{p+1}}$.
\\To bound the second integral we use the upper bound from Lemma \ref{monotonicity}. Moreover, we observe that $\bigl(  |Du_\varepsilon|_{\gamma(x_0)}-1  \bigr)_+ \le 4 \bigl(  |\xi|_{\gamma(x_0)}-1 \bigr)$ due to our assumption on $\xi$. Together with Young's inequality, this yields
\begin{align*}
    |II| \le & \  c \int_{B_s(z_0)} |\mathbf{A}_\varepsilon \bigl( x,Du_\varepsilon \bigr)- \mathbf{A}_\varepsilon(x, \xi  )| |\nabla \eta| |w| \eta \ \mathrm{d}x\\
    \le & \ c \int_{B_s(z_0)}\bigl[  \varepsilon+ \bigl(  |\xi|_{\gamma(x)}-1 \bigr)^{p-2} \bigr] |Dw||\nabla \eta| |w| \eta \ \mathrm{d}x\\
    \le & \ \hat{c} \int_{B_s(z_0)} \biggl[  \dfrac{\varepsilon}{2} + \dfrac{\lambda}{4} \dfrac{\bigl(  |\xi|_{\gamma(x)}-1 \bigr)^p}{|\xi|_{\gamma(x)}^2}  \biggr] \eta^2 |Dw|^2 \ \mathrm{d} x\\
    & +  c \int_{B_s(z_0)}\bigl[  \varepsilon+ |\xi|^2_{\gamma(x)}\bigl(  |\xi|_{\gamma(x)}-1 \bigr)^{p-4} \bigr] |\nabla \eta|^2 |w|^2  \ \mathrm{d}x
\end{align*}
for a constant $c:=c(p,C_0,C_1)$. 
\\Now, we take care of the integral $III$. Using the third formula in  Lemma \ref{monotonicity}, we derive that
\begin{align*}
    |III| \le & \  c \int_{B_s(z_0)} \biggl[\varepsilon+ \dfrac{\bigl[  |\xi|_{\gamma(x)}-1+ \bigl( |\xi|_{\gamma(x_0)}-1 \bigr)\bigr]^{p-1}}{|\xi|_{\gamma(x)}-1}  \biggr]|x-x_0| |\xi| |D\varphi| \ \mathrm{d}x \\
    \le &  \ c \rho \int_{B_s(z_0)} \biggl[\varepsilon+ \dfrac{\bigl[  |\xi|_{\gamma(x)}-1+ \bigl( |\xi|_{\gamma(x_0)}-1 \bigr)\bigr]^{p-1}}{|\xi|_{\gamma(x)}-1}  \biggr]|\xi| \eta^2 |Dw| \ \mathrm{d}x \\ 
     & + c \rho \int_{B_s(z_0)} \biggl[\varepsilon+ \dfrac{\bigl[  |\xi|_{\gamma(x)}-1+ \bigl( |\xi|_{\gamma(x_0)}-1 \bigr)\bigr]^{p-1}}{|\xi|_{\gamma(x)}-1}  \biggr]|\xi| \eta |\nabla \eta| |w| \ \mathrm{d}x  .
\end{align*}
Observe that $ |\xi|_{\gamma(x_0)}-1   \le 4 \bigl(  |\xi|_{\gamma(x)}-1 \bigr)$ due to our assumption on $\xi$. This together with the Young's inequality allows us to estimate
\begin{align*}
    |III| \le & \ c \rho \int_{B_s(z_0)}\bigl[  \varepsilon+ \bigl(  |\xi|_{\gamma(x)}-1 \bigr)^{p-2} \bigr] |\xi| \eta^2  |Dw|  \ \mathrm{d}x \\
    & + \ c \rho \int_{B_s(z_0)}\bigl[  \varepsilon+ \bigl(  |\xi|_{\gamma(x)}-1 \bigr)^{p-2} \bigr] |\xi| \eta |\nabla \eta| |w| \ \mathrm{d}x \\
    \le & \ \hat{c} \int_{B_s(z_0)} \biggl[  \dfrac{\varepsilon}{2} + \dfrac{\lambda}{4} \dfrac{\bigl(  |\xi|_{\gamma(x)}-1 \bigr)^p}{|\xi|_{\gamma(x)}^2}  \biggr] \eta^2 |Dw|^2 \ \mathrm{d} x \\
     & + c \rho^2  \int_{B_s(z_0)}\bigl[  \varepsilon+ |\xi|^2_{\gamma(x)}\bigl(  |\xi|_{\gamma(x)}-1 \bigr)^{p-4} \bigr] |\xi|^2 \eta^2  \ \mathrm{d}x\\
     & +  c  \int_{B_s(z_0)}\bigl[  \varepsilon+ |\xi|^2_{\gamma(x)}\bigl(  |\xi|_{\gamma(x)}-1 \bigr)^{p-4} \bigr] |w|^2 |\nabla \eta|^2  \ \mathrm{d}x\\
     & +  c \rho^2  \int_{B_s(z_0)}\bigl[  \varepsilon+ \bigl(  |\xi|_{\gamma(x)}-1 \bigr)^{p} \bigr] \eta^2  \ \mathrm{d}x.
\end{align*}
Joining the preceding estimates, re-absorbing the terms containing $Dw$ in the right-hand side into the left-hand side, we find that
\begin{align*}
    \dfrac{\hat{c} \lambda}{2} \int_{B_s(z_0)} &
     \dfrac{\bigl(  |\xi|_{\gamma(x)}-1 \bigr)^p}{|\xi|_{\gamma(x)}^2}   \eta^2 |Dw|^2 \ \mathrm{d} x\\
     \le  & \  c  \int_{B_s(z_0)}\bigl[  \varepsilon+ |\xi|^2_{\gamma(x)}\bigl(  |\xi|_{\gamma(x)}-1 \bigr)^{p-4} \bigr] |w|^2 |\nabla \eta|^2  \ \mathrm{d}x\\
     & + c \rho^2  \int_{B_s(z_0)}\bigl[  \varepsilon+ |\xi|^2_{\gamma(x)}\bigl(  |\xi|_{\gamma(x)}-1 \bigr)^{p-4} \bigr] |\xi|^2 \eta^2  \ \mathrm{d}x\\
    & +  c \rho^2  \int_{B_s(z_0)}\bigl[  \varepsilon+ \bigl(  |\xi|_{\gamma(x)}-1 \bigr)^{p} \bigr] \eta^2  \ \mathrm{d}x,
\end{align*}
where the constants $\hat{c}$ and $c$ depend on $p,C_0,C_1$. Due to our assumption on $\xi$ and the properties of $\eta$, we have
$$ \int_{B_{s/2}(z_0)} |Dw|^2 \ \mathrm{d} x \le \ \dfrac{c}{s^2} \int_{B_s(z_0)} |w|^2 \ \mathrm{d}x + c \rho^2 |B_s|, $$
for a constant $c:=c(p,M,\delta,C_0,C_1)$. Now, with an application of H\"older's inequality and Sobolev-Poincarè inequality, we derive a reverse H\"older's inequality of the form
\begin{align*}
    \fint_{B_{s/2}(z_0)} |Dw|^2 \ \mathrm{d} x \le &  \ \dfrac{c}{s^2} \fint_{B_s(z_0)} |w|^2 \ \mathrm{d}x + c \rho^2 
    \le  \ c \biggl[  \fint_{B_s(z_0)} |Dw|^{\frac{2n}{n+2}} \ \mathrm{d}x \biggr]^{\frac{n+2}{n}} +c \rho^2,
\end{align*}
for a constant $c:=c(n,p,M,\delta,C_0,C_1)$.
Now, since $Dw= Du_\varepsilon - \xi$, the higher integrability follows by Gehring's lemma \cite[Theorem 6.6]{Giusti}.
\endproof

\subsection{Comparison with a linear system}
In this section we will consider the weak solution $v \in u_\varepsilon + W^{1,2}_0 \bigl( B_{\rho/2}(x_0) , \mathbb{R}^N \bigr)$ of the linear elliptic system
\begin{equation}\label{linear}
    \int_{B_{\rho/2}(x_0)} \mathcal{B}_\varepsilon \bigl(  x, ( Du_\varepsilon )_{x_0,\rho/2} \bigr) (Dv,D\varphi) \ \mathrm{d}x=0
\end{equation}
for any $\varphi \in W^{1,2}_0 \bigl( B_{\rho/2}(x_0) , \mathbb{R}^N \bigr) $.

\begin{lem}\label{comparison}
    Under the general assumptions of Section \ref{secND}, assume that for every $x \in B_\rho(x_0)$
    \begin{equation}\label{boundmean}
        1+\delta+\dfrac{1}{4} \mu \le \big| ( Du_\varepsilon )_{x_0,\rho/2}  \big|_{\gamma(x)} \le 1+\delta+\mu.
    \end{equation}
    Then, there exist $\vartheta:= \vartheta (n,p,M,\delta,C_0,C_1) \in \bigl(0 , \frac{1}{2} \bigr]$ and $c:= c (n,p,M,\delta,C_0,C_1)>0$ such that
    $$ \fint_{B_{\rho/2}} |Du_\varepsilon-Dv|^2 \ \mathrm{d}x \le \ c \biggr[  \dfrac{ \mathbf{\Phi} (x_0,\rho)}{\mu^2} \biggl]^\vartheta  \mathbf{\Phi} (x_0,\rho) +c \rho^2.$$
\end{lem}
\proof Throughout the proof we omit the reference to the center $x_0$ and we write $B_\rho$ instead of $B_\rho(x_0)$. Moreover, we denote $\xi := (Du_\varepsilon)_{\rho/2}$. Using the weak formulation \eqref{regsystem2} of the regularized elliptic system, we obtain
\begin{align*}
    0 = & \int_{B_{\rho/2}} \langle \mathbf{A}_\varepsilon \bigl( x, Du_\varepsilon \bigr)  , D \varphi \rangle_{\gamma(x)}\ \mathrm{d} x \\
    = & \int_{B_{\rho/2}} \langle \mathbf{A}_\varepsilon \bigl( x, Du_\varepsilon \bigr) - \mathbf{A}_\varepsilon (x, \xi) , D \varphi \rangle_{\gamma(x)} \ \mathrm{d} x \\
    & + \int_{B_{\rho/2}} \bigl[\langle \mathbf{A}_\varepsilon ( x, \xi) , D \varphi \rangle_{\gamma(x)}- \langle \mathbf{A}_\varepsilon (x_0, \xi) , D \varphi \rangle_{\gamma(x_0)} \bigr] \mathrm{d} x 
\end{align*}
for every  $\varphi \in W^{1,2}_0 \bigl( B_{\rho/2}(x_0) , \mathbb{R}^N \bigr) $. Using also the fact that $v$ is a weak solution of the linear elliptic system \eqref{linear} we find that
\begin{align*}
    \int_{B_{\rho/2}} & \mathcal{B}_\varepsilon  (  x, \xi ) (Du_\varepsilon-Dv,D\varphi) \ \mathrm{d}x \\
    =&  \int_{B_{\rho/2}} \mathcal{B}_\varepsilon  (  x, \xi )(Du_\varepsilon,D\varphi) \ \mathrm{d}x \\
     =&  \int_{B_{\rho/2}} \mathcal{B}_\varepsilon  (  x, \xi )(Du_\varepsilon-\xi,D\varphi)  \ \mathrm{d}x + \int_{B_{\rho/2}} \bigl[ \mathcal{B}_\varepsilon  (  x, \xi ) (\xi,D\varphi)- \mathcal{B}_\varepsilon  (  x_0, \xi ) (\xi,D\varphi) \bigr]\ \mathrm{d}x\\
     =&  \int_{B_{\rho/2}} \bigl[\mathcal{B}_\varepsilon  (  x, \xi )(Du_\varepsilon-\xi,D\varphi) -  \langle \mathbf{A}_\varepsilon \bigl( x, Du_\varepsilon \bigr) - \mathbf{A}_\varepsilon (x, \xi) , D \varphi \rangle_{\gamma(x)} \bigr] \ \mathrm{d}x \\
     &+  \int_{B_{\rho/2}} \bigl[\langle \mathbf{A}_\varepsilon ( x, \xi) , D \varphi \rangle_{\gamma(x)}- \langle \mathbf{A}_\varepsilon (x_0, \xi) , D \varphi \rangle_{\gamma(x_0)} \bigr] \mathrm{d} x \\
     &+ \int_{B_{\rho/2}} \bigl[ \mathcal{B}_\varepsilon  (  x, \xi ) (\xi,D\varphi)- \mathcal{B}_\varepsilon  (  x_0, \xi ) (\xi,D\varphi) \bigr]\ \mathrm{d}x.
\end{align*}
In order to bound the right-hand side of the preceding equality, we  exploit Lemma \ref{linearization}, the third inequality in Lemma \ref{monotonicity} and Lemma \ref{bilformB}. This is possible thanks to assumptions \eqref{bounded} and \eqref{boundmean}. Hence, we obtain
\begin{align}\label{estimatecomp}
    \int_{B_{\rho/2}(x_0)} & \mathcal{B}_\varepsilon  (  x, \xi ) (Du_\varepsilon-Dv,D\varphi) \ \mathrm{d}x \notag \\
    \le &  \ c \mu^{p-3} \int_{B_{\rho/2}(x_0)} |Du_\varepsilon-\xi|^2 |D\varphi| \ \mathrm{d}x  \notag \\
    &+ c \int_{B_{\rho/2}} \biggl[  \varepsilon + \dfrac{\bigl[|\xi|_{\gamma(x)}-1+ \bigl( |\xi|_{\gamma(x_0)}-1  \bigr) \bigr]^{p-1}}{|\xi|_{\gamma(x)}-1}  \biggr] |x-x_0| |\xi||D \varphi| \ \mathrm{d}x \notag \\
     &+ c (\varepsilon+\mu^{p-2})\int_{B_{\rho/2}} |x-x_0| |\xi|^2|D \varphi| \ \mathrm{d}x. 
\end{align}
Since $u_\varepsilon-v \in W^{1,p}_0 \bigl(  B_{\rho/2}, \mathbb{R}^N \bigr)$, we may choose $\varphi:= u_\varepsilon-v $ as testing function in \eqref{estimatecomp}. Together with the bound from below from Lemma \ref{ellipticity} and H\"older's inequality this leads to
\begin{align*}
    & \hat{c}   \mu^{p-2}  \int_{B_{\rho/2}} |Du_\varepsilon-Dv|^2 \ \mathrm{d}x \\
    & \le   \ c \mu^{p-3} \biggl( \int_{B_{\rho/2}} |Du_\varepsilon-\xi|^4  \ \mathrm{d}x  \biggr)^{\frac{1}{2}} \biggl( \int_{B_{\rho/2}} |Du_\varepsilon-Dv|^4  \ \mathrm{d}x  \biggr)^{\frac{1}{2}}+ c \rho\int_{B_{\rho/2}} |Du_\varepsilon-Dv|  \ \mathrm{d}x \\
    & \le   \ c \mu^{p-3} \biggl( \int_{B_{\rho/2}} |Du_\varepsilon-\xi|^4 \ \mathrm{d}x  \biggr)^{\frac{1}{2}} \biggl( \int_{B_{\rho/2}} |Du_\varepsilon-Dv|^4  \ \mathrm{d}x  \biggr)^{\frac{1}{2}}
    + c \rho^{\frac{n+2}{n}} \biggl(\int_{B_{\rho/2}} |Du_\varepsilon-Dv|^2  \ \mathrm{d}x  \biggr)^\frac{1}{2},
\end{align*}
for constants $\hat{c}$ and $c$ depending at most on $p,M,\delta ,C_0$ and $C_1$.
Next, we divide both sides by $$  \hat{c}  \mu^{p-2}  \biggl(\int_{B_{\rho/2}} |Du_\varepsilon-Dv|^2 \ \mathrm{d}x \biggr)^\frac{1}{2}$$
square the result and take means. This yields 
$$  \fint_{B_{\rho/2}} |Du_\varepsilon-Dv|^2 \ \mathrm{d}x \le \dfrac{c}{\mu^2}  \fint_{B_{\rho/2}} |Du_\varepsilon-\xi|^4 \ \mathrm{d}x +  c \rho^2$$
for a constant $c:=c (p,M,\delta ,C_0,C_1)$. Thanks to assumptions \eqref{bounded} and \eqref{mu}, $|Du_\varepsilon|_{\gamma(x)}$ and $|\xi|_{\gamma(x)}$ are bounded by $\bigl(  1+ \frac{2}{\delta} \bigr) \mu $. Then, an application of Lemma \ref{higherint} and the fact that $\rho \le 1$ yield
\begin{align*}
    \fint_{B_{\rho/2}}  |Du_\varepsilon-\xi|^4 \ \mathrm{d}x \le & \ \dfrac{c}{\mu^{2 \vartheta}}  \fint_{B_{\rho/2}} |Du_\varepsilon-\xi|^{2(1+\vartheta)} \ \mathrm{d}x \\
    \le &  \ \dfrac{c}{\mu^{2 \vartheta}}  \biggl[  \fint_{B_{\rho/2}} |Du_\varepsilon-\xi|^{2} \ \mathrm{d}x    \biggr]^\frac{1}{1+\vartheta}+c\rho^{2(1+\vartheta)}      
    \le  \ c \biggr[  \dfrac{ \mathbf{\Phi} (x_0,\rho)}{\mu^2} \biggl]^\vartheta  \mathbf{\Phi} (x_0,\rho) +c \rho^{2},
\end{align*}
for a constant $c:=c(n,p,M,\delta,C_0,C_1)$. Inserting this inequality above finishes the proof of the lemma.
\endproof
The following a priori estimate for solutions to linear elliptic system can be inferred from \cite{Campanato} once the ellipticity conditions for the quadratic form $\mathcal{B}_\varepsilon \bigl(  x, (Du_\varepsilon)_{x_0,\rho/2} \bigr)$ are established.
\begin{lem}\label{apriori}
    Let the general assumptions of the Section \ref{secND} be in force and assume that \eqref{boundmean} holds true. Then, the weak solutions  $v \in u_\varepsilon + W^{1,2}_0 \bigl( B_{\rho/2}(x_0) , \mathbb{R}^N \bigr)$ of \eqref{linear} satisfies  $v \in  W^{2,2} \bigl( B_{\rho/2}(x_0) , \mathbb{R}^N \bigr)$ and there exists a constant $\tilde{c}_0:=\tilde{c}_0(n,N,p,\delta,C_0,C_1)$ such that for any $\tau \in \big( 0, \frac{1}{2} \bigr]$ we have
    $$  \fint_{B_{\tau \rho}(x_0)} |Dv-  (Dv)_{x_0,\tau \rho}|^2 \ \mathrm{d}x \le \tilde{c}_0 \tau^2 \fint_{B_{\rho/2}(x_0)} |Dv-  (Dv)_{x_0,\rho/2}|^2 \ \mathrm{d}x . $$
\end{lem}
\proof Notice that the ellipticity constant and the upper bound of the quadratic form
\\$\mathcal{B}_\varepsilon \bigl(  x, (Du_\varepsilon)_{x_0,\rho/2} \bigr)$ depend only on $p$ and $\delta$ due to Lemma \ref{ellipticity} and assumption \eqref{boundmean}.
\endproof

\subsection{Exploiting the measure theoretic information}

The aim of this section is to convert the measure theoretic information \eqref{ND} into a lower bound for the mean value of $Du_\varepsilon$ and a smallness of the excess.

\begin{lem}\label{lemexcess}
    Let the general assumptions of Section \ref{secND} be in force. Moreover, assume that \eqref{measureND} holds for some $\nu \in \bigl(0, \frac{1}{8} \bigr]$. Then, there exists a constant $c:=c(n,p,M,\delta,C_0,C_1,C_2)$ such that for any $\tau \in \bigl[ \frac{1}{2}, 1 \bigr)$ there holds 
    $$\mathbf{\Phi}(x_0,\tau \rho) \le c \mu^2 \biggl[ \dfrac{\nu^{\frac{2}{n}}}{(1-\tau)^2}  + \dfrac{\rho^2}{\nu} \biggr].$$
\end{lem}
\proof Throughout the proof we omit the reference to the center $x_0$ and write $B_\rho$ instead of $B_\rho(x_0)$.
We observe that, for every $x \in \Omega$, $\mathcal{G}(x, \cdot ): \mathbb{R}^{nN} \setminus \{ |\zeta|_{\gamma(x)} \le 1  \} \rightarrow \mathbb{R}^{nN} \setminus \{0\} $ is a one to one mapping whose inverse is given by
$$\mathcal{G}^{-1}(x, \tilde{\zeta}):= \dfrac{|\tilde{\zeta}|_{\gamma(x)}+1}{|\tilde{\zeta}|_{\gamma(x)}} \tilde{\zeta}.$$
For every $x \in \Omega$ we define $\xi_x \in \mathbb{R}^{nN}$ by 
$$|\xi_x|^{p-1}_{\gamma(x)}\xi_x:= \bigl(  |\mathcal{G}(y,Du_\varepsilon)|^{p-1}_{\gamma(y)}  \mathcal{G}(y,Du_\varepsilon)\bigr)_{\tau\rho}$$
and let
$$\tilde{\xi}_x:= \mathcal{G}^{-1}(x, \xi_x).$$
Note that $|\xi_x| \le c(\delta+\mu)$ by assumption \eqref{bounded} and $|\tilde{\xi}_x| \le c (1+\delta+\mu)$. Due to the minimality of the integral average $(Du_\varepsilon)_{\tau \rho}$ with respect to the mapping $z \mapsto \int_{B_{\tau \rho}} |Du_\varepsilon-z|^2 \ \mathrm{d}x$, we have
\begin{align*}
    \mathbf{\Phi}(\tau \rho)  \le & \fint_{B_{\tau \rho}} |Du_\varepsilon - \tilde{\xi}_x| \ \mathrm{d}x \\
    =& \ \dfrac{1}{|B_{\tau \rho}|} \int_{E^\nu_{\tau \rho}} |Du_\varepsilon - \tilde{\xi}_x| \ \mathrm{d}x + \dfrac{1}{|B_{\tau \rho}|} \int_{B_{\tau \rho} \setminus E^\nu_{\tau \rho}} |Du_\varepsilon - \tilde{\xi}_x| \ \mathrm{d}x =: I+II.
\end{align*}
We recall that $|Du_\varepsilon| , |\tilde{\xi}_x| \le c(1+\delta+\mu)$ and hence by \eqref{mu} we have $|Du_\varepsilon| , |\tilde{\xi}_x| \le c \bigl(  2+ \frac{1}{\delta} \bigr) \mu$. Due to assumption \eqref{measureND}, we therefore obtain for the integral $II$
$$II \le \dfrac{c(\delta)\mu^2}{|B_{\tau\rho}|} |B_{\tau\rho} \setminus E^\nu_{\tau\rho}| \le \dfrac{c(\delta)\mu^2 \nu}{\tau^n }.$$
For the estimate of $I$ we first note that $|Du_\varepsilon|_{\gamma(x)} \ge 1+ \frac{3}{2} \delta$ on $E^\nu_{\tau\rho}$ since $\nu \in \bigl(0, \frac{1}{8} \bigr]$ and $\mu \ge \delta$. Therefore, an application of Lemma \ref{lem3} yields
\begin{align*}
    I \le  \ \dfrac{c(\delta,C_0,C_1)}{|B_{\tau\rho}|} \int_{E^\nu_{\tau \rho}}  |\mathcal{G}(x,Du_\varepsilon)-\mathcal{G}(x, \tilde{\xi}_x)|^2 \ \mathrm{d} x
    =   \ \dfrac{c(\delta,C_0,C_1)}{|B_{\tau\rho}|} \int_{E^\nu_{\tau \rho}}  |\mathcal{G}(x,Du_\varepsilon)-\xi_x|^2 \ \mathrm{d} x.
\end{align*}
Next, we note that on $E^\nu_{\tau\rho}$ it holds
$$|\mathcal{G}(x,Du_\varepsilon)|_{\gamma(x)} + |\xi_x|_{\gamma(x)} \ge |\mathcal{G}(x,Du_\varepsilon)|_{\gamma(x)} \ge \bigl(  |Du_\varepsilon|_{\gamma(x)} -1 \bigr)_+ > \delta + (1-\nu)\mu >\dfrac{1}{2}\mu. $$
Using this information, Lemma \ref{lem2}, the definition of $\xi_x$ and Poincarè's inequality, we obtain
\begin{align*}
    I \le &    \ \dfrac{c}{\mu^{2p-2}|B_{\tau\rho}|} \int_{E^\nu_{\tau \rho} } 
    \bigl( |\mathcal{G}(x,Du_\varepsilon)|_{\gamma(x)}+|\xi_x|_{\gamma(x)}  \bigr)^{2p-2} |\mathcal{G}(x,Du_\varepsilon)-\xi_x|^2 \ \mathrm{d} x\\
    \le & \  \dfrac{c}{\mu^{2p-2} }\fint_{B_{\tau \rho} } 
    \big| |\mathcal{G}(x,Du_\varepsilon)|^{p-1}_{\gamma(x)} \mathcal{G}(x,Du_\varepsilon) -|\xi_x|^{p-1}_{\gamma(x)} \xi_x  \big|^{2} \ \mathrm{d} x\\
    \le & \ \dfrac{c \rho^2}{\mu^{2p-2}} \biggl[ 
    \fint_{B_{\tau \rho} } 
    \big| D \bigl[|\mathcal{G}(x,Du_\varepsilon)|^{p-1}_{\gamma(x)} \mathcal{G}(x,Du_\varepsilon) \bigr] \big|^\frac{2n}{n+2} \ \mathrm{d} x
    \biggr]^\frac{n+2}{n}\\
    = & \  \dfrac{c \rho^2}{\mu^{2p-2}} \biggl[ 
    \fint_{B_{\tau \rho} } 
    \big| D \bigl[g \bigl( |Du_\varepsilon|_{\gamma} \bigr) D u_\varepsilon\bigr] \big|^\frac{2n}{n+2} \ \mathrm{d} x
    \biggr]^\frac{n+2}{n}
\end{align*}
for a constant $c:=c(p,\delta,C_0,C_1)$. We once again decompose the domain of integration into $E^\nu_{\tau \rho}$ and $B_{\tau\rho} \setminus E^\nu_{\tau \rho} $. Subsequently, applying H\"older's inequality and taking into account assumption \eqref{measureND} leads us to
\begin{align*}
    I \le  \dfrac{c \rho^2}{\mu^{2p-2}|B_{\tau\rho}|} \biggl[ 
    \int_{E^\nu_{\tau \rho} } 
    \big| D \bigl[g \bigl( |Du_\varepsilon|_{\gamma} \bigr) D u_\varepsilon\bigr] \big|^2 \ \mathrm{d} x + \nu^\frac{2}{n} 
    \int_{B_{\tau\rho} \setminus E^\nu_{\tau \rho} } 
    \big| D \bigl[g \bigl( |Du_\varepsilon|_{\gamma} \bigr) D u_\varepsilon\bigr] \big|^2 \ \mathrm{d} x
    \biggr].
\end{align*}
We note that $\tau \ge \frac{1}{2}$ and hence $|B_{\tau \rho}| \ge c(n) \rho^n$. For the first integral we use Lemma \ref{enineq2} while for the second one we use Lemma \ref{enineq1} and the assumption $\mu \ge \delta$. In this way we obtain 
\begin{align*}
    I \le \  \dfrac{c \mu^2 \nu}{(1-\tau)^2} + \dfrac{c \rho^2}{\nu \mu^{2p-2}} + \dfrac{c \nu^\frac{2}{n}\mu^2}{(1-\tau)^2}+ \dfrac{c \rho^2\nu^\frac{2}{\nu}}{\mu^{2p-2}(1-\tau)^2}
    \le  \ c \mu^2 \biggl[  \dfrac{ \nu^\frac{2}{n}}{(1-\tau)^2}+\dfrac{\rho^2}{n}  \biggr]
\end{align*}
for a constant $c:=c(n,p,M,\delta,C_0,C_1,C_2)$. Inserting this estimate above yields the desired inequality.
\endproof

\begin{lem}\label{lem5.5}
    Let the general assumptions of Section \ref{secND} be in force. Then, for any $\sigma \in \bigl( 0, \frac{1}{64C_1} \bigr]$ there exist $\nu \in \bigl( 0, \frac{1}{8} \bigr] $ and $\rho_0 \in (0,1]$, depending at most on $n,p,M,\delta,\sigma,C_0,C_1$ and $C_2$, such that the smallness assumption $\rho \le \rho_0 $ and the measure theoretic hypothesis \eqref{measureND} imply
    \begin{equation}\label{5.4}
        |  (Du_\varepsilon)_{x_0,\rho}  |_{\gamma(x)} \ge 1+ \delta + \dfrac{1}{2} \mu \quad \text{and} \quad \mathbf{\Phi}(x_0, \rho) \le \sigma \mu^2
    \end{equation}
    for any $x \in B_{\rho_0}(x_0)$.
\end{lem}
\proof We let $\tau \in \bigl[ \frac{1}{2}, 1 \bigr)$, $\nu \in \bigl(0, \frac{1}{4} \bigr]$ and $\rho_0 \in (0,1]$. Consider $B_\rho(x_0) \subset B_R$ with $\rho \le \rho_0$. In the following we omit the reference to the center $x_0$. Using the minimality of $(Du_\varepsilon)_\rho$ with respect to the mapping $z \mapsto \int_{B_\rho}|Du_\varepsilon-z|^2 \ \mathrm{d}x$ and decomposing the domain of integration into $B_{\tau\rho}$ and $B_\rho \setminus B_{\tau\rho}$, we obtain
\begin{align*}
    \mathbf{\Phi}(\rho) \le & \fint_{B_\rho} |Du_\varepsilon-(Du_\varepsilon)_{\tau\rho}|^2 \ \mathrm{d}x \\
    =& \ \dfrac{1}{|B_\rho|} \biggl[  \int_{B_{\tau\rho}} |Du_\varepsilon-(Du_\varepsilon)_{\tau\rho}|^2 \ \mathrm{d}x +  \int_{B_\rho \setminus B_{\tau\rho}} |Du_\varepsilon-(Du_\varepsilon)_{\tau\rho}|^2 \ \mathrm{d}x  \biggr]=: I+II.
\end{align*}
For the first integral we use Lemma \ref{lemexcess} and obtain
$$I =  \tau^n \mathbf{\Phi}(\tau\rho) \le c \mu^2 \biggl[  \dfrac{\nu^\frac{2}{n}}{(1-\tau)^2}+ \dfrac{\rho^2}{\nu} \biggr]$$
where $c:=c(n,p,M,\delta,C_0,C_1,C_2)$. For the integral $II$ we use the inequality $|Du_\varepsilon|_\gamma \le 1+\delta+\mu \le c(\delta) \mu$ and get
$$II \le \dfrac{4(1+\delta+\mu)^2|B_\rho \setminus B_{\tau\rho}|}{|B_\rho|} \le c(\delta)\mu^2 (1-\tau)^n \le c(n,\delta) \mu^2 (1-\tau)$$
so that
$$\mathbf{\Phi}(\rho) \le c \mu^2 \biggl[  \dfrac{\nu^\frac{2}{n}}{(1-\tau)^2}+ 1-\tau + \dfrac{\rho^2}{\nu}   \biggr]$$
for a constant $c$ depending at most on $n,p,M,\delta,C_0,C_1$ and $C_2$. Now, we first choose $\tau \in \bigl[  \frac{1}{2},1 \bigr)$ in dependence on $n,p,M,\delta$ and $\sigma$ in such a way that $c (1-\tau) \le \frac{1}{3} \sigma$. Subsequently, we choose $\nu \in \bigl(0, \frac{1}{4} \bigr]$ in dependence on $n,p,M,\delta$ and $\sigma$ such that
$$\nu \le \min \biggl\{  \biggl(  \dfrac{\sigma (1-\tau)^2}{3c} \biggr)^\frac{n}{2} , \dfrac{\delta}{8(1+\delta)} \biggr\}.$$
Finally, we choose $\rho_0 \in (0,1]$ such that
$$\rho_0^2 \le \min \biggl\{  \dfrac{\nu\sigma}{3c}, \biggl(\dfrac{\delta}{8c(1+\delta)}\biggr)^2  \biggr\}.$$
In this way we obtain
\begin{equation}
    \mathbf{\Phi}( \rho) \le \sigma \mu^2. \label{secondineq}
\end{equation}
Now, we prove the first inequality of \eqref{5.4}. We observe that the measure theoretic assumption \eqref{measureND} yields
$$|E^\nu_\rho|> (1-\nu)|B_\rho|.$$
Hence, due to the definition of the set $E^\nu_\rho$, we obtain
\begin{align*}
    \int_{B_\rho} |Du_\varepsilon|_\gamma \ \mathrm{d}x \ge \int_{E^\nu_\rho} |Du_\varepsilon|_\gamma \ \mathrm{d}x \ge (1+\delta+(1-\nu)\mu)|E^\nu_\rho|
    \ge  \ (1-\nu )(1+\delta+(1-\nu)\mu)|B_\rho|.
\end{align*}
On the other hand, due to \eqref{bounded},  \eqref{secondineq} and the Lipschitz continuity of $\gamma_{\alpha\beta}$, for every $y \in B_\rho$ we estimate 
\begin{align*}
    \bigg|  \fint_{B_\rho}  &|Du_\varepsilon|_\gamma \ \mathrm{d}x - |(Du_\varepsilon)_\rho|_{\gamma(y)}  \bigg|\\
    \le & \fint_{B_\rho}  |Du_\varepsilon  - (Du_\varepsilon)_\rho|_{\gamma} \ \mathrm{d}x + \fint_{B_\rho}  \big||(Du_\varepsilon)_\rho |_\gamma -| (Du_\varepsilon)_\rho|_{\gamma(y)} \big| \ \mathrm{d}x \\
    \le & \ \sqrt{C_1\mathbf{\Phi}(\rho)} + c \rho(1+\delta+\mu)
    \le  \ \mu \sqrt{C_1\sigma} +c \rho_0(1+\delta+\mu)
\end{align*}
for any $\rho \le \rho_0$ with $c:=c(n,p,M,\delta,C_0,C_1,C_2)$, which implies
$$  |(Du_\varepsilon)_\rho|_{\gamma(y)} \ge (1-\nu )(1+\delta+(1-\nu)\mu) - \mu \sqrt{C_1\sigma} - c \rho_0(1+\delta+\mu)$$
for every $y \in B_\rho$. Due to the choice of $\nu$ and $\rho_0$, we have $\nu \le \frac{\delta}{8(1+\delta)}$ and  $\rho_0 \le \frac{\delta}{8c(1+\delta)} $. Together with the assumptions $\mu \ge \delta$ and $\sigma \le \frac{1}{64C_1}$, we obtain
\begin{align*}
    (1-\nu) & (1+\delta+(1-\nu)\mu)-\mu \sqrt{C_1\sigma} -c \rho_0(1+\delta+\mu)-1-\delta-\dfrac{1}{2} \mu  \\
    =& -\nu (1+\delta) +(1-\nu)^2 \mu - \mu \sqrt{C_1\sigma}  -c \rho_0(1+\delta+\mu)-\dfrac{1}{2}\mu\\
    \ge & -\dfrac{\delta}{4} + \biggl(  \dfrac{1}{2}-2 \nu- \sqrt{C_1\sigma} -c \rho_0 \biggr) \mu  \\
    \ge & \biggl(  \dfrac{1}{4}-2 \nu- \sqrt{C_1\sigma} -c \rho_0 \biggr) \mu 
    \ge  \biggl(  \dfrac{1}{8}-2 \nu -c \rho_0 \biggr) \mu \ge 0
\end{align*}
after diminishing the values of $\nu$ and $\rho_0$ if necessary. Inserting this above yields the first inequality of \eqref{5.4} and finishes the proof of the lemma.
\endproof

\subsection{Proof of Proposition \ref{ND}}\label{proofND}

The aim of this section is to prove Proposition \ref{ND}. We start with
a decay estimate for the excess $\mathbf{\Phi}(x_0,\rho) $ of $Du_\varepsilon$.

\begin{lem}\label{lem5.6}
    Assume that the general assumptions of the Section \ref{secND} are in force. Let $\tau \in \bigl( 0, \frac{1}{2} \bigr]$ and $\vartheta \in \bigl(0, \frac{1}{2} \bigr]$ be the exponent from Lemma \ref{comparison}. If
    \begin{equation}\label{5.5}
        |(Du_\varepsilon)_{x_0,\rho}  |_{\gamma(x)} \ge 1+ \delta + \dfrac{1}{4} \mu \quad \text{and} \quad \mathbf{\Phi}(x_0, \rho) \le  \tau^\frac{n+2}{\vartheta} \mu^2
    \end{equation}
   hold true for every $x \in B_{\rho}(x_0)$, then we have 
   $$\mathbf{\Phi}(x_0,\tau\rho)\le c_* \bigl[  \tau^2 \mathbf{\Phi}(x_0,\rho) +\tau^{-n} \rho^2\mu^2 \bigr]$$
   with a constant $c_*:=c_*(n,p,M,\delta,C_0,C_1)$.
\end{lem}
\proof Throughout the proof we omit the reference to the center $x_0$ and write $B_\rho$ instead of $B_\rho(x_0)$. By $v \in u_\varepsilon +  W^{1,2}_0 \bigl(  B_\rho, \mathbb{R}^N\bigr)$ we denote the unique weak solution of the linear elliptic system \eqref{linear}. For $\tau \in \bigl( 0, \frac{1}{2} \bigr]$ we have
\begin{align*}
    \mathbf{\Phi}(\rho) \le & \fint_{B_{\tau \rho}} |D u_\varepsilon-(Dv)_{\tau \rho}|^2 \ \mathrm{d}x
    \le  \ 2 \fint_{B_{\tau \rho}} |D u_\varepsilon-Dv|^2 \ \mathrm{d}x 
    + 2 \fint_{B_{\tau \rho}} |D v-(Dv)_{\tau \rho}|^2 \ \mathrm{d}x.
\end{align*}
In view of Lemma \ref{apriori} we deduce
\begin{align*}
     \fint_{B_{\tau \rho}} |D v    -(Dv)_{\tau \rho}|^2 \ \mathrm{d}x 
    \le &  \ c \tau^2 \fint_{B_{\rho/2}} |D v-(Dv)_{ \rho/2}|^2 \ \mathrm{d}x\\
    \le & \ c \tau^2 \fint_{B_{ \rho/2}} |D v-(Du_\varepsilon)_{\tau \rho}|^2 \ \mathrm{d}x \\
    \le & \ c \tau^2 \fint_{B_{ \rho/2}} |D u_\varepsilon-D v|^2 \ \mathrm{d}x + c \tau^2 \fint_{B_{ \rho/2}} |D u_\varepsilon-(Du_\varepsilon)_{ \rho/2}|^2 \ \mathrm{d}x \\
    \le & \ c \tau^2  \fint_{B_{ \rho/2}} |D u_\varepsilon-D v|^2 \ \mathrm{d}x + c \tau^2 \mathbf{\Phi}(\rho)
\end{align*}
where $c:=c(n,N,p,\delta,C_0,C_1)$. Inserting this above and applying Lemma \ref{comparison} and the second inequality of \eqref{5.5}, we end up with
\begin{align*}
    \mathbf{\Phi}(\tau \rho) \le & \ \dfrac{c}{\tau^n} \fint_{B_{ \rho/2}} |D u_\varepsilon-D v|^2 \ \mathrm{d}x + c \tau^2 \mathbf{\Phi}(\rho) \\
    \le & \ c \biggl[\dfrac{1}{\tau^n} \biggl( \dfrac{\mathbf{\Phi}(\rho)}{\mu^2}  \biggr)^\vartheta +\tau^2 \biggr] \mathbf{\Phi}(\rho)+ c \tau^{-n} \rho^2
    \le  \ c_* \bigl(  \tau^2 \mathbf{\Phi}(\rho)+\tau^{-n} \rho^2 \mu^2 \bigr)
\end{align*}
with a constant $c_*:=c_*(n,N,p,\delta,M,C_0,C_1)$.
\endproof

\vspace{0.8cm}
\begin{proof}[{\it Proof of Proposition \ref{ND}.}] By $\vartheta \in \bigl( 0, \frac{1}{2} \bigr] $ we denote the constant from Lemma \ref{comparison} and by $c_*$ the one from Lemma \ref{lem5.6}. For $\beta \in (0,1)$ we choose $\tau \in \bigl(0, \frac{1}{8} \bigr]$ such that
$$\tau  \le \min \biggl\{  \dfrac{1}{8}, 2^{-\frac{1}{\beta}}, (10c_*)^{- \frac{1}{2(1-\beta)} }\biggr\}$$
and let 
$$\sigma := \min \biggl\{  \tau^{\frac{n+2}{\theta}}, \dfrac{\tau^{n+2}}{C_1} \biggr\}.$$
For the particular choice of $\sigma$, we let $\rho_0 \in (0,1]$ be the radius from Lemma \ref{lem5.5}. Moreover, we define 
$$\hat{\rho}:= \min \biggl\{  \rho_0, \biggl( \dfrac{\sigma \tau^n}{2c_*} \biggr)^\frac{1}{2\beta} \tau  \biggr\}.$$
In the following we consider a ball $B_{2\rho}(x_0) \subset B_{r_1}$ with $\rho \le \hat{\rho}$. We omit the reference to the center $x_0$ and write $B_\rho$ instead of $B_\rho(x_0)$. By $\nu \in \bigl(0, \frac{1}{8} \bigr]$ we denote the constant from Lemma \ref{lem5.5} and assume that \eqref{measureND} is satisfied for this particular choice of $\nu$.
From Lemma \ref{lem5.5} applied with $\sigma$ we infer that
\begin{equation}\label{formula4.9}
        |  (Du_\varepsilon)_{\rho}  |_{\gamma(x)} \ge 1+ \delta + \dfrac{1}{4} \mu \quad \text{and} \quad \mathbf{\Phi}( \rho) \le \sigma \mu^2
    \end{equation}
    for any $x \in B_{\rho_0}$. By induction we shall prove that for any $i \in \mathbb{N}$ we have
\begin{equation*}
    \mathbf{\Phi}( \tau^i \rho) \le \sigma \tau^{2\beta i}\mu^2 \eqno{\rm{{ (I)_i}}}
\end{equation*}
 and 
 \begin{equation*}
     |  (Du_\varepsilon)_{\tau^i \rho}  |_{\gamma(x)} \ge 1+ \delta  + \biggl[   \dfrac{1}{2}- \dfrac{1}{8} \sum^{i-1}_{j=0} 2^{-j}  \biggr]\mu , \quad \forall x \in B_{\tau^i \rho}.  \eqno{\rm{{ (II)_i}}}
 \end{equation*}
For $i=1$ we can apply Lemma \ref{lem5.6}, since \eqref{formula4.9} ensures that the assumptions of the lemma are satisfied. Then, from Lemma \ref{lem5.6}, the second inequality of \eqref{formula4.9} and the choice of $\tau, \sigma$ and $\hat{\rho}$, we have
\begin{align*}
    \mathbf{\Phi}(\tau \rho) \le & \ c_* \bigl[  \tau^2 \mathbf{\Phi}( \rho) + \tau^{-n} \rho^2 \mu^2   \bigr]\\
    \le & \ c_* \bigl[  \tau^2 \mathbf{\Phi}( \rho) + \tau^{-n} \rho^{2\beta} \mu^2   \bigr]
    \le \ \dfrac{1}{2} \tau^{2 \beta}  \mathbf{\Phi}(\rho) + \dfrac{c_*}{\tau^n} \rho^{2 \beta} \mu^2
    \le  \ \sigma \tau^{2 \beta} \biggl[ \dfrac{1}{2} + \dfrac{c_* \hat{\rho}^{2 \beta}}{\sigma \tau^{n+2\beta}}  \biggr] \mu^2
    \le  \ \sigma \tau^{2 \beta} \mu^2.
\end{align*}
This proves $(I)_i$. For the proof of $(II)_i$, we use the second inequality of \eqref{formula4.9} and the choice of $\sigma$. For every $x \in B_{\tau^i \rho}$, we derive
\begin{align*}
    |(D u_\varepsilon)_{\tau \rho}-(D u_\varepsilon)_{\rho}|^2_{\gamma(x)}
    \le & \ C_1 \fint_{B_{\tau \rho}} |Du_\varepsilon-(Du_\varepsilon)_\rho|^2 \ \mathrm{d} y 
    \le  \ C_1 \tau^{-n} \mathbf{\Phi}(\rho) \le \tau^2 \mu^2
\end{align*}
so that
$$ |(D u_\varepsilon)_{\tau \rho}-(D u_\varepsilon)_{\rho}|_{\gamma(x)} \le \tau \mu \le \dfrac{1}{8} \mu. $$
Together with the first inequality of \eqref{formula4.9} this implies $(II)_i$.

Now, we consider $i >1$ and prove $(I_i)$ and $(II)_i$ assuming that $(I)_{i-1}$ and $(II)_{i-1}$ hold. From $(I)_{i-1}$ and $(II)_{i-1}$ we observe that the assumptions of Lemma \ref{lem5.6} as formulated in \eqref{5.5} are satisfied on $B_{\tau^i \rho}$. Therefore, applying the lemma with $\tau^i \rho$ instead of $\rho$, recalling the choices of $\tau$ and $\hat{\rho}$ and joining the result with $(I)_{i-1}$ yield
\begin{align*}
    \mathbf{\Phi}(\tau^i \rho) \le & \ c_* \bigl[  \tau^2 \mathbf{\Phi}( \tau^{i-1}\rho) + \tau^{-n} (\tau^{i-1}\rho)^{2 \beta} \mu^2   \bigr]\\
    \le & \ \dfrac{1}{2} \tau^{2 \beta}  \mathbf{\Phi}(\tau^{i-1}\rho) + \dfrac{c_*}{\tau^n} (\tau^{i-1}\rho)^{2 \beta} \mu^2
    \le  \ \sigma \tau^{2 \beta i} \biggl[ \dfrac{1}{2} + \dfrac{c_* \hat{\rho}^{2 \beta}}{\sigma \tau^{n+2\beta}}  \biggr] \mu^2
    \le  \ \sigma \tau^{2 \beta i} \mu^2.
\end{align*}
This proves $(I)_i$. Moreover, from $(I)_{i-1}$ and the choice of $\sigma$
\begin{align*}
    |(D u_\varepsilon)_{\tau^i \rho}-(D u_\varepsilon)_{\tau^{i-1}\rho}|^2_{\gamma(x)}
    \le & \ C_1 \fint_{B_{\tau^i \rho}} |Du_\varepsilon-(Du_\varepsilon)_{\tau^{i-1}\rho}|^2 \ \mathrm{d} y \\
    \le &  \ C_1 \tau^{-n} \mathbf{\Phi}(\tau^{i-1}\rho) 
    \le   \ C_1 \tau^{-n} \sigma  \tau^{2 \beta (i-1)} \mu^2 \le \ \tau^2 \tau^{2 \beta (i-1)} \mu^2
\end{align*}
so that
$$ |(D u_\varepsilon)_{\tau^i \rho}-(D u_\varepsilon)_{\tau^{i-1}\rho}|_{\gamma(x)} \le \ \tau \tau^{ \beta (i-1)} \mu \le \ \dfrac{1}{8} 
     2^{-(i-1)}   \mu$$
by our choice of $\tau$. Together with $(II)_{i-1}$ this proves $(II)_i$.

Now, we come to the proof of \eqref{limit} and \eqref{decay}. For $i \in \mathbb{N}$ we obtain from the minimizing property of the mean value, Lemma \ref{lem3}, $(I)_i$, \eqref{formula4.9} and our choices of $\tau$ and $\sigma$ 
\begin{align}\label{formula4.10}
    \mathbf{\Psi}_{2 \delta}(\tau^i \rho) := & \fint_{B_{\tau^i \rho}} 
    |\mathcal{G}_{2 \delta}(x,Du_\varepsilon)- \bigl( \mathcal{G}_{2 \delta}(x,Du_\varepsilon) \bigr)_{\tau^i \rho}|^2 \ \mathrm{d}x \notag\\
    \le & \  \fint_{B_{\tau^i \rho}} 
    |\mathcal{G}_{2 \delta}(x,Du_\varepsilon)- \mathcal{G}_{2 \delta}\bigl(x,(Du_\varepsilon)_{\tau^i \rho} \bigr)|^2 \ \mathrm{d}x \notag\\
    \le & \ c \mathbf{\Phi}(\tau^i \rho) \le  c \sigma \tau^{2 \beta i} \mu^2 \le   c  \tau^\frac{n+2}{\vartheta} \tau^{2 \beta i} \mu^2 \le  c \tau^{n+2} \tau^{2 \beta i} \mu^2
\end{align}
 for a constant $c:=c(C_0,C_1)$, where in the last inequality we used that $\tau^\frac{n+2}{\vartheta} \le \tau^{n+2}$.
This allow us to compute
\begin{align*}
| \bigl( \mathcal{G}_{2 \delta} & (x,Du_\varepsilon) \bigr)_{\tau^i \rho}
- \bigl( \mathcal{G}_{2 \delta}(x,Du_\varepsilon) \bigr)_{\tau^{i-1} \rho}|^2 \\
     \le & \ \fint_{B_{\tau^i \rho}} 
    |\mathcal{G}_{2 \delta}(x,Du_\varepsilon)- \bigl( \mathcal{G}_{2 \delta}(x,Du_\varepsilon) \bigr)_{\tau^{i-1} \rho}|^2 \ \mathrm{d}x \\
    = & \  \tau^{-n}  \mathbf{\Psi}_{2 \delta}(\tau^{i-1} \rho)  \le 
    c \tau^{n+2} \tau^{2 \beta (i-1)} \mu^2
\end{align*}
for a constant $c:=c(C_0,C_1)$ independent of $i$. Given $j < k$, we use the preceding inequality to conclude that 
\begin{align}\label{formula5.8}
    | \bigl( \mathcal{G}_{2 \delta} & (x,Du_\varepsilon) \bigr)_{\tau^j \rho}
- \bigl( \mathcal{G}_{2 \delta}(x,Du_\varepsilon) \bigr)_{\tau^{k} \rho}| \notag\\
\le &  \sum_{i=j+1}^k | \bigl( \mathcal{G}_{2 \delta}  (x,Du_\varepsilon) \bigr)_{\tau^i \rho}
- \bigl( \mathcal{G}_{2 \delta}(x,Du_\varepsilon) \bigr)_{\tau^{i-1} \rho}|^2   \notag\\
\le & \ c \tau^{\frac{n}{2}+1} \sum_{i=j+1}^k \tau^{\beta(i-1)} \mu 
\le  \ c \tau^{\frac{n}{2}+1} \dfrac{\tau^{\beta j}}{1- \tau^{\beta}} \mu \le \ 2 c \tau^{\frac{n}{2}+1} \tau^{\beta j} \mu,
\end{align}
where the constant $c$ is independent of $j$ and $k$. This shows that $\bigl( ( \mathcal{G}_{2 \delta}  (x,Du_\varepsilon) )_{\tau^i \rho}  \bigr)_{i \in \mathbb{N}}$ is a Cauchy sequence and therefore the limit
$$\Gamma_{x_0} := \lim_{i \rightarrow + \infty} \bigl( \mathcal{G}_{2 \delta}  (x,Du_\varepsilon) \bigr)_{\tau^i \rho}$$
exists. Due to the assumptions \eqref{bounded} and the property $(iii)$ of $\gamma$, we deduce that
\begin{align*}
    \big| \bigl(\mathcal{G}_{2 \delta}  (x,Du_\varepsilon) \bigr)_{\tau^i \rho}  \big|_{\gamma(x_0)} \le & \fint_{B_{\tau^i \rho}} 
    \big| \bigl(\mathcal{G}_{2 \delta}  (x,Du_\varepsilon) \bigr)_{\tau^i \rho}  \big|_{\gamma(x_0)} \ \mathrm{d}x \\
    \le &  \fint_{B_{\tau^i \rho}} 
    \bigl(\big| \bigl(\mathcal{G}_{2 \delta}  (x,Du_\varepsilon) \bigr)_{\tau^i \rho}  \big|_{\gamma(x_0)} - \big| \bigl(\mathcal{G}_{2 \delta}  (x,Du_\varepsilon) \bigr)_{\tau^i \rho}  \big|_{\gamma(x)} \bigr) \ \mathrm{d}x \\
    & +  \fint_{B_{\tau^i \rho}} 
    \big| \bigl(\mathcal{G}_{2 \delta}  (x,Du_\varepsilon) \bigr)_{\tau^i \rho}  \big|_{\gamma(x_0)} \ \mathrm{d}x \\
    \le & \  (1+c\tau^i \rho) \mu,
\end{align*}
for any $i \in \mathbb{N}$ and for a constant $c:=c(C_0,C_1,C_2)$. By letting $i \rightarrow + \infty$, this implies
$$ |\Gamma_{x_0}|_{\gamma(x_0)} \le \mu $$
which yields the assertion \eqref{gammanormx0}.
Moreover, passing to the limit as $k \rightarrow + \infty$ in \eqref{formula5.8} yields
$$ | \bigl( \mathcal{G}_{2 \delta}  (x,Du_\varepsilon) \bigr)_{\tau^i \rho}
- \Gamma_{x_0} | \le 2 c \tau^{\frac{n}{2}+1} \tau^{\beta j} \mu, \quad \text{for any} \ j \in \mathbb{ N}. $$
Joining this with \eqref{formula4.10}, we find
\begin{align*}
    \fint_{B_{\tau^j \rho}}   |  \mathcal{G}_{2 \delta}  (x,Du_\varepsilon) 
- \Gamma_{x_0} |^2 \ \mathrm{d}x \le & \ 2 \mathbf{\Psi}_{2 \delta}(\tau^j \rho) + 2 | \bigl( \mathcal{G}_{2 \delta}  (x,Du_\varepsilon) \bigr)_{\tau^i \rho}
- \Gamma_{x_0} |^2 
\le  \ c \tau^{n+2} \tau^{2 \beta j} \mu^2.
\end{align*}
For $r \in (0, \rho]$ there exists $j \in \mathbb{N}_0$ such that $\tau^{j+1}  \rho < r \le \tau^j \rho $. Then, we obtain from the last inequality
\begin{align*}
    \fint_{B_r}  | \mathcal{G}_{2 \delta}  (x,Du_\varepsilon) 
- \Gamma_{x_0} |^2 \ \mathrm{d}x \le &  \ \tau^{-n} \fint_{B_{\tau^j \rho}}   |  \mathcal{G}_{2 \delta}  (x,Du_\varepsilon) 
- \Gamma_{x_0} |^2 \ \mathrm{d}x\\ 
\le &  \  c \tau^{2} \tau^{2 \beta j} \mu^2 \le \ c \tau^{2(1 -\beta)} \biggl(  \dfrac{r}{\rho}\biggr)^{2\beta} \mu^2.
\end{align*}
Recalling that the constant $c$ depends only on $C_0$ and $C_1$, $\tau$ can always be diminished in such a way $c \tau^{2(1 -\beta)} < 1$. Therefore, we get 
$$\fint_{B_r}  | \mathcal{G}_{2 \delta}  (x,Du_\varepsilon) 
- \Gamma_{x_0} |^2 \ \mathrm{d}x \le \biggl(  \dfrac{r}{\rho}\biggr)^{2\beta} \mu^2 .$$
This implies
$$\big| \bigl(\mathcal{G}_{2 \delta}  (x,Du_\varepsilon) \bigr)_r 
- \Gamma_{x_0} |^2 \le \fint_{B_r}  \big| \mathcal{G}_{2 \delta}  (x,Du_\varepsilon) 
- \Gamma_{x_0} |^2 \ \mathrm{d}x \le \biggl(  \dfrac{r}{\rho}\biggr)^{2\beta} \mu^2 $$
so that also
$$\Gamma_{x_0} := \lim_{r \rightarrow 0^+}  \bigl(\mathcal{G}_{2 \delta}  (x,Du_\varepsilon) \bigr)_{r}.$$
This completes the proof of Proposition \ref{ND}.
\end{proof}

\section{The degenerate regime}\label{secD}
\noindent In this section we give the proof of Proposition \ref{D}. 

For given $\varepsilon \in (0,1]$ we denote by $u_\varepsilon \in W^{1,p} \bigl( B_R, \mathbb{R}^N \bigr)$ the unique weak solution of the regularized system \eqref{regsystem}. We assume that \eqref{bounded} is in force for some $\mu, \delta >0$ on some ball $B_{2\rho}(x_0) \subset B_{r_1} \Subset B_R$. Let $U_\varepsilon := \bigl( |Du_\varepsilon| -1 -\delta \bigr)^2_+$ denote the function defined in \eqref{Ue}. Note that \eqref{bounded} implies
 $$\sup_{B_{2\rho}(x_0)} U_\varepsilon \le \mu^2.$$



Next, we prove a De Giorgi type lemma for the function $U_\varepsilon$, which can be found for instance in \cite{BoDuGiPdN,DiBenedetto}.

\begin{lem}\label{reducingthesupremum}
    Assume that the general assumptions of the Section \ref{secD} are in force and let $\theta \in (0,1)$. Then, there exists $\tilde{\nu}:= \tilde{\nu}(n,p,M,\delta,C_0,C_1,C_2)$ such that the measure theoretic assumption 
    \begin{equation}\label{MTA}
        |\{   x \in B_\rho(x_0): U_\varepsilon(x) > (1-\theta) \mu^2 \}| < \tilde{\nu} |B_\rho(x_0)|
    \end{equation}
    implies that either $\mu^2 < \rho/\theta$ or
    $$U_\varepsilon \le \biggl( 1- \dfrac{1}{2} \theta \biggr) \mu^2$$
    hold true.
\end{lem}

The following lemma can be found in \cite{BoDuGiPdN,DiBenedetto}.
\begin{lem}\label{lemma6.4}
     Assume that the general assumptions of the Section \ref{secD} are in force and assume that \eqref{measureD} is satisfied for some $\nu \in (0,1)$. Then, for any $i_* \in  \mathbb{N}$ we either have
     $$\mu^2 < \dfrac{2^{i_*}\rho}{\nu}$$
     or
     $$ \big| \bigl\{  x \in B_\rho(x_0) : U_\varepsilon(x) > \bigl( 1-2^{-i_*} \nu  \bigr) \mu^2 \bigr\}   \big| < \dfrac{c_*}{\nu \sqrt{i_*}} |B_\rho(x_0)|$$
     for a constant $c_*:=c_* (n,p,M,\delta,C_0,C_1,C_2)$.
\end{lem}

Now, we are able to prove Proposition \ref{D}.
\begin{proof}[{\it Proof of Proposition \ref{D}.}] 
Let $\tilde{\nu} \in (0,1)$ and $c_*$ be the constants from Lemmas \ref{reducingthesupremum} and \ref{lemma6.4}. Note that both depend on $n,p,M,\delta,C_0,C_1$ and $C_2$. Choose $i_* \in \mathbb{N}$ such that 
$$i_* \ge \biggl(\dfrac{c_*}{\nu \tilde{\nu}} \biggr)^2.$$
Then, $i_*$ depends on $n,p,M,\delta,C_0,C_1,C_2$ and $\nu$. Lemma \ref{lemma6.4} implies either $\mu^2 < \frac{2^{i_*}\rho}{\nu}$ or 
$$ \big| \bigl\{  x \in B_\rho(x_0) : U_\varepsilon(x) > \bigl( 1-2^{-i_*} \nu  \bigr) \mu^2 \bigr\}   \big| < \dfrac{c_*}{\nu \sqrt{i_*}} |B_\rho(x_0)| \le \tilde{\nu} |B_\rho(x_0)|.$$
In the first case, the proposition is proved with $\hat{c}:= 2^{i_*}/ \nu$, while in the second case we may apply Lemma \ref{reducingthesupremum} with $\theta :=2^{i_*} \nu$. therefore, either $\mu^2 < 2^{i_*} \rho/ \nu$ or
$$U_\varepsilon \le (1-2^{-(i_*+1)}\nu) \mu^2 \quad \text{in} \ B_{\rho/2}(x_0).$$
The first alternative coincides with the first alternative above, while the second one implies 
$$\sup_{B_{\rho/2}(x_0)} \big| \mathcal{G}_\delta (x,Du_\varepsilon) \big|_{\gamma(x)} \le k \mu $$
for any $k \ge \sqrt{1-2^{-(i_*+1)}\nu}$, since $U_\varepsilon(x)= \big| \mathcal{G}_\delta (x,Du_\varepsilon) \big|^2_{\gamma(x)}$. Therefore, we may choose $k \in [2^{-\frac{1}{2}},1)$.
\end{proof}

\vskip0.5cm
\noindent \textbf{Acknowledgements.} 
The author wishes to thank Prof.\ V.\ B\"ogelein for fruitful discussions. The author acknowledges the hospitality of the University of Salzburg, where most of this work has been done. The author has been partially supported by the Gruppo Nazionale per l’Analisi Matematica, la Probabilità e le loro Applicazioni (GNAMPA) of the Istituto Nazionale di Alta Matematica (INdAM) through INdAM-GNAMPA
project (CUP\_E53C22001930001).

\begin{footnotesize}

\end{footnotesize}


\begin{thebibliography}{}

\bibitem{BoDuGiPdN}
\newblock B\"ogelein, V., Duzaar, F., Giova, R., Passarelli di Napoli, A., 
\newblock \textit{Higher regularity in congested traffic dynamics}, 
\newblock Mathematische Annalen.
https://doi.org/10.1007/s00208-022-02375-y

\bibitem{BoDuLiSc}
\newblock B\"ogelein, V., Duzaar, F., Liao, N., Scheven, C.,
\newblock \textit{Gradient Hölder regularity for degenerate parabolic systems},
\newblock Nonlinear Analysis 225 (2022) 113119

\bibitem{brasco1}
\newblock Brasco, L.,
\newblock \textit{Global $L^\infty$-gradient estimates for solutions to a certain degenerate elliptic equation}, \newblock Nonlinear
Anal. 74(2), 516--531 (2011)

\bibitem{brasco2}
\newblock Brasco, L., Carlier, G.,
\newblock \textit{On certain anisotropic elliptic equations arising in congested optimal transport:
local gradient bounds},
\newblock Adv. Calc. Var. 7, 379--407 (2014)


\bibitem{brasco}
\newblock Brasco, L., Carlier, G., Santambrogio, F.,
\newblock \textit{Congested traffic dynamics, weak flows and very degenerate
elliptic equations},
\newblock J. Math. Pures Appl. (9) 93(6), 652--671 (2010)



\bibitem{Campanato}
\newblock Campanato, S.,
\newblock \textit{Equazioni ellittiche del II ordine e spazi $\mathcal{L}^{2,\lambda}$},
\newblock Ann. Mat. Pura Appl. 69, 321--381 (1965)

\bibitem{carlier}
\newblock Carlier, G., Jimenez, C., Santambrogio, F.,
\newblock \textit{Optimal transportation with traffic congestion and Wardrop
equilibria},
\newblock SIAM J. Control Optim. 47(3), 1330--1350 (2008)


\bibitem{carlier1}
\newblock Carlier, G., Santambrogio, F., \newblock \textit{A continuous theory of traffic congestion and Wardrop equilibria},
\newblock Jo. Math. Sci. 181(6), 792--804 (2012)


\bibitem{chipot}
\newblock Chipot, M., Evans, L.C.,
\newblock \textit{Linearization at infinity and Lipschitz estimates for certain problems in
calculus of variations},
\newblock Proc. R. Soc. Edinburgh Sect. A 102, 291--303 (1986)


\bibitem{clop}
\newblock Clop, A., Giova, R., Hatami, F., Passarelli di Napoli, A.,
\newblock \textit{Very degenerate elliptic equations under almost critical Sobolev regularity}, 
\newblock Forum Math. 32 (2020) 1515--1537


\bibitem{colombo}
\newblock Colombo, M., Figalli, A.,
\newblock \textit{An excess-decay result for a class of degenerate elliptic equations},
\newblock Discrete Contin. Dyn. Syst., Ser. S 7 (2014) 631--652.



\bibitem{DiBenedetto}
\newblock Di Benedetto, E.,
\newblock \textit{Partial Differential Equations}, \newblock Cornerstones, 2nd edn. Birkh\"auser Boston Ltd, Boston
(2010)

\bibitem{evans}
\newblock Evans, L.C.,
\newblock \textit{Partial Differential Equations},
\newblock American Mathematical Society, Providence, RI (2010)


\bibitem{fonseca}
\newblock Fonseca, I., Fusco, N., Marcellini, P.,
\newblock \textit{An existence result for a nonconvex variational problem via regularity},
\newblock ESAIM Control Optim. Calc. Var. 7 (2002) 69--95



\bibitem{Giusti}
\newblock Giusti, E.,
\newblock \textit{Direct Methods in the Calculus of Variations},
\newblock World Scientific Publishing Company, Tuck Link (2003)


\bibitem{mons}
\newblock Mons, L.,
\newblock \textit{Higher regularity for minimizers of very degenerate integral functionals},
\newblock J. Math. Anal. Appl.518(2023)126717



\bibitem{santambrogio}
\newblock Santambrogio, F., Vespri, V.,
\newblock \textit{Continuity in two dimensions for a very degenerate elliptic equation},
\newblock Nonlinear Anal. 73, 3832--3841 (2010)

\bibitem{tolk}
\newblock Tolksdorf, P.,
\newblock \textit{Everywhere-regularity for some quasilinear systems with a lack of ellipticity},
\newblock Ann. Mat.
Pura Appl. 134(4), 241--266 (1983)


\bibitem{uh}
\newblock Uhlenbeck, K.,
\newblock \textit{Regularity for a class of nonlinear elliptic systems},
\newblock Acta Math. 138, 219--240 (1977)


\bibitem{ur}
\newblock Uraltseva, N.N.,
\newblock \textit{Degenerate quasilinear elliptic systems},
\newblock Zap. Nau\^cn Sem. Leningrad. Otdel. Mat.
Inst. Steklov. (LOMI) 7, 184--222 (1968)

\end{thebibliography}
\end{document}